
\catcode`!=11 
 
  

\def\PiC{P\kern-.12em\lower.5ex\hbox{I}\kern-.075emC}
\def\PiCTeX{\PiC\kern-.11em\TeX}

\def\!ifnextchar#1#2#3{%
  \let\!testchar=#1%
  \def\!first{#2}%
  \def\!second{#3}%
  \futurelet\!nextchar\!testnext}
\def\!testnext{%
  \ifx \!nextchar \!spacetoken 
    \let\!next=\!skipspacetestagain
  \else
    \ifx \!nextchar \!testchar
      \let\!next=\!first
    \else 
      \let\!next=\!second 
    \fi 
  \fi
  \!next}
\def\\{\!skipspacetestagain} 
  \expandafter\def\\ {\futurelet\!nextchar\!testnext} 
\def\\{\let\!spacetoken= } \\  

\def\!tfor#1:=#2\do#3{%
  \edef\!fortemp{#2}%
  \ifx\!fortemp\!empty 
    \else
    \!tforloop#2\!nil\!nil\!!#1{#3}%
  \fi}
\def\!tforloop#1#2\!!#3#4{%
  \def#3{#1}%
  \ifx #3\!nnil
    \let\!nextwhile=\!fornoop
  \else
    #4\relax
    \let\!nextwhile=\!tforloop
  \fi 
  \!nextwhile#2\!!#3{#4}}

\def\!etfor#1:=#2\do#3{%
  \def\!!tfor{\!tfor#1:=}%
  \edef\!!!tfor{#2}%
  \expandafter\!!tfor\!!!tfor\do{#3}}

\def\!cfor#1:=#2\do#3{%
  \edef\!fortemp{#2}%
  \ifx\!fortemp\!empty 
  \else
    \!cforloop#2,\!nil,\!nil\!!#1{#3}%
  \fi}
\def\!cforloop#1,#2\!!#3#4{%
  \def#3{#1}%
  \ifx #3\!nnil
    \let\!nextwhile=\!fornoop 
  \else
    #4\relax
    \let\!nextwhile=\!cforloop
  \fi
  \!nextwhile#2\!!#3{#4}}

\def\!ecfor#1:=#2\do#3{%
  \def\!!cfor{\!cfor#1:=}%
  \edef\!!!cfor{#2}%
  \expandafter\!!cfor\!!!cfor\do{#3}}

\def\!empty{}
\def\!nnil{\!nil}
\def\!fornoop#1\!!#2#3{}

\def\!ifempty#1#2#3{%
  \edef\!emptyarg{#1}%
  \ifx\!emptyarg\!empty
    #2%
  \else
    #3%
  \fi}
 
\def\!getnext#1\from#2{%
  \expandafter\!gnext#2\!#1#2}%
\def\!gnext\\#1#2\!#3#4{%
  \def#3{#1}%
  \def#4{#2\\{#1}}%
  \ignorespaces}

%
\def\!getnextvalueof#1\from#2{%
  \expandafter\!gnextv#2\!#1#2}%
\def\!gnextv\\#1#2\!#3#4{%
  #3=#1%
  \def#4{#2\\{#1}}%
  \ignorespaces}

\def\!copylist#1\to#2{%
  \expandafter\!!copylist#1\!#2}
\def\!!copylist#1\!#2{%
  \def#2{#1}\ignorespaces}

\def\!wlet#1=#2{%
  \let#1=#2 
  \wlog{\string#1=\string#2}}
 
\def\!listaddon#1#2{%
  \expandafter\!!listaddon#2\!{#1}#2}
\def\!!listaddon#1\!#2#3{%
  \def#3{#1\\#2}}
 

\def\!rightappend#1\withCS#2\to#3{\expandafter\!!rightappend#3\!#2{#1}#3}
\def\!!rightappend#1\!#2#3#4{\def#4{#1#2{#3}}}

\def\!leftappend#1\withCS#2\to#3{\expandafter\!!leftappend#3\!#2{#1}#3}
\def\!!leftappend#1\!#2#3#4{\def#4{#2{#3}#1}}

\def\!lop#1\to#2{\expandafter\!!lop#1\!#1#2}
\def\!!lop\\#1#2\!#3#4{\def#4{#1}\def#3{#2}}



\def\!loop#1\repeat{\def\!body{#1}\!iterate}
\def\!iterate{\!body\let\!next=\!iterate\else\let\!next=\relax\fi\!next}
 
\def\!!loop#1\repeat{\def\!!body{#1}\!!iterate}
\def\!!iterate{\!!body\let\!!next=\!!iterate\else\let\!!next=\relax\fi\!!next}
 
\def\!removept#1#2{\edef#2{\expandafter\!!removePT\the#1}}
{\catcode`p=12 \catcode`t=12 \gdef\!!removePT#1pt{#1}}

\def\placevalueinpts of <#1> in #2 {%
  \!removept{#1}{#2}}
 
\def\!mlap#1{\hbox to 0pt{\hss#1\hss}}
\def\!vmlap#1{\vbox to 0pt{\vss#1\vss}}
 
\def\!not#1{%
  #1\relax
    \!switchfalse
  \else
    \!switchtrue
  \fi
  \if!switch
  \ignorespaces}


 

\let\!!!wlog=\wlog              
\def\wlog#1{}    

\newskip\headingtoplotskip      
\newskip\linethickness         
\newskip\longticklength         
\newskip\plotsymbolspacing     
\newskip\shortticklength        
\newskip\stackleading           
\newskip\tickstovaluesleading   
\newskip\totalarclength        
\newskip\valuestolabelleading   

\newbox\!boxA                   
\newbox\!boxB                   
\newbox\!picbox                 
\newbox\!plotsymbol             
\newbox\!putobject              
\newbox\!shadesymbol            

\newcount\!countA               
\newcount\!countB               
\newcount\!countC               
\newcount\!countD               
\newcount\!countE               
\newcount\!countF               
\newcount\!countG               
\newcount\!fiftypt              
\newcount\!intervalno           
\newcount\!npoints              
\newcount\!nsegments            
\newcount\!ntemp                
\newcount\!parity               
\newcount\!scalefactor          
\newcount\!tickcase             

\newskip\!Xleft                
\newskip\!Xright               
\newskip\!Xsave                
\newskip\!Ybot                 
\newskip\!Ysave                
\newskip\!Ytop                 
\newskip\!angle                
\newskip\!arclength            
\newskip\!areabloc             
\newskip\!arealloc             
\newskip\!arearloc             
\newskip\!areatloc             
\newskip\!bshrinkage           
\newskip\!checkbot             
\newskip\!checkleft            
\newskip\!checkright           
\newskip\!checktop             
\newdimen\!dimenA               
\newdimen\!dimenB               
\newdimen\!dimenC               
\newdimen\!dimenD               
\newdimen\!dimenE               
\newdimen\!dimenF               
\newdimen\!dimenG               
\newdimen\!dimenH               
\newdimen\!dimenI               
\newdimen\!distacross           
\newdimen\!downlength           
\newdimen\!dp                   
\newdimen\!dshade               
\newdimen\!dxpos                
\newdimen\!dxprime              
\newdimen\!dypos                
\newdimen\!dyprime              
\newdimen\!ht                   
\newdimen\!leaderlength         
\newdimen\!lshrinkage           
\newdimen\!midarclength         
\newdimen\!offset               
\newdimen\!plotheadingoffset    
\newdimen\!plotsymbolxshift     
\newdimen\!plotsymbolyshift     
\newdimen\!plotxorigin          
\newdimen\!plotyorigin          
\newdimen\!rshrinkage           
\newdimen\!shadesymbolxshift    
\newdimen\!shadesymbolyshift    
\newdimen\!tshrinkage           
\newdimen\!uplength             
\newdimen\!wd                   
\newdimen\!xB                   
\newdimen\!xC                   
\newdimen\!xE                   
\newdimen\!xM                   
\newdimen\!xS                   
\newdimen\!xaxislength          
\newdimen\!xdiff                
\newdimen\!xleft                
\newdimen\!xloc                 
\newdimen\!xorigin              
\newdimen\!xpivot               
\newdimen\!xpos                 
\newdimen\!xprime               
\newdimen\!xright               
\newdimen\!xshade               
\newdimen\!xshift               
\newdimen\!xtemp                
\newdimen\!xunit                
\newdimen\!xxE                  
\newdimen\!xxM                  
\newdimen\!xxS                  
\newdimen\!xxloc                
\newdimen\!yB                   
\newdimen\!yC                   
\newdimen\!yE                   
\newdimen\!yM                   
\newdimen\!yS                   
\newdimen\!yaxislength          
\newdimen\!ybot                 
\newdimen\!ydiff                
\newdimen\!yloc                 
\newdimen\!yorigin              
\newdimen\!ypivot               
\newdimen\!ypos                 
\newdimen\!yprime               
\newdimen\!yshade               
\newdimen\!yshift               
\newdimen\!ytemp                
\newdimen\!ytop                 
\newdimen\!yunit                
\newdimen\!yyE                  
\newdimen\!yyM                  
\newdimen\!yyS                  
\newdimen\!yyloc                

\newif\if!axisvisible           
\newif\if!gridlinestoo          
\newif\if!keepPO                
\newif\if!placeaxislabel        
\newif\if!switch                
\newif\if!xswitch               

\newtoks\!axisLaBeL             
\newtoks\!keywordtoks           

\newwrite\!replotfile           

\newhelp\!keywordhelp{The keyword mentioned in the error message in unknown. 
Replace NEW KEYWORD in the indicated response by the keyword that 
should have been specified.}    

\!wlet\!!origin=\!xM                   
\!wlet\!!unit=\!uplength               
\!wlet\!Lresiduallength=\!dimenG       
\!wlet\!Rresiduallength=\!dimenF       
\!wlet\!axisLength=\!distacross        
\!wlet\!axisend=\!ydiff                
\!wlet\!axisstart=\!xdiff              
\!wlet\!axisxlevel=\!arclength         
\!wlet\!axisylevel=\!downlength        
\!wlet\!beta=\!dimenE                  
\!wlet\!gamma=\!dimenF                 
\!wlet\!shadexorigin=\!plotxorigin     
\!wlet\!shadeyorigin=\!plotyorigin     
\!wlet\!ticklength=\!xS                
\!wlet\!ticklocation=\!xE              
\!wlet\!ticklocationincr=\!yE          
\!wlet\!tickwidth=\!yS                 
\!wlet\!totalleaderlength=\!dimenE     
\!wlet\!xone=\!xprime                  
\!wlet\!xtwo=\!dxprime                 
\!wlet\!ySsave=\!yM                    
\!wlet\!ybB=\!yB                       
\!wlet\!ybC=\!yC                       
\!wlet\!ybE=\!yE                       
\!wlet\!ybM=\!yM                       
\!wlet\!ybS=\!yS                       
\!wlet\!ybpos=\!yyloc                  
\!wlet\!yone=\!yprime                  
\!wlet\!ytB=\!xB                       
\!wlet\!ytC=\!xC                       
\!wlet\!ytE=\!downlength               
\!wlet\!ytM=\!arclength                
\!wlet\!ytS=\!distacross               
\!wlet\!ytpos=\!xxloc                  
\!wlet\!ytwo=\!dyprime                 

\def\!zpt{0pt}                              
\!xunit=1pt
\!yunit=1pt
\!arearloc=\!xunit
\!areatloc=\!yunit
\!dshade=5pt
\!leaderlength=24in
\def\!tfs{256}                              
\def\!wmax{5.3pt}                           
\def\!wmin{2.7pt}                           
\!xaxislength=\!xunit
\!xpivot=\!zpt
\!yaxislength=\!yunit 
\!ypivot=\!zpt
\plotsymbolspacing=.4pt
  \!dimenA=50pt \def\!fiftypt{\the\!dimenA}     

\def\!rootten{3.162278pt}                   
\def\!tenAa{8.690286pt}                     
\def\!tenAc{2.773839pt}                     
\def\!tenAe{2.543275pt}                     

\def\!cosrotationangle{1}      
\def\!sinrotationangle{0}      
\def\!xpivotcoord{0}           
\def\!xref{0}                  
\def\!xshadesave{0}            
\def\!ypivotcoord{0}           
\def\!yref{0}                  
\def\!yshadesave{0}            
\def\!zero{0}                  

\let\wlog=\!!!wlog
%
  
\def\normalgraphs{%
  \longticklength=.4\baselineskip
  \shortticklength=.25\baselineskip
  \tickstovaluesleading=.25\baselineskip
  \valuestolabelleading=.8\baselineskip
  \linethickness=.4pt
  \stackleading=.17\baselineskip
  \headingtoplotskip=1.5\baselineskip
  \visibleaxes
  \ticksout
  \nogridlines
  \unloggedticks}
%
\def\setplotarea x from #1 to #2, y from #3 to #4 {%
  \!arealloc=\!M{#1}\!xunit \advance \!arealloc -\!xorigin
  \!areabloc=\!M{#3}\!yunit \advance \!areabloc -\!yorigin
  \!arearloc=\!M{#2}\!xunit \advance \!arearloc -\!xorigin
  \!areatloc=\!M{#4}\!yunit \advance \!areatloc -\!yorigin
  \!initinboundscheck
  \!xaxislength=\!arearloc  \advance\!xaxislength -\!arealloc
  \!yaxislength=\!areatloc  \advance\!yaxislength -\!areabloc
  \!plotheadingoffset=\!zpt
  \!dimenput {{\setbox0=\hbox{}\wd0=\!xaxislength\ht0=\!yaxislength\box0}}
     [bl] (\!arealloc,\!areabloc)}
%
\def\visibleaxes{%
  \def\!axisvisibility{\!axisvisibletrue}}

%

\def\!fixkeyword#1{%
  \errhelp=\!keywordhelp
  \errmessage{Unrecognized keyword `#1': \the\!keywordtoks{NEW KEYWORD}'}}

\!keywordtoks={enter `i\fixkeyword}

\def\fixkeyword#1{%
  \!nextkeyword#1 }


\def\axis {%
  \def\!nextkeyword##1 {%
    \expandafter\ifx\csname !axis##1\endcsname \relax
      \def\!next{\!fixkeyword{##1}}%
    \else
      \def\!next{\csname !axis##1\endcsname}%
    \fi
    \!next}%
  \!offset=\!zpt
  \!axisvisibility
  \!placeaxislabelfalse
  \!nextkeyword}

\def\!axisbottom{%
  \!axisylevel=\!areabloc
  \def\!tickxsign{0}%
  \def\!tickysign{-}%
  \def\!axissetup{\!axisxsetup}%
  \def\!axislabeltbrl{t}%
  \!nextkeyword}

\def\!axistop{%
  \!axisylevel=\!areatloc
  \def\!tickxsign{0}%
  \def\!tickysign{+}%
  \def\!axissetup{\!axisxsetup}%
  \def\!axislabeltbrl{b}%
  \!nextkeyword}

\def\!axisleft{%
  \!axisxlevel=\!arealloc
  \def\!tickxsign{-}%
  \def\!tickysign{0}%
  \def\!axissetup{\!axisysetup}%
  \def\!axislabeltbrl{r}%
  \!nextkeyword}

\def\!axisright{%
  \!axisxlevel=\!arearloc
  \def\!tickxsign{+}%
  \def\!tickysign{0}%
  \def\!axissetup{\!axisysetup}%
  \def\!axislabeltbrl{l}%
  \!nextkeyword}

\def\!axisshiftedto#1=#2 {%
  \if 0\!tickxsign
    \!axisylevel=\!M{#2}\!yunit
    \advance\!axisylevel -\!yorigin
  \else
    \!axisxlevel=\!M{#2}\!xunit
    \advance\!axisxlevel -\!xorigin
  \fi
  \!nextkeyword}

\def\!axisvisible{%
  \!axisvisibletrue  
  \!nextkeyword}

\def\!axisinvisible{%
  \!axisvisiblefalse
  \!nextkeyword}

\def\!axislabel#1 {%
  \!axisLaBeL={#1}%
  \!placeaxislabeltrue
  \!nextkeyword}

\expandafter\def\csname !axis/\endcsname{%
  \!axissetup 
  \if!placeaxislabel
    \!placeaxislabel
  \fi
  \if +\!tickysign 
    \!dimenA=\!axisylevel
    \advance\!dimenA \!offset 
    \advance\!dimenA -\!areatloc 
    \ifdim \!dimenA>\!plotheadingoffset
      \!plotheadingoffset=\!dimenA 
    \fi
  \fi}

\def\grid #1 #2 {%
  \!countA=#1\advance\!countA 1
  \axis bottom invisible ticks length <\!zpt> andacross quantity {\!countA} /
  \!countA=#2\advance\!countA 1
  \axis left   invisible ticks length <\!zpt> andacross quantity {\!countA} / }

\def\plotheading#1 {%
  \advance\!plotheadingoffset \headingtoplotskip
  \!dimenput {#1} [B] <.5\!xaxislength,\!plotheadingoffset>
    (\!arealloc,\!areatloc)}

\def\!axisxsetup{%
  \!axisxlevel=\!arealloc
  \!axisstart=\!arealloc
  \!axisend=\!arearloc
  \!axisLength=\!xaxislength
  \!!origin=\!xorigin
  \!!unit=\!xunit
  \!xswitchtrue
  \if!axisvisible 
    \!makeaxis
  \fi}

\def\!axisysetup{%
  \!axisylevel=\!areabloc
  \!axisstart=\!areabloc
  \!axisend=\!areatloc
  \!axisLength=\!yaxislength
  \!!origin=\!yorigin
  \!!unit=\!yunit
  \!xswitchfalse
  \if!axisvisible
    \!makeaxis
  \fi}

\def\!makeaxis{%
  \setbox\!boxA=\hbox{
    \beginpicture
      \!setdimenmode
      \setcoordinatesystem point at {\!zpt} {\!zpt}   
      \putrule from {\!zpt} {\!zpt} to
        {\!tickysign\!tickysign\!axisLength} 
        {\!tickxsign\!tickxsign\!axisLength}
    \endpicturesave <\!Xsave,\!Ysave>}%
    \wd\!boxA=\!zpt
    \!placetick\!axisstart}

\def\!placeaxislabel{%
  \advance\!offset \valuestolabelleading
  \if!xswitch
    \!dimenput {\the\!axisLaBeL} [\!axislabeltbrl]
      <.5\!axisLength,\!tickysign\!offset> (\!axisxlevel,\!axisylevel)
    \advance\!offset \!dp  
    \advance\!offset \!ht  
  \else
    \!dimenput {\the\!axisLaBeL} [\!axislabeltbrl]
      <\!tickxsign\!offset,.5\!axisLength> (\!axisxlevel,\!axisylevel)
  \fi
  \!axisLaBeL={}}

%


\def\arrow <#1> [#2,#3]{%
  \!ifnextchar<{\!arrow{#1}{#2}{#3}}{\!arrow{#1}{#2}{#3}<\!zpt,\!zpt> }}

\def\!arrow#1#2#3<#4,#5> from #6 #7 to #8 #9 {%
%
  \!xloc=\!M{#8}\!xunit   
  \!yloc=\!M{#9}\!yunit
  \!dxpos=\!xloc  \!dimenA=\!M{#6}\!xunit  \advance \!dxpos -\!dimenA
  \!dypos=\!yloc  \!dimenA=\!M{#7}\!yunit  \advance \!dypos -\!dimenA
  \let\!MAH=\!M
  \!setdimenmode
  \!xshift=#4\relax  \!yshift=#5\relax
  \!reverserotateonly\!xshift\!yshift
  \advance\!xshift\!xloc  \advance\!yshift\!yloc
%
  \!xS=-\!dxpos  \advance\!xS\!xshift
  \!yS=-\!dypos  \advance\!yS\!yshift
  \!start (\!xS,\!yS)
  \!ljoin (\!xshift,\!yshift)
%
  \!Pythag\!dxpos\!dypos\!arclength
  \!divide\!dxpos\!arclength\!dxpos  
  \!dxpos=32\!dxpos  \!removept\!dxpos\!!cos
  \!divide\!dypos\!arclength\!dypos  
  \!dypos=32\!dypos  \!removept\!dypos\!!sin
%
  \!halfhead{#1}{#2}{#3}
  \!halfhead{#1}{-#2}{-#3}
  \let\!M=\!MAH
  \ignorespaces}
%
  \def\!halfhead#1#2#3{%
    \!dimenC=-#1%
    \divide \!dimenC 2 
    \!dimenD=#2\!dimenC
    \!rotate(\!dimenC,\!dimenD)by(\!!cos,\!!sin)to(\!xM,\!yM)
    \!dimenC=-#1
    \!dimenD=#3\!dimenC
    \!dimenD=.5\!dimenD
    \!rotate(\!dimenC,\!dimenD)by(\!!cos,\!!sin)to(\!xE,\!yE)
    \!start (\!xshift,\!yshift)
    \advance\!xM\!xshift  \advance\!yM\!yshift
    \advance\!xE\!xshift  \advance\!yE\!yshift
    \!qjoin (\!xM,\!yM) (\!xE,\!yE) 
    \ignorespaces}

\def\betweenarrows #1#2 from #3 #4 to #5 #6 {%
  \!xloc=\!M{#3}\!xunit  \!xxloc=\!M{#5}\!xunit%
  \!yloc=\!M{#4}\!yunit  \!yyloc=\!M{#6}\!yunit%
  \!dxpos=\!xxloc  \advance\!dxpos by -\!xloc
  \!dypos=\!yyloc  \advance\!dypos by -\!yloc
  \advance\!xloc .5\!dxpos
  \advance\!yloc .5\!dypos
  \let\!MBA=\!M
  \!setdimenmode
  \ifdim\!dypos=\!zpt
    \ifdim\!dxpos<\!zpt \!dxpos=-\!dxpos \fi
    \put {\!lrarrows{\!dxpos}{#1}}#2{} at {\!xloc} {\!yloc}
  \else
    \ifdim\!dxpos=\!zpt
      \ifdim\!dypos<\!zpt \!dypos=-\!zpt \fi
      \put {\!udarrows{\!dypos}{#1}}#2{} at {\!xloc} {\!yloc}
    \fi
  \fi
  \let\!M=\!MBA
  \ignorespaces}

\def\!lrarrows#1#2{
  {\setbox\!boxA=\hbox{$\mkern-2mu\mathord-\mkern-2mu$}%
   \setbox\!boxB=\hbox{$\leftarrow$}\!dimenE=\ht\!boxB
   \setbox\!boxB=\hbox{}\ht\!boxB=2\!dimenE
   \hbox to #1{$\mathord\leftarrow\mkern-6mu
     \cleaders\copy\!boxA\hfil
     \mkern-6mu\mathord-$%
     \kern.4em $\vcenter{\box\!boxB}$$\vcenter{\hbox{#2}}$\kern.4em
     $\mathord-\mkern-6mu
     \cleaders\copy\!boxA\hfil
     \mkern-6mu\mathord\rightarrow$}}}

\def\!udarrows#1#2{
  {\setbox\!boxB=\hbox{#2}%
   \setbox\!boxA=\hbox to \wd\!boxB{\hss$\vert$\hss}%
   \!dimenE=\ht\!boxA \advance\!dimenE \dp\!boxA \divide\!dimenE 2
   \vbox to #1{\offinterlineskip
      \vskip .05556\!dimenE
      \hbox to \wd\!boxB{\hss$\mkern.4mu\uparrow$\hss}\vskip-\!dimenE
      \cleaders\copy\!boxA\vfil
      \vskip-\!dimenE\copy\!boxA
      \vskip\!dimenE\copy\!boxB\vskip.4em
      \copy\!boxA\vskip-\!dimenE
      \cleaders\copy\!boxA\vfil
      \vskip-\!dimenE \hbox to \wd\!boxB{\hss$\mkern.4mu\downarrow$\hss}
      \vskip .05556\!dimenE}}}

%

\def\putbar#1breadth <#2> from #3 #4 to #5 #6 {%
  \!xloc=\!M{#3}\!xunit  \!xxloc=\!M{#5}\!xunit%
  \!yloc=\!M{#4}\!yunit  \!yyloc=\!M{#6}\!yunit%
  \!dypos=\!yyloc  \advance\!dypos by -\!yloc
  \!dimenI=#2  
  \ifdim \!dimenI=\!zpt 
    \putrule#1from {#3} {#4} to {#5} {#6} 
  \else 
    \let\!MBar=\!M
    \!setdimenmode 
    \divide\!dimenI 2
    \ifdim \!dypos=\!zpt             
      \advance \!yloc -\!dimenI 
      \advance \!yyloc \!dimenI
    \else
      \advance \!xloc -\!dimenI 
      \advance \!xxloc \!dimenI
    \fi
    \putrectangle#1corners at {\!xloc} {\!yloc} and {\!xxloc} {\!yyloc}
    \let\!M=\!MBar 
  \fi
  \ignorespaces}

\def\setbars#1breadth <#2> baseline at #3 = #4 {%
  \edef\!barshift{#1}%
  \edef\!barbreadth{#2}%
  \edef\!barorientation{#3}%
  \edef\!barbaseline{#4}%
  \def\!bardobaselabel{\!bardoendlabel}%
  \def\!bardoendlabel{\!barfinish}%
  \let\!drawcurve=\!barcurve
  \!setbars}
\def\!setbars{%
  \futurelet\!nextchar\!!setbars}
\def\!!setbars{%
  \if b\!nextchar
    \def\!!!setbars{\!setbarsbget}%
  \else 
    \if e\!nextchar
      \def\!!!setbars{\!setbarseget}%
    \else
      \def\!!!setbars{\relax}%
    \fi
  \fi
  \!!!setbars}
\def\!setbarsbget baselabels (#1) {%
  \def\!barbaselabelorientation{#1}%
  \def\!bardobaselabel{\!!bardobaselabel}%
  \!setbars}
\def\!setbarseget endlabels (#1) {%
  \edef\!barendlabelorientation{#1}%
  \def\!bardoendlabel{\!!bardoendlabel}%
  \!setbars}

\def\!barcurve #1 #2 {%
  \if y\!barorientation
    \def\!basexarg{#1}%
    \def\!baseyarg{\!barbaseline}%
  \else
    \def\!basexarg{\!barbaseline}%
    \def\!baseyarg{#2}%
  \fi
  \expandafter\putbar\!barshift breadth <\!barbreadth> from {\!basexarg}
    {\!baseyarg} to {#1} {#2}
  \def\!endxarg{#1}%
  \def\!endyarg{#2}%
  \!bardobaselabel}

\def\!!bardobaselabel "#1" {%
  \put {#1}\!barbaselabelorientation{} at {\!basexarg} {\!baseyarg}
  \!bardoendlabel}
 
\def\!!bardoendlabel "#1" {%
  \put {#1}\!barendlabelorientation{} at {\!endxarg} {\!endyarg}
  \!barfinish}

\def\!barfinish{%
  \!ifnextchar/{\!finish}{\!barcurve}}

%
%
%
\def\putrectangle{%
  \!ifnextchar<{\!putrectangle}{\!putrectangle<\!zpt,\!zpt> }}
\def\!putrectangle<#1,#2> corners at #3 #4 and #5 #6 {%
%
  \!xone=\!M{#3}\!xunit  \!xtwo=\!M{#5}\!xunit%
  \!yone=\!M{#4}\!yunit  \!ytwo=\!M{#6}\!yunit%
  \ifdim \!xtwo<\!xone
    \!dimenI=\!xone  \!xone=\!xtwo  \!xtwo=\!dimenI
  \fi
  \ifdim \!ytwo<\!yone
    \!dimenI=\!yone  \!yone=\!ytwo  \!ytwo=\!dimenI
  \fi
  \!dimenI=#1\relax  \advance\!xone\!dimenI  \advance\!xtwo\!dimenI
  \!dimenI=#2\relax  \advance\!yone\!dimenI  \advance\!ytwo\!dimenI
  \let\!MRect=\!M
  \!setdimenmode
%
  \!shaderectangle
%
  \!dimenI=.5\linethickness
  \advance \!xone  -\!dimenI
  \advance \!xtwo   \!dimenI
  \putrule from {\!xone} {\!yone} to {\!xtwo} {\!yone} 
  \putrule from {\!xone} {\!ytwo} to {\!xtwo} {\!ytwo} 
%
  \advance \!xone   \!dimenI
  \advance \!xtwo  -\!dimenI%
  \advance \!yone  -\!dimenI
  \advance \!ytwo   \!dimenI
  \putrule from {\!xone} {\!yone} to {\!xone} {\!ytwo} 
  \putrule from {\!xtwo} {\!yone} to {\!xtwo} {\!ytwo} 
  \let\!M=\!MRect
  \ignorespaces}
 

\def\shaderectanglesoff{%
  \def\!shaderectangle{}%
  \ignorespaces}

\shaderectanglesoff
 
\def\!!shaderectangle{%
  \!dimenA=\!xtwo  \advance \!dimenA -\!xone
  \!dimenB=\!ytwo  \advance \!dimenB -\!yone
  \ifdim \!dimenA<\!dimenB
    \!startvshade (\!xone,\!yone,\!ytwo)
    \!lshade      (\!xtwo,\!yone,\!ytwo)
  \else
    \!starthshade (\!yone,\!xone,\!xtwo)
    \!lshade      (\!ytwo,\!xone,\!xtwo)
  \fi
  \ignorespaces}
  
\def\frame{%
  \!ifnextchar<{\!frame}{\!frame<\!zpt> }}
\long\def\!frame<#1> #2{%
  \beginpicture
    \setcoordinatesystem units <1pt,1pt> point at 0 0 
    \put {#2} [Bl] at 0 0 
    \!dimenA=#1\relax
    \!dimenB=\!wd \advance \!dimenB \!dimenA
    \!dimenC=\!ht \advance \!dimenC \!dimenA
    \!dimenD=\!dp \advance \!dimenD \!dimenA
    \let\!MFr=\!M
    \!setdimenmode
    \putrectangle corners at {-\!dimenA} {-\!dimenD} and {\!dimenB} {\!dimenC}
    \!setcoordmode
    \let\!M=\!MFr
  \endpicture
  \ignorespaces}
 
\def\rectangle <#1> <#2> {%
  \setbox0=\hbox{}\wd0=#1\ht0=#2\frame {\box0}}

%

\def\plot{%
  \!ifnextchar"{\!plotfromfile}{\!drawcurve}}
\def\!plotfromfile"#1"{%
  \expandafter\!drawcurve \input #1 /}

\def\setquadratic{%
  \let\!drawcurve=\!qcurve
  \let\!!Shade=\!!qShade
  \let\!!!Shade=\!!!qShade}

\def\setlinear{%
  \let\!drawcurve=\!lcurve
  \let\!!Shade=\!!lShade
  \let\!!!Shade=\!!!lShade}

\def\sethistograms{%
  \let\!drawcurve=\!hcurve}

\def\!qcurve #1 #2 {%
  \!start (#1,#2)
  \!Qjoin}
\def\!Qjoin#1 #2 #3 #4 {%
  \!qjoin (#1,#2) (#3,#4)             
  \!ifnextchar/{\!finish}{\!Qjoin}}

\def\!lcurve #1 #2 {%
  \!start (#1,#2)
  \!Ljoin}
\def\!Ljoin#1 #2 {%
  \!ljoin (#1,#2)                    
  \!ifnextchar/{\!finish}{\!Ljoin}}

\def\!finish/{\ignorespaces}

\def\!hcurve #1 #2 {%
  \edef\!hxS{#1}%
  \edef\!hyS{#2}%
  \!hjoin}
\def\!hjoin#1 #2 {%
  \putrectangle corners at {\!hxS} {\!hyS} and {#1} {#2}
  \edef\!hxS{#1}%
  \!ifnextchar/{\!finish}{\!hjoin}}

\def\vshade #1 #2 #3 {%
  \!startvshade (#1,#2,#3)
  \!Shadewhat}

\def\hshade #1 #2 #3 {%
  \!starthshade (#1,#2,#3)
  \!Shadewhat}

\def\!Shadewhat{%
  \futurelet\!nextchar\!Shade}
\def\!Shade{%
  \if <\!nextchar
    \def\!nextShade{\!!Shade}%
  \else
    \if /\!nextchar
      \def\!nextShade{\!finish}%
    \else
      \def\!nextShade{\!!!Shade}%
    \fi
  \fi
  \!nextShade}
\def\!!lShade<#1> #2 #3 #4 {%
  \!lshade <#1> (#2,#3,#4)                 
  \!Shadewhat}
\def\!!!lShade#1 #2 #3 {%
  \!lshade (#1,#2,#3)
  \!Shadewhat} 
\def\!!qShade<#1> #2 #3 #4 #5 #6 #7 {%
  \!qshade <#1> (#2,#3,#4) (#5,#6,#7)      
  \!Shadewhat}
\def\!!!qShade#1 #2 #3 #4 #5 #6 {%
  \!qshade (#1,#2,#3) (#4,#5,#6)
  \!Shadewhat} 

\setlinear

\def\setdashpattern <#1>{%
  \def\!Flist{}\def\!Blist{}\def\!UDlist{}%
  \!countA=0
  \!ecfor\!item:=#1\do{%
    \!dimenA=\!item\relax
    \expandafter\!rightappend\the\!dimenA\withCS{\\}\to\!UDlist%
    \advance\!countA  1
    \ifodd\!countA
      \expandafter\!rightappend\the\!dimenA\withCS{\!Rule}\to\!Flist%
      \expandafter\!leftappend\the\!dimenA\withCS{\!Rule}\to\!Blist%
    \else 
      \expandafter\!rightappend\the\!dimenA\withCS{\!Skip}\to\!Flist%
      \expandafter\!leftappend\the\!dimenA\withCS{\!Skip}\to\!Blist%
    \fi}%
  \!leaderlength=\!zpt
  \def\!Rule##1{\advance\!leaderlength  ##1}%
  \def\!Skip##1{\advance\!leaderlength  ##1}%
  \!Flist%
  \ifdim\!leaderlength>\!zpt 
  \else
    \def\!Flist{\!Skip{24in}}\def\!Blist{\!Skip{24in}}\ignorespaces
    \def\!UDlist{\\{\!zpt}\\{24in}}\ignorespaces
    \!leaderlength=24in
  \fi
  \!dashingon}

\def\!dashingon{%
  \def\!advancedashing{\!!advancedashing}%
  \def\!drawlinearsegment{\!lineardashed}%
  \def\!puthline{\!putdashedhline}%
  \def\!putvline{\!putdashedvline}%
  \ignorespaces}%
\def\!dashingoff{%
  \def\!advancedashing{\relax}%
  \def\!drawlinearsegment{\!linearsolid}%
  \def\!puthline{\!putsolidhline}%
  \def\!putvline{\!putsolidvline}%
  \ignorespaces}

\def\setdots{%
  \!ifnextchar<{\!setdots}{\!setdots<5pt>}}
\def\!setdots<#1>{%
  \!dimenB=#1\advance\!dimenB -\plotsymbolspacing
  \ifdim\!dimenB<\!zpt
    \!dimenB=\!zpt
  \fi
\setdashpattern <\plotsymbolspacing,\!dimenB>}
 
\def\setdotsnear <#1> for <#2>{%
  \!dimenB=#2\relax  \advance\!dimenB -.05pt  
  \!dimenC=#1\relax  \!countA=\!dimenC 
  \!dimenD=\!dimenB  \advance\!dimenD .5\!dimenC  \!countB=\!dimenD
  \divide \!countB  \!countA
  \ifnum 1>\!countB 
    \!countB=1
  \fi
  \divide\!dimenB  \!countB
  \setdots <\!dimenB>}
 
\def\setdashes{%
  \!ifnextchar<{\!setdashes}{\!setdashes<5pt>}}
\def\!setdashes<#1>{\setdashpattern <#1,#1>}
 
\def\setdashesnear <#1> for <#2>{%
  \!dimenB=#2\relax  
  \!dimenC=#1\relax  \!countA=\!dimenC 
  \!dimenD=\!dimenB  \advance\!dimenD .5\!dimenC  \!countB=\!dimenD
  \divide \!countB  \!countA
  \ifodd \!countB 
  \else 
    \advance \!countB  1
  \fi
  \divide\!dimenB  \!countB
  \setdashes <\!dimenB>}
 
\def\setsolid{%
  \def\!Flist{\!Rule{24in}}\def\!Blist{\!Rule{24in}}%
  \def\!UDlist{\\{24in}\\{\!zpt}}%
  \!dashingoff}  
\setsolid


 
  
 
\def\!divide#1#2#3{%
  \!dimenB=#1
  \!dimenC=#2
  \!dimenD=\!dimenB
  \divide \!dimenD \!dimenC
  \!dimenA=\!dimenD
  \multiply\!dimenD \!dimenC
  \advance\!dimenB -\!dimenD
  \!dimenD=\!dimenC
    \ifdim\!dimenD<\!zpt \!dimenD=-\!dimenD 
  \fi
  \ifdim\!dimenD<64pt
    \!divstep[\!tfs]\!divstep[\!tfs]%
  \else 
    \!!divide
  \fi
  #3=\!dimenA\ignorespaces}

\def\!!divide{%
  \ifdim\!dimenD<256pt
    \!divstep[64]\!divstep[32]\!divstep[32]%
  \else 
    \!divstep[8]\!divstep[8]\!divstep[8]\!divstep[8]\!divstep[8]%
    \!dimenA=2\!dimenA
  \fi}

\def\!divstep[#1]{
  \!dimenB=#1\!dimenB
  \!dimenD=\!dimenB
    \divide \!dimenD by \!dimenC
  \!dimenA=#1\!dimenA
    \advance\!dimenA by \!dimenD%
  \multiply\!dimenD by \!dimenC
    \advance\!dimenB by -\!dimenD}
 
\def\Divide <#1> by <#2> forming <#3> {%
  \!divide{#1}{#2}{#3}}

 
 

 

\def\ellipticalarc axes ratio #1:#2 #3 degrees from #4 #5 center at #6 #7 {%
  \!angle=#3pt\relax
  \ifdim\!angle>\!zpt 
    \def\!sign{}
  \else 
    \def\!sign{-}\!angle=-\!angle
  \fi
  \!xxloc=\!M{#6}\!xunit
  \!yyloc=\!M{#7}\!yunit     
  \!xxS=\!M{#4}\!xunit
  \!yyS=\!M{#5}\!yunit
  \advance\!xxS -\!xxloc
  \advance\!yyS -\!yyloc
  \!divide\!xxS{#1pt}\!xxS 
  \!divide\!yyS{#2pt}\!yyS 
  \let\!MC=\!M
  \!setdimenmode
  \!xS=#1\!xxS  \advance\!xS\!xxloc
  \!yS=#2\!yyS  \advance\!yS\!yyloc
  \!start (\!xS,\!yS)%
  \!loop\ifdim\!angle>14.9999pt
    \!rotate(\!xxS,\!yyS)by(\!cos,\!sign\!sin)to(\!xxM,\!yyM) 
    \!rotate(\!xxM,\!yyM)by(\!cos,\!sign\!sin)to(\!xxE,\!yyE)
    \!xM=#1\!xxM  \advance\!xM\!xxloc  \!yM=#2\!yyM  \advance\!yM\!yyloc
    \!xE=#1\!xxE  \advance\!xE\!xxloc  \!yE=#2\!yyE  \advance\!yE\!yyloc
    \!qjoin (\!xM,\!yM) (\!xE,\!yE)
    \!xxS=\!xxE  \!yyS=\!yyE 
    \advance \!angle -15pt
  \repeat
  \ifdim\!angle>\!zpt
    \!angle=100.53096\!angle
    \divide \!angle 360 
    \!sinandcos\!angle\!!sin\!!cos
    \!rotate(\!xxS,\!yyS)by(\!!cos,\!sign\!!sin)to(\!xxM,\!yyM) 
    \!rotate(\!xxM,\!yyM)by(\!!cos,\!sign\!!sin)to(\!xxE,\!yyE)
    \!xM=#1\!xxM  \advance\!xM\!xxloc  \!yM=#2\!yyM  \advance\!yM\!yyloc
    \!xE=#1\!xxE  \advance\!xE\!xxloc  \!yE=#2\!yyE  \advance\!yE\!yyloc
    \!qjoin (\!xM,\!yM) (\!xE,\!yE)
  \fi
  \let\!M=\!MC
  \ignorespaces}

\def\!rotate(#1,#2)by(#3,#4)to(#5,#6){%
  \!dimenA=#3#1\advance \!dimenA -#4#2
  \!dimenB=#3#2\advance \!dimenB  #4#1
  \divide \!dimenA 32  \divide \!dimenB 32 
  #5=\!dimenA  #6=\!dimenB
  \ignorespaces}
\def\!sin{4.17684}
\def\!cos{31.72624}

\def\!sinandcos#1#2#3{%
 \!dimenD=#1
 \!dimenA=\!dimenD
 \!dimenB=32pt
 \!removept\!dimenD\!value
 \!dimenC=\!dimenD
 \!dimenC=\!value\!dimenC \divide\!dimenC by 64 
 \advance\!dimenB by -\!dimenC
 \!dimenC=\!value\!dimenC \divide\!dimenC by 96 
 \advance\!dimenA by -\!dimenC
 \!dimenC=\!value\!dimenC \divide\!dimenC by 128 
 \advance\!dimenB by \!dimenC%
 \!removept\!dimenA#2
 \!removept\!dimenB#3
 \ignorespaces}




\def\putrule#1from #2 #3 to #4 #5 {%
  \!xloc=\!M{#2}\!xunit  \!xxloc=\!M{#4}\!xunit%
  \!yloc=\!M{#3}\!yunit  \!yyloc=\!M{#5}\!yunit%
  \!dxpos=\!xxloc  \advance\!dxpos by -\!xloc
  \!dypos=\!yyloc  \advance\!dypos by -\!yloc
  \ifdim\!dypos=\!zpt
    \def\!!Line{\!puthline{#1}}\ignorespaces
  \else
    \ifdim\!dxpos=\!zpt
      \def\!!Line{\!putvline{#1}}\ignorespaces
    \else 
       \def\!!Line{}
    \fi
  \fi
  \let\!ML=\!M
  \!setdimenmode
  \!!Line%
  \let\!M=\!ML
  \ignorespaces}

\def\!putsolidhline#1{%
  \ifdim\!dxpos>\!zpt 
    \put{\!hline\!dxpos}#1[l] at {\!xloc} {\!yloc}
  \else 
    \put{\!hline{-\!dxpos}}#1[l] at {\!xxloc} {\!yyloc}
  \fi
  \ignorespaces}
 
\def\!putsolidvline#1{%
  \ifdim\!dypos>\!zpt 
    \put{\!vline\!dypos}#1[b] at {\!xloc} {\!yloc}
  \else 
    \put{\!vline{-\!dypos}}#1[b] at {\!xxloc} {\!yyloc}
  \fi
  \ignorespaces}
 
\def\!hline#1{\hbox to #1{\leaders \hrule height\linethickness\hfill}}
\def\!vline#1{\vbox to #1{\leaders \vrule width\linethickness\vfill}}

\def\!putdashedhline#1{%
  \ifdim\!dxpos>\!zpt 
    \!DLsetup\!Flist\!dxpos
    \put{\hbox to \!totalleaderlength{\!hleaders}\!hpartialpattern\!Rtrunc}
      #1[l] at {\!xloc} {\!yloc} 
  \else 
    \!DLsetup\!Blist{-\!dxpos}
    \put{\!hpartialpattern\!Ltrunc\hbox to \!totalleaderlength{\!hleaders}}
      #1[r] at {\!xloc} {\!yloc} 
  \fi
  \ignorespaces}
 
\def\!putdashedvline#1{%
  \!dypos=-\!dypos
  \ifdim\!dypos>\!zpt 
    \!DLsetup\!Flist\!dypos 
    \put{\vbox{\vbox to \!totalleaderlength{\!vleaders}
      \!vpartialpattern\!Rtrunc}}#1[t] at {\!xloc} {\!yloc} 
  \else 
    \!DLsetup\!Blist{-\!dypos}
    \put{\vbox{\!vpartialpattern\!Ltrunc
      \vbox to \!totalleaderlength{\!vleaders}}}#1[b] at {\!xloc} {\!yloc} 
  \fi
  \ignorespaces}

\def\!DLsetup#1#2{
  \let\!RSlist=#1
  \!countB=#2
  \!countA=\!leaderlength
  \divide\!countB by \!countA
  \!totalleaderlength=\!countB\!leaderlength
  \!Rresiduallength=#2%
  \advance \!Rresiduallength by -\!totalleaderlength
  \!Lresiduallength=\!leaderlength
  \advance \!Lresiduallength by -\!Rresiduallength
  \ignorespaces}
 
\def\!hleaders{%
  \def\!Rule##1{\vrule height\linethickness width##1}%
  \def\!Skip##1{\hskip##1}%
  \leaders\hbox{\!RSlist}\hfill}
 
\def\!hpartialpattern#1{%
  \!dimenA=\!zpt \!dimenB=\!zpt 
  \def\!Rule##1{#1{##1}\vrule height\linethickness width\!dimenD}%
  \def\!Skip##1{#1{##1}\hskip\!dimenD}%
  \!RSlist}
 
\def\!vleaders{%
  \def\!Rule##1{\hrule width\linethickness height##1}%
  \def\!Skip##1{\vskip##1}%
  \leaders\vbox{\!RSlist}\vfill}
 
\def\!vpartialpattern#1{%
  \!dimenA=\!zpt \!dimenB=\!zpt 
  \def\!Rule##1{#1{##1}\hrule width\linethickness height\!dimenD}%
  \def\!Skip##1{#1{##1}\vskip\!dimenD}%
  \!RSlist}
 
\def\!Rtrunc#1{\!trunc{#1}>\!Rresiduallength}
\def\!Ltrunc#1{\!trunc{#1}<\!Lresiduallength}
 
\def\!trunc#1#2#3{%
  \!dimenA=\!dimenB         
  \advance\!dimenB by #1%
  \!dimenD=\!dimenB  \ifdim\!dimenD#2#3\!dimenD=#3\fi
  \!dimenC=\!dimenA  \ifdim\!dimenC#2#3\!dimenC=#3\fi
  \advance \!dimenD by -\!dimenC}

\def\!start (#1,#2){%
  \!plotxorigin=\!xorigin  \advance \!plotxorigin by \!plotsymbolxshift
  \!plotyorigin=\!yorigin  \advance \!plotyorigin by \!plotsymbolyshift
  \!xS=\!M{#1}\!xunit \!yS=\!M{#2}\!yunit
  \!rotateaboutpivot\!xS\!yS
  \!copylist\!UDlist\to\!!UDlist
  \!getnextvalueof\!downlength\from\!!UDlist
  \!distacross=\!zpt
  \!intervalno=0 
  \global\totalarclength=\!zpt
  \ignorespaces}

\def\!ljoin (#1,#2){%
  \advance\!intervalno by 1
  \!xE=\!M{#1}\!xunit \!yE=\!M{#2}\!yunit
  \!rotateaboutpivot\!xE\!yE
  \!xdiff=\!xE \advance \!xdiff by -\!xS
  \!ydiff=\!yE \advance \!ydiff by -\!yS
  \!Pythag\!xdiff\!ydiff\!arclength
  \global\advance \totalarclength by \!arclength%
  \!drawlinearsegment
  \!xS=\!xE \!yS=\!yE
  \ignorespaces}

\def\!linearsolid{%
  \!npoints=\!arclength
  \!countA=\plotsymbolspacing
  \divide\!npoints by \!countA
  \ifnum \!npoints<1 
    \!npoints=1 
  \fi
  \divide\!xdiff by \!npoints
  \divide\!ydiff by \!npoints
  \!xpos=\!xS \!ypos=\!yS
  \loop\ifnum\!npoints>-1
    \!plotifinbounds
    \advance \!xpos by \!xdiff
    \advance \!ypos by \!ydiff
    \advance \!npoints by -1
  \repeat
  \ignorespaces}

\def\!lineardashed{%
  \ifdim\!distacross>\!arclength
    \advance \!distacross by -\!arclength  
  \else
    \loop\ifdim\!distacross<\!arclength
      \!divide\!distacross\!arclength\!dimenA
      \!removept\!dimenA\!t
      \!xpos=\!t\!xdiff \advance \!xpos by \!xS
      \!ypos=\!t\!ydiff \advance \!ypos by \!yS
      \!plotifinbounds
      \advance\!distacross by \plotsymbolspacing
      \!advancedashing
    \repeat  
    \advance \!distacross by -\!arclength
  \fi
  \ignorespaces}

\def\!!advancedashing{%
  \advance\!downlength by -\plotsymbolspacing
  \ifdim \!downlength>\!zpt
  \else
    \advance\!distacross by \!downlength
    \!getnextvalueof\!uplength\from\!!UDlist
    \advance\!distacross by \!uplength
    \!getnextvalueof\!downlength\from\!!UDlist
  \fi}

\def\inboundscheckoff{%
  \def\!plotifinbounds{\!plot(\!xpos,\!ypos)}%
  \def\!initinboundscheck{\relax}\ignorespaces}
 
\inboundscheckoff
 
\def\!!plotifinbounds{%
  \ifdim \!xpos<\!checkleft
  \else
    \ifdim \!xpos>\!checkright
    \else
      \ifdim \!ypos<\!checkbot
      \else
         \ifdim \!ypos>\!checktop
         \else
           \!plot(\!xpos,\!ypos)
         \fi 
      \fi
    \fi
  \fi}

\def\!!initinboundscheck{%
  \!checkleft=\!arealloc     \advance\!checkleft by \!xorigin
  \!checkright=\!arearloc    \advance\!checkright by \!xorigin
  \!checkbot=\!areabloc      \advance\!checkbot by \!yorigin
  \!checktop=\!areatloc      \advance\!checktop by \!yorigin}

%


\def\!logten#1#2{%
  \expandafter\!!logten#1\!nil
  \!removept\!dimenF#2%
  \ignorespaces}

\def\!!logten#1#2\!nil{%
  \if -#1%
    \!dimenF=\!zpt
    \def\!next{\ignorespaces}%
  \else
    \if +#1%
      \def\!next{\!!logten#2\!nil}%
    \else
      \if .#1%
        \def\!next{\!!logten0.#2\!nil}%
      \else
        \def\!next{\!!!logten#1#2..\!nil}%
      \fi
    \fi
  \fi
  \!next}

\def\!!!logten#1#2.#3.#4\!nil{%
  \!dimenF=1pt 
  \if 0#1%
    \!!logshift#3pt 
  \else 
    \!logshift#2/
    \!dimenE=#1.#2#3pt 
  \fi 
  \ifdim \!dimenE<\!rootten
    \multiply \!dimenE 10 
    \advance  \!dimenF -1pt
  \fi
  \!dimenG=\!dimenE
    \advance\!dimenG 10pt
  \advance\!dimenE -10pt 
  \multiply\!dimenE 10 
  \!divide\!dimenE\!dimenG\!dimenE
  \!removept\!dimenE\!t
  \!dimenG=\!t\!dimenE
  \!removept\!dimenG\!tt
  \!dimenH=\!tt\!tenAe
    \divide\!dimenH 100
  \advance\!dimenH \!tenAc
  \!dimenH=\!tt\!dimenH
    \divide\!dimenH 100   
  \advance\!dimenH \!tenAa
  \!dimenH=\!t\!dimenH
    \divide\!dimenH 100 
  \advance\!dimenF \!dimenH}

\def\!logshift#1{%
  \if #1/%
    \def\!next{\ignorespaces}%
  \else
    \advance\!dimenF 1pt 
    \def\!next{\!logshift}%
  \fi 
  \!next}
 
 \def\!!logshift#1{%
   \advance\!dimenF -1pt
   \if 0#1%
     \def\!next{\!!logshift}%
   \else
     \if p#1%
       \!dimenF=1pt
       \def\!next{\!dimenE=1p}%
     \else
       \def\!next{\!dimenE=#1.}%
     \fi
   \fi
   \!next}

\def\beginpicture{%
  \setbox\!picbox=\hbox\bgroup%
  \!xleft=\maxdimen  
  \!xright=-\maxdimen
  \!ybot=\maxdimen
  \!ytop=-\maxdimen}
 
\def\endpicture{%
  \ifdim\!xleft=\maxdimen
    \!xleft=\!zpt \!xright=\!zpt \!ybot=\!zpt \!ytop=\!zpt 
  \fi
  \global\!Xleft=\!xleft \global\!Xright=\!xright
  \global\!Ybot=\!ybot \global\!Ytop=\!ytop
  \egroup%
  \ht\!picbox=\!Ytop  \dp\!picbox=-\!Ybot
  \ifdim\!Ybot>\!zpt
  \else 
    \ifdim\!Ytop<\!zpt
      \!Ybot=\!Ytop
    \else
      \!Ybot=\!zpt
    \fi
  \fi
  \hbox{\kern-\!Xleft\lower\!Ybot\box\!picbox\kern\!Xright}}
 
\def\endpicturesave <#1,#2>{%
  \endpicture \global #1=\!Xleft \global #2=\!Ybot \ignorespaces}

\def\setcoordinatesystem{%
  \!ifnextchar{u}{\!getlengths }
    {\!getlengths units <\!xunit,\!yunit>}}
\def\!getlengths units <#1,#2>{%
  \!xunit=#1\relax
  \!yunit=#2\relax
  \!ifcoordmode 
    \let\!SCnext=\!SCccheckforRP
  \else
    \let\!SCnext=\!SCdcheckforRP
  \fi
  \!SCnext}
\def\!SCccheckforRP{%
  \!ifnextchar{p}{\!cgetreference }
    {\!cgetreference point at {\!xref} {\!yref} }}
\def\!cgetreference point at #1 #2 {%
  \edef\!xref{#1}\edef\!yref{#2}%
  \!xorigin=\!xref\!xunit  \!yorigin=\!yref\!yunit  
  \!initinboundscheck 
  \ignorespaces}
\def\!SCdcheckforRP{%
  \!ifnextchar{p}{\!dgetreference}%
    {\ignorespaces}}
\def\!dgetreference point at #1 #2 {%
  \!xorigin=#1\relax  \!yorigin=#2\relax
  \ignorespaces}

\long\def\put#1#2 at #3 #4 {%
  \!setputobject{#1}{#2}%
  \!xpos=\!M{#3}\!xunit  \!ypos=\!M{#4}\!yunit  
  \!rotateaboutpivot\!xpos\!ypos%
  \advance\!xpos -\!xorigin  \advance\!xpos -\!xshift
  \advance\!ypos -\!yorigin  \advance\!ypos -\!yshift
  \kern\!xpos\raise\!ypos\box\!putobject\kern-\!xpos%
  \!doaccounting\ignorespaces}
 
\long\def\multiput #1#2 at {%
  \!setputobject{#1}{#2}%
  \!ifnextchar"{\!putfromfile}{\!multiput}}
\def\!putfromfile"#1"{%
  \expandafter\!multiput \input #1 /}
\def\!multiput{%
  \futurelet\!nextchar\!!multiput}
\def\!!multiput{%
  \if *\!nextchar
    \def\!nextput{\!alsoby}%
  \else
    \if /\!nextchar
      \def\!nextput{\!finishmultiput}%
    \else
      \def\!nextput{\!alsoat}%
    \fi
  \fi
  \!nextput}
\def\!finishmultiput/{%
  \setbox\!putobject=\hbox{}%
  \ignorespaces}
 
\def\!alsoat#1 #2 {%
  \!xpos=\!M{#1}\!xunit  \!ypos=\!M{#2}\!yunit  
  \!rotateaboutpivot\!xpos\!ypos%
  \advance\!xpos -\!xorigin  \advance\!xpos -\!xshift
  \advance\!ypos -\!yorigin  \advance\!ypos -\!yshift
  \kern\!xpos\raise\!ypos\copy\!putobject\kern-\!xpos%
  \!doaccounting
  \!multiput}
 
\def\!alsoby*#1 #2 #3 {%
  \!dxpos=\!M{#2}\!xunit \!dypos=\!M{#3}\!yunit 
  \!rotateonly\!dxpos\!dypos
  \!ntemp=#1%
  \!!loop\ifnum\!ntemp>0
    \advance\!xpos by \!dxpos  \advance\!ypos by \!dypos
    \kern\!xpos\raise\!ypos\copy\!putobject\kern-\!xpos%
    \advance\!ntemp by -1
  \repeat
  \!doaccounting 
  \!multiput}
 
\def\accountingon{\def\!doaccounting{\!!doaccounting}\ignorespaces}

\accountingon
\def\!!doaccounting{%
  \!xtemp=\!xpos  
  \!ytemp=\!ypos
  \ifdim\!xtemp<\!xleft 
     \!xleft=\!xtemp 
  \fi
  \advance\!xtemp by  \!wd 
  \ifdim\!xright<\!xtemp 
    \!xright=\!xtemp
  \fi
  \advance\!ytemp by -\!dp
  \ifdim\!ytemp<\!ybot  
    \!ybot=\!ytemp
  \fi
  \advance\!ytemp by  \!dp
  \advance\!ytemp by  \!ht 
  \ifdim\!ytemp>\!ytop  
    \!ytop=\!ytemp  
  \fi}
 
\long\def\!setputobject#1#2{%
  \setbox\!putobject=\hbox{#1}%
  \!ht=\ht\!putobject  \!dp=\dp\!putobject  \!wd=\wd\!putobject
  \wd\!putobject=\!zpt
  \!xshift=.5\!wd   \!yshift=.5\!ht   \advance\!yshift by -.5\!dp
  \edef\!putorientation{#2}%
  \expandafter\!SPOreadA\!putorientation[]\!nil%
  \expandafter\!SPOreadB\!putorientation<\!zpt,\!zpt>\!nil\ignorespaces}
 
\def\!SPOreadA#1[#2]#3\!nil{\!etfor\!orientation:=#2\do\!SPOreviseshift}
 
\def\!SPOreadB#1<#2,#3>#4\!nil{\advance\!xshift by -#2\advance\!yshift by -#3}
 
\def\!SPOreviseshift{%
  \if l\!orientation 
    \!xshift=\!zpt
  \else 
    \if r\!orientation 
      \!xshift=\!wd
    \else 
      \if b\!orientation
        \!yshift=-\!dp
      \else 
        \if B\!orientation 
          \!yshift=\!zpt
        \else 
          \if t\!orientation 
            \!yshift=\!ht
          \fi 
        \fi
      \fi
    \fi
  \fi}

\long\def\!dimenput#1#2(#3,#4){%
  \!setputobject{#1}{#2}%
  \!xpos=#3\advance\!xpos by -\!xshift
  \!ypos=#4\advance\!ypos by -\!yshift
  \kern\!xpos\raise\!ypos\box\!putobject\kern-\!xpos%
  \!doaccounting\ignorespaces}

\def\!setdimenmode{%
  \let\!M=\!M!!\ignorespaces}
\def\!setcoordmode{%
  \let\!M=\!M!\ignorespaces}
\def\!ifcoordmode{%
  \ifx \!M \!M!}
\def\!ifdimenmode{%
  \ifx \!M \!M!!}
\def\!M!#1#2{#1#2} 
\def\!M!!#1#2{#1}
\!setcoordmode
\let\setdimensionmode=\!setdimenmode
\let\setcoordinatemode=\!setcoordmode




\def\!stack[#1]{%
  \let\!lglue=\hfill \let\!rglue=\hfill
  \expandafter\let\csname !#1glue\endcsname=\relax
  \!ifnextchar<{\!!stack}{\!!stack<\stackleading>}}
\def\!!stack<#1>#2{%
  \vbox{\def\!valueslist{}\!ecfor\!value:=#2\do{%
    \expandafter\!rightappend\!value\withCS{\\}\to\!valueslist}%
    \!lop\!valueslist\to\!value
    \let\\=\cr\lineskiplimit=\maxdimen\lineskip=#1%
    \baselineskip=-1000pt\halign{\!lglue##\!rglue\cr \!value\!valueslist\cr}}%
  \ignorespaces}


\def\!lines[#1]#2{%
  \let\!lglue=\hfill \let\!rglue=\hfill
  \expandafter\let\csname !#1glue\endcsname=\relax
  \vbox{\halign{\!lglue##\!rglue\cr #2\crcr}}%
  \ignorespaces}


\def\!Lines[#1]#2{%
  \let\!lglue=\hfill \let\!rglue=\hfill
  \expandafter\let\csname !#1glue\endcsname=\relax
  \vtop{\halign{\!lglue##\!rglue\cr #2\crcr}}%
  \ignorespaces}

 
 
 
\def\setplotsymbol(#1#2){%
  \!setputobject{#1}{#2}
  \setbox\!plotsymbol=\box\!putobject%
  \!plotsymbolxshift=\!xshift 
  \!plotsymbolyshift=\!yshift 
  \ignorespaces}
 
\font\fiverm=cmr5
\setplotsymbol({\fiverm .})

 
\def\!!plot(#1,#2){%
  \!dimenA=-\!plotxorigin \advance \!dimenA by #1
  \!dimenB=-\!plotyorigin \advance \!dimenB by #2
  \kern\!dimenA\raise\!dimenB\copy\!plotsymbol\kern-\!dimenA%
  \ignorespaces}
 
\def\!!!plot(#1,#2){%
  \!dimenA=-\!plotxorigin \advance \!dimenA by #1
  \!dimenB=-\!plotyorigin \advance \!dimenB by #2
  \kern\!dimenA\raise\!dimenB\copy\!plotsymbol\kern-\!dimenA%
  \!countE=\!dimenA
  \!countF=\!dimenB
  \immediate\write\!replotfile{\the\!countE,\the\!countF.}%
  \ignorespaces}

\def\savelinesandcurves on "#1" {%
  \immediate\closeout\!replotfile
  \immediate\openout\!replotfile=#1%
  \let\!plot=\!!!plot}

\def\dontsavelinesandcurves {%
  \let\!plot=\!!plot}
\dontsavelinesandcurves

{\catcode`\%=11\xdef\!Commentsignal{
\def\writesavefile#1 {%
  \immediate\write\!replotfile{\!Commentsignal #1}%
  \ignorespaces}

\def\replot"#1" {%
  \expandafter\!replot\input #1 /}
\def\!replot#1,#2. {%
  \!dimenA=#1sp
  \kern\!dimenA\raise#2sp\copy\!plotsymbol\kern-\!dimenA
  \futurelet\!nextchar\!!replot}
\def\!!replot{%
  \if /\!nextchar 
    \def\!next{\!finish}%
  \else
    \def\!next{\!replot}%
  \fi
  \!next}


 
 
\def\!Pythag#1#2#3{%
  \!dimenE=#1\relax                                     
  \ifdim\!dimenE<\!zpt 
    \!dimenE=-\!dimenE 
  \fi
  \!dimenF=#2\relax
  \ifdim\!dimenF<\!zpt 
    \!dimenF=-\!dimenF 
  \fi
  \advance \!dimenF by \!dimenE
  \ifdim\!dimenF=\!zpt 
    \!dimenG=\!zpt
  \else 
    \!divide{8\!dimenE}\!dimenF\!dimenE
    \advance\!dimenE by -4pt
      \!dimenE=2\!dimenE
    \!removept\!dimenE\!!t
    \!dimenE=\!!t\!dimenE
    \advance\!dimenE by 64pt
    \divide \!dimenE by 2
    \!dimenH=7pt
    \!!Pythag\!!Pythag\!!Pythag
    \!removept\!dimenH\!!t
    \!dimenG=\!!t\!dimenF
    \divide\!dimenG by 8
  \fi
  #3=\!dimenG
  \ignorespaces}

\def\!!Pythag{
  \!divide\!dimenE\!dimenH\!dimenI
  \advance\!dimenH by \!dimenI
    \divide\!dimenH by 2}

\def\placehypotenuse for <#1> and <#2> in <#3> {%
  \!Pythag{#1}{#2}{#3}}

 
 
 
\def\!qjoin (#1,#2) (#3,#4){%
  \advance\!intervalno by 1
  \!ifcoordmode
    \edef\!xmidpt{#1}\edef\!ymidpt{#2}%
  \else
    \!dimenA=#1\relax \edef\!xmidpt{\the\!dimenA}%
    \!dimenA=#2\relax \edef\!ymidpt{\the\!dimenA}%
  \fi
  \!xM=\!M{#1}\!xunit  \!yM=\!M{#2}\!yunit   \!rotateaboutpivot\!xM\!yM
  \!xE=\!M{#3}\!xunit  \!yE=\!M{#4}\!yunit   \!rotateaboutpivot\!xE\!yE
%
  \!dimenA=\!xM  \advance \!dimenA by -\!xS
  \!dimenB=\!xE  \advance \!dimenB by -\!xM
  \!xB=3\!dimenA \advance \!xB by -\!dimenB
  \!xC=2\!dimenB \advance \!xC by -2\!dimenA
%
  \!dimenA=\!yM  \advance \!dimenA by -\!yS%
  \!dimenB=\!yE  \advance \!dimenB by -\!yM%
  \!yB=3\!dimenA \advance \!yB by -\!dimenB%
  \!yC=2\!dimenB \advance \!yC by -2\!dimenA%
%
  \!xprime=\!xB  \!yprime=\!yB
  \!dxprime=.5\!xC  \!dyprime=.5\!yC
  \!getf \!midarclength=\!dimenA
  \!getf \advance \!midarclength by 4\!dimenA
  \!getf \advance \!midarclength by \!dimenA
  \divide \!midarclength by 12
%
  \!arclength=\!dimenA
  \!getf \advance \!arclength by 4\!dimenA
  \!getf \advance \!arclength by \!dimenA
  \divide \!arclength by 12
  \advance \!arclength by \!midarclength
  \global\advance \totalarclength by \!arclength
%
%
  \ifdim\!distacross>\!arclength 
    \advance \!distacross by -\!arclength
  \else
    \!initinverseinterp
    \loop\ifdim\!distacross<\!arclength
      \!inverseinterp
      \!xpos=\!t\!xC \advance\!xpos by \!xB
        \!xpos=\!t\!xpos \advance \!xpos by \!xS
      \!ypos=\!t\!yC \advance\!ypos by \!yB
        \!ypos=\!t\!ypos \advance \!ypos by \!yS
      \!plotifinbounds
      \advance\!distacross \plotsymbolspacing
      \!advancedashing
    \repeat  
    \advance \!distacross by -\!arclength
  \fi
  \!xS=\!xE
  \!yS=\!yE
  \ignorespaces}

\def\!getf{\!Pythag\!xprime\!yprime\!dimenA%
  \advance\!xprime by \!dxprime
  \advance\!yprime by \!dyprime}

\def\!initinverseinterp{%
  \ifdim\!arclength>\!zpt
    \!divide{8\!midarclength}\!arclength\!dimenE
    \ifdim\!dimenE<\!wmin \!setinverselinear
    \else 
      \ifdim\!dimenE>\!wmax \!setinverselinear
      \else
        \def\!inverseinterp{\!inversequad}\ignorespaces
%
%
         \!removept\!dimenE\!Ew
         \!dimenF=-\!Ew\!dimenE
         \advance\!dimenF by 32pt
         \!dimenG=8pt 
         \advance\!dimenG by -\!dimenE
         \!dimenG=\!Ew\!dimenG
         \!divide\!dimenF\!dimenG\!beta
         \!gamma=1pt
         \advance \!gamma by -\!beta
      \fi
    \fi
  \fi
  \ignorespaces}

\def\!inversequad{%
  \!divide\!distacross\!arclength\!dimenG
  \!removept\!dimenG\!v
  \!dimenG=\!v\!gamma
  \advance\!dimenG by \!beta
  \!dimenG=\!v\!dimenG
  \!removept\!dimenG\!t}

\def\!setinverselinear{%
  \def\!inverseinterp{\!inverselinear}%
  \divide\!dimenE by 8 \!removept\!dimenE\!t
  \!countC=\!intervalno \multiply \!countC 2
  \!countB=\!countC     \advance \!countB -1
  \!countA=\!countB     \advance \!countA -1
  \wlog{\the\!countB th point (\!xmidpt,\!ymidpt) being plotted 
    doesn't lie in the}%
  \wlog{ middle third of the arc between the \the\!countA th 
    and \the\!countC th points:}%
  \wlog{ [arc length \the\!countA\space to \the\!countB]/[arc length 
    \the \!countA\space to \the\!countC]=\!t.}%
  \ignorespaces}
 
\def\!inverselinear{%
  \!divide\!distacross\!arclength\!dimenG
  \!removept\!dimenG\!t}

 

\def\startrotation{%
  \let\!rotateaboutpivot=\!!rotateaboutpivot
  \let\!rotateonly=\!!rotateonly
  \!ifnextchar{b}{\!getsincos }%
    {\!getsincos by {\!cosrotationangle} {\!sinrotationangle} }}
\def\!getsincos by #1 #2 {%
  \edef\!cosrotationangle{#1}%
  \edef\!sinrotationangle{#2}%
  \!ifcoordmode 
    \let\!ROnext=\!ccheckforpivot
  \else
    \let\!ROnext=\!dcheckforpivot
  \fi
  \!ROnext}
\def\!ccheckforpivot{%
  \!ifnextchar{a}{\!cgetpivot}%
    {\!cgetpivot about {\!xpivotcoord} {\!ypivotcoord} }}
\def\!cgetpivot about #1 #2 {%
  \edef\!xpivotcoord{#1}%
  \edef\!ypivotcoord{#2}%
  \!xpivot=#1\!xunit  \!ypivot=#2\!yunit
  \ignorespaces}
\def\!dcheckforpivot{%
  \!ifnextchar{a}{\!dgetpivot}{\ignorespaces}}
\def\!dgetpivot about #1 #2 {%
  \!xpivot=#1\relax  \!ypivot=#2\relax
  \ignorespaces}

\def\stoprotation{%
  \let\!rotateaboutpivot=\!!!rotateaboutpivot
  \let\!rotateonly=\!!!rotateonly
  \ignorespaces}
 
\def\!!rotateaboutpivot#1#2{%
  \!dimenA=#1\relax  \advance\!dimenA -\!xpivot
  \!dimenB=#2\relax  \advance\!dimenB -\!ypivot
  \!dimenC=\!cosrotationangle\!dimenA
    \advance \!dimenC -\!sinrotationangle\!dimenB
  \!dimenD=\!cosrotationangle\!dimenB
    \advance \!dimenD  \!sinrotationangle\!dimenA
  \advance\!dimenC \!xpivot  \advance\!dimenD \!ypivot
  #1=\!dimenC  #2=\!dimenD
  \ignorespaces}

\def\!!rotateonly#1#2{%
  \!dimenA=#1\relax  \!dimenB=#2\relax 
  \!dimenC=\!cosrotationangle\!dimenA
    \advance \!dimenC -\!rotsign\!sinrotationangle\!dimenB
  \!dimenD=\!cosrotationangle\!dimenB
    \advance \!dimenD  \!rotsign\!sinrotationangle\!dimenA
  #1=\!dimenC  #2=\!dimenD
  \ignorespaces}
\def\!rotsign{}
\def\!!!rotateaboutpivot#1#2{\relax}
\def\!!!rotateonly#1#2{\relax}
\stoprotation

\def\!reverserotateonly#1#2{%
  \def\!rotsign{-}%
  \!rotateonly{#1}{#2}%
  \def\!rotsign{}%
  \ignorespaces}

\def\!getspan span <#1>{%
  \!dshade=#1\relax
  \!ifcoordmode 
    \let\!GRnext=\!GRccheckforAP
  \else
    \let\!GRnext=\!GRdcheckforAP
  \fi
  \!GRnext}
\def\!GRccheckforAP{%
  \!ifnextchar{p}{\!cgetanchor }
    {\!cgetanchor point at {\!xshadesave} {\!yshadesave} }}
\def\!cgetanchor point at #1 #2 {%
  \edef\!xshadesave{#1}\edef\!yshadesave{#2}%
  \!xshade=\!xshadesave\!xunit  \!yshade=\!yshadesave\!yunit
  \ignorespaces}
\def\!GRdcheckforAP{%
  \!ifnextchar{p}{\!dgetanchor}%
    {\ignorespaces}}
\def\!dgetanchor point at #1 #2 {%
  \!xshade=#1\relax  \!yshade=#2\relax
  \ignorespaces}

\def\setshadesymbol{%
  \!ifnextchar<{\!setshadesymbol}{\!setshadesymbol<,,,> }}

\def\!setshadesymbol <#1,#2,#3,#4> (#5#6){%
  \!setputobject{#5}{#6}%
  \setbox\!shadesymbol=\box\!putobject%
  \!shadesymbolxshift=\!xshift \!shadesymbolyshift=\!yshift
%
  \!dimenA=\!xshift \advance\!dimenA \!smidge
  \!override\!dimenA{#1}\!lshrinkage%
  \!dimenA=\!wd \advance \!dimenA -\!xshift
    \advance\!dimenA \!smidge
    \!override\!dimenA{#2}\!rshrinkage
  \!dimenA=\!dp \advance \!dimenA \!yshift
    \advance\!dimenA \!smidge
    \!override\!dimenA{#3}\!bshrinkage
  \!dimenA=\!ht \advance \!dimenA -\!yshift
    \advance\!dimenA \!smidge
    \!override\!dimenA{#4}\!tshrinkage
  \ignorespaces}
\def\!smidge{-.2pt}%

\def\!override#1#2#3{%
  \edef\!!override{#2}%
  \ifx \!!override\empty
    #3=#1\relax
  \else
    \if z\!!override
      #3=\!zpt
    \else
      \ifx \!!override\!blankz
        #3=\!zpt
      \else
        #3=#2\relax
      \fi
    \fi
  \fi
  \ignorespaces}
\def\!blankz{ z}

\setshadesymbol ({\fiverm .})

\def\!startvshade#1(#2,#3,#4){%
  \let\!!xunit=\!xunit%
  \let\!!yunit=\!yunit%
  \let\!!xshade=\!xshade%
  \let\!!yshade=\!yshade%
  \def\!getshrinkages{\!vgetshrinkages}%
  \let\!setshadelocation=\!vsetshadelocation%
  \!xS=\!M{#2}\!!xunit
  \!ybS=\!M{#3}\!!yunit
  \!ytS=\!M{#4}\!!yunit
  \!shadexorigin=\!xorigin  \advance \!shadexorigin \!shadesymbolxshift
  \!shadeyorigin=\!yorigin  \advance \!shadeyorigin \!shadesymbolyshift
  \ignorespaces}
 
\def\!starthshade#1(#2,#3,#4){%
  \let\!!xunit=\!yunit%
  \let\!!yunit=\!xunit%
  \let\!!xshade=\!yshade%
  \let\!!yshade=\!xshade%
  \def\!getshrinkages{\!hgetshrinkages}%
  \let\!setshadelocation=\!hsetshadelocation%
  \!xS=\!M{#2}\!!xunit
  \!ybS=\!M{#3}\!!yunit
  \!ytS=\!M{#4}\!!yunit
  \!shadexorigin=\!xorigin  \advance \!shadexorigin \!shadesymbolxshift
  \!shadeyorigin=\!yorigin  \advance \!shadeyorigin \!shadesymbolyshift
  \ignorespaces}

\def\!lattice#1#2#3#4#5{%
  \!dimenA=#1
  \!dimenB=#2
  \!countB=\!dimenB
%
  \!dimenC=#3
  \advance\!dimenC -\!dimenA
  \!countA=\!dimenC
  \divide\!countA \!countB
  \ifdim\!dimenC>\!zpt
    \!dimenD=\!countA\!dimenB
    \ifdim\!dimenD<\!dimenC
      \advance\!countA 1 
    \fi
  \fi
  \!dimenC=\!countA\!dimenB
    \advance\!dimenC \!dimenA
  #4=\!countA
  #5=\!dimenC
  \ignorespaces}

\def\!qshade#1(#2,#3,#4)#5(#6,#7,#8){%
  \!xM=\!M{#2}\!!xunit
  \!ybM=\!M{#3}\!!yunit
  \!ytM=\!M{#4}\!!yunit
  \!xE=\!M{#6}\!!xunit
  \!ybE=\!M{#7}\!!yunit
  \!ytE=\!M{#8}\!!yunit
  \!getcoeffs\!xS\!ybS\!xM\!ybM\!xE\!ybE\!ybB\!ybC
  \!getcoeffs\!xS\!ytS\!xM\!ytM\!xE\!ytE\!ytB\!ytC
  \def\!getylimits{\!qgetylimits}%
  \!shade{#1}\ignorespaces}
 
\def\!lshade#1(#2,#3,#4){%
  \!xE=\!M{#2}\!!xunit
  \!ybE=\!M{#3}\!!yunit
  \!ytE=\!M{#4}\!!yunit
  \!dimenE=\!xE  \advance \!dimenE -\!xS
  \!dimenC=\!ytE \advance \!dimenC -\!ytS
  \!divide\!dimenC\!dimenE\!ytB
  \!dimenC=\!ybE \advance \!dimenC -\!ybS
  \!divide\!dimenC\!dimenE\!ybB
  \def\!getylimits{\!lgetylimits}%
  \!shade{#1}\ignorespaces}
 
\def\!getcoeffs#1#2#3#4#5#6#7#8{%
  \!dimenC=#4\advance \!dimenC -#2
  \!dimenE=#3\advance \!dimenE -#1
  \!divide\!dimenC\!dimenE\!dimenF
  \!dimenC=#6\advance \!dimenC -#4
  \!dimenH=#5\advance \!dimenH -#3
  \!divide\!dimenC\!dimenH\!dimenG
  \advance\!dimenG -\!dimenF
  \advance \!dimenH \!dimenE
  \!divide\!dimenG\!dimenH#8
  \!removept#8\!t
  #7=-\!t\!dimenE
  \advance #7\!dimenF
  \ignorespaces}

\def\!shade#1{%
  \!getshrinkages#1<,,,>\!nil
  \advance \!dimenE \!xS
  \!lattice\!!xshade\!dshade\!dimenE
    \!parity\!xpos
  \!dimenF=-\!dimenF
    \advance\!dimenF \!xE
  \!loop\!not{\ifdim\!xpos>\!dimenF}
    \!shadecolumn%
    \advance\!xpos \!dshade
    \advance\!parity 1
  \repeat
  \!xS=\!xE
  \!ybS=\!ybE
  \!ytS=\!ytE
  \ignorespaces}

\def\!vgetshrinkages#1<#2,#3,#4,#5>#6\!nil{%
  \!override\!lshrinkage{#2}\!dimenE
  \!override\!rshrinkage{#3}\!dimenF
  \!override\!bshrinkage{#4}\!dimenG
  \!override\!tshrinkage{#5}\!dimenH
  \ignorespaces}
\def\!hgetshrinkages#1<#2,#3,#4,#5>#6\!nil{%
  \!override\!lshrinkage{#2}\!dimenG
  \!override\!rshrinkage{#3}\!dimenH
  \!override\!bshrinkage{#4}\!dimenE
  \!override\!tshrinkage{#5}\!dimenF
  \ignorespaces}

\def\!shadecolumn{%
  \!dxpos=\!xpos
  \advance\!dxpos -\!xS
  \!removept\!dxpos\!dx
  \!getylimits
  \advance\!ytpos -\!dimenH
  \advance\!ybpos \!dimenG
  \!yloc=\!!yshade
  \ifodd\!parity 
     \advance\!yloc \!dshade
  \fi
  \!lattice\!yloc{2\!dshade}\!ybpos%
    \!countA\!ypos
  \!dimenA=-\!shadexorigin \advance \!dimenA \!xpos
  \loop\!not{\ifdim\!ypos>\!ytpos}
    \!setshadelocation
    \!rotateaboutpivot\!xloc\!yloc%
    \!dimenA=-\!shadexorigin \advance \!dimenA \!xloc
    \!dimenB=-\!shadeyorigin \advance \!dimenB \!yloc
    \kern\!dimenA \raise\!dimenB\copy\!shadesymbol \kern-\!dimenA
    \advance\!ypos 2\!dshade
  \repeat
  \ignorespaces}
 
\def\!qgetylimits{%
  \!dimenA=\!dx\!ytC              
  \advance\!dimenA \!ytB
  \!ytpos=\!dx\!dimenA
  \advance\!ytpos \!ytS
  \!dimenA=\!dx\!ybC              
  \advance\!dimenA \!ybB
  \!ybpos=\!dx\!dimenA
  \advance\!ybpos \!ybS}
 
\def\!lgetylimits{%
  \!ytpos=\!dx\!ytB
  \advance\!ytpos \!ytS
  \!ybpos=\!dx\!ybB
  \advance\!ybpos \!ybS}
 
\def\!vsetshadelocation{
  \!xloc=\!xpos
  \!yloc=\!ypos}
\def\!hsetshadelocation{
  \!xloc=\!ypos
  \!yloc=\!xpos}





\def\!axisticks {%
  \def\!nextkeyword##1 {%
    \expandafter\ifx\csname !ticks##1\endcsname \relax
      \def\!next{\!fixkeyword{##1}}%
    \else
      \def\!next{\csname !ticks##1\endcsname}%
    \fi
    \!next}%
  \!axissetup
    \def\!axissetup{\relax}%
  \edef\!ticksinoutsign{\!ticksinoutSign}%
  \!ticklength=\longticklength
  \!tickwidth=\linethickness
  \!gridlinestatus
  \!setticktransform
  \!maketick
  \!tickcase=0
  \def\!LTlist{}%
  \!nextkeyword}

\def\ticksout{%
  \def\!ticksinoutSign{+}}

\ticksout

\def\nogridlines{%
  \def\!gridlinestatus{\!gridlinestoofalse}}
\nogridlines

\def\loggedticks{%
  \def\!setticktransform{\let\!ticktransform=\!logten}}
\def\unloggedticks{%
  \def\!setticktransform{\let\!ticktransform=\!donothing}}
\def\!donothing#1#2{\def#2{#1}}
\unloggedticks

\expandafter\def\csname !ticks/\endcsname{%
  \!not {\ifx \!LTlist\empty}
    \!placetickvalues
  \fi
  \def\!tickvalueslist{}%
  \def\!LTlist{}%
  \expandafter\csname !axis/\endcsname}

\def\!maketick{%
  \setbox\!boxA=\hbox{%
    \beginpicture
      \!setdimenmode
      \setcoordinatesystem point at {\!zpt} {\!zpt}   
      \linethickness=\!tickwidth
      \ifdim\!ticklength>\!zpt
        \putrule from {\!zpt} {\!zpt} to
          {\!ticksinoutsign\!tickxsign\!ticklength}
          {\!ticksinoutsign\!tickysign\!ticklength}
      \fi
      \if!gridlinestoo
        \putrule from {\!zpt} {\!zpt} to
          {-\!tickxsign\!xaxislength} {-\!tickysign\!yaxislength}
      \fi
    \endpicturesave <\!Xsave,\!Ysave>}%
    \wd\!boxA=\!zpt}
  
\def\!ticksin{%
  \def\!ticksinoutsign{-}%
  \!maketick
  \!nextkeyword}

\def\!ticksout{%
  \def\!ticksinoutsign{+}%
  \!maketick
  \!nextkeyword}

\def\!tickslength<#1> {%
  \!ticklength=#1\relax
  \!maketick
  \!nextkeyword}

\def\!tickslong{%
  \!tickslength<\longticklength> }

\def\!ticksshort{%
  \!tickslength<\shortticklength> }

\def\!tickswidth<#1> {%
  \!tickwidth=#1\relax
  \!maketick
  \!nextkeyword}

\def\!ticksandacross{%
  \!gridlinestootrue
  \!maketick
  \!nextkeyword}

\def\!ticksbutnotacross{%
  \!gridlinestoofalse
  \!maketick
  \!nextkeyword}

\def\!tickslogged{%
  \let\!ticktransform=\!logten
  \!nextkeyword}

\def\!ticksunlogged{%
  \let\!ticktransform=\!donothing
  \!nextkeyword}

\def\!ticksunlabeled{%
  \!tickcase=0
  \!nextkeyword}

\def\!ticksnumbered{%
  \!tickcase=1
  \!nextkeyword}

\def\!tickswithvalues#1/ {%
  \edef\!tickvalueslist{#1! /}%
  \!tickcase=2
  \!nextkeyword}

\def\!ticksquantity#1 {%
  \ifnum #1>1
    \!updatetickoffset
    \!countA=#1\relax
    \advance \!countA -1
    \!ticklocationincr=\!axisLength
      \divide \!ticklocationincr \!countA
    \!ticklocation=\!axisstart
    \loop \!not{\ifdim \!ticklocation>\!axisend}
      \!placetick\!ticklocation
      \ifcase\!tickcase
          \relax 
        \or
          \relax 
        \or
          \expandafter\!gettickvaluefrom\!tickvalueslist
          \edef\!tickfield{{\the\!ticklocation}{\!value}}%
          \expandafter\!listaddon\expandafter{\!tickfield}\!LTlist%
      \fi
      \advance \!ticklocation \!ticklocationincr
    \repeat
  \fi
  \!nextkeyword}

\def\!ticksat#1 {%
  \!updatetickoffset
  \edef\!Loc{#1}%
  \if /\!Loc
    \def\next{\!nextkeyword}%
  \else
    \!ticksincommon
    \def\next{\!ticksat}%
  \fi
  \next}    
      
\def\!ticksfrom#1 to #2 by #3 {%
  \!updatetickoffset
  \edef\!arg{#3}%
  \expandafter\!separate\!arg\!nil
  \!scalefactor=1
  \expandafter\!countfigures\!arg/
  \edef\!arg{#1}%
  \!scaleup\!arg by\!scalefactor to\!countE
  \edef\!arg{#2}%
  \!scaleup\!arg by\!scalefactor to\!countF
  \edef\!arg{#3}%
  \!scaleup\!arg by\!scalefactor to\!countG
  \loop \!not{\ifnum\!countE>\!countF}
    \ifnum\!scalefactor=1
      \edef\!Loc{\the\!countE}%
    \else
      \!scaledown\!countE by\!scalefactor to\!Loc
    \fi
    \!ticksincommon
    \advance \!countE \!countG
  \repeat
  \!nextkeyword}

\def\!updatetickoffset{%
  \!dimenA=\!ticksinoutsign\!ticklength
  \ifdim \!dimenA>\!offset
    \!offset=\!dimenA
  \fi}

\def\!placetick#1{%
  \if!xswitch
    \!xpos=#1\relax
    \!ypos=\!axisylevel
  \else
    \!xpos=\!axisxlevel
    \!ypos=#1\relax
  \fi
  \advance\!xpos \!Xsave
  \advance\!ypos \!Ysave
  \kern\!xpos\raise\!ypos\copy\!boxA\kern-\!xpos
  \ignorespaces}

\def\!gettickvaluefrom#1 #2 /{%
  \edef\!value{#1}%
  \edef\!tickvalueslist{#2 /}%
  \ifx \!tickvalueslist\!endtickvaluelist
    \!tickcase=0
  \fi}
\def\!endtickvaluelist{! /}

\def\!ticksincommon{%
  \!ticktransform\!Loc\!t
  \!ticklocation=\!t\!!unit
  \advance\!ticklocation -\!!origin
  \!placetick\!ticklocation
  \ifcase\!tickcase
    \relax 
  \or 
    \ifdim\!ticklocation<-\!!origin
      \edef\!Loc{$\!Loc$}%
    \fi
    \edef\!tickfield{{\the\!ticklocation}{\!Loc}}%
    \expandafter\!listaddon\expandafter{\!tickfield}\!LTlist%
  \or 
    \expandafter\!gettickvaluefrom\!tickvalueslist
    \edef\!tickfield{{\the\!ticklocation}{\!value}}%
    \expandafter\!listaddon\expandafter{\!tickfield}\!LTlist%
  \fi}

\def\!separate#1\!nil{%
  \!ifnextchar{-}{\!!separate}{\!!!separate}#1\!nil}
\def\!!separate-#1\!nil{%
  \def\!sign{-}%
  \!!!!separate#1..\!nil}
\def\!!!separate#1\!nil{%
  \def\!sign{+}%
  \!!!!separate#1..\!nil}
\def\!!!!separate#1.#2.#3\!nil{%
  \def\!arg{#1}%
  \ifx\!arg\!empty
    \!countA=0
  \else
    \!countA=\!arg
  \fi
  \def\!arg{#2}%
  \ifx\!arg\!empty
    \!countB=0
  \else
    \!countB=\!arg
  \fi}
 
\def\!countfigures#1{%
  \if #1/%
    \def\!next{\ignorespaces}%
  \else
    \multiply\!scalefactor 10
    \def\!next{\!countfigures}%
  \fi
  \!next}

\def\!scaleup#1by#2to#3{%
  \expandafter\!separate#1\!nil
  \multiply\!countA #2\relax
  \advance\!countA \!countB
  \if -\!sign
    \!countA=-\!countA
  \fi
  #3=\!countA
  \ignorespaces}

\def\!scaledown#1by#2to#3{%
  \!countA=#1\relax
  \ifnum \!countA<0 
    \def\!sign{-}
    \!countA=-\!countA
  \else
    \def\!sign{}%
  \fi
  \!countB=\!countA
  \divide\!countB #2\relax
  \!countC=\!countB
    \multiply\!countC #2\relax
  \advance \!countA -\!countC
  \edef#3{\!sign\the\!countB.}
  \!countC=\!countA 
  \ifnum\!countC=0 
    \!countC=1
  \fi
  \multiply\!countC 10
  \!loop \ifnum #2>\!countC
    \edef#3{#3\!zero}%
    \multiply\!countC 10
  \repeat
  \edef#3{#3\the\!countA}
  \ignorespaces}

\def\!placetickvalues{%
  \advance\!offset \tickstovaluesleading
  \if!xswitch
    \setbox\!boxA=\hbox{%
      \def\\##1##2{%
        \!dimenput {##2} [B] (##1,\!axisylevel)}%
      \beginpicture 
        \!LTlist
      \endpicturesave <\!Xsave,\!Ysave>}%
    \!dimenA=\!axisylevel
      \advance\!dimenA -\!Ysave
      \advance\!dimenA \!tickysign\!offset
      \if -\!tickysign
        \advance\!dimenA -\ht\!boxA
      \else
        \advance\!dimenA  \dp\!boxA
      \fi
    \advance\!offset \ht\!boxA 
      \advance\!offset \dp\!boxA
    \!dimenput {\box\!boxA} [Bl] <\!Xsave,\!Ysave> (\!zpt,\!dimenA)
  \else
    \setbox\!boxA=\hbox{%
      \def\\##1##2{%
        \!dimenput {##2} [r] (\!axisxlevel,##1)}%
      \beginpicture 
        \!LTlist
      \endpicturesave <\!Xsave,\!Ysave>}%
    \!dimenA=\!axisxlevel
      \advance\!dimenA -\!Xsave
      \advance\!dimenA \!tickxsign\!offset
      \if -\!tickxsign
        \advance\!dimenA -\wd\!boxA
      \fi
    \advance\!offset \wd\!boxA
    \!dimenput {\box\!boxA} [Bl] <\!Xsave,\!Ysave> (\!dimenA,\!zpt)
  \fi}

\normalgraphs
\catcode`!=12 


\def\limsup{\mathop{\overline{\hbox{lim}}}}
\def\liminf{\mathop{\underline{\hbox{lim}}}}




\def\abs#1{\left\vert #1 \right\vert}      
\def\grint#1{\left\lfloor#1\right\rfloor}  

\def\set#1{\left\{#1\right\}}
\def\theset#1#2{\left\{#1:#2\right\}}            



\def\union{\cup}

\def\Union{\bigcup}


\def\R{{\bf R}}

\def\Z{{\bf Z}}



\def\Prob#1{P\left(#1\right)}

\def\Expect#1{E\left[#1\right]}

\def\Var#1{\hbox{Var}\left(#1\right)}
\def\Cov#1#2{\hbox{Cov}\left(#1,#2\right)}


 
\def\cvgindistn{\;\mathop{\rightarrow^{\mkern-18mu D}\,}\;} 


\def\endofproof{~\hfill\vrule height8pt width3.5pt depth-1pt}

\def\narrowmargins{\hoffset-7mm
                   \hsize178mm
                   \voffset0mm
                   \vsize230mm}

\def\headroom#1{\vphantom{\raise3pt\hbox{#1}\lower3pt\hbox{#1}}#1}


\font\titlefont=cmbx10 scaled 1400
\font\largefont=cmr10 scaled 1200
\font\smallfont=cmr8
\font\smallit=cmti8
\font\smallbf=cmbx8

\def\papertitle#1{\centerline{\titlefont#1}\vskip12pt}
\def\paperauthor#1{\centerline{\largefont#1}\vskip5pt}
\def\paperabstract#1{{\vskip30pt\leftskip=50pt\rightskip=50pt%
       {\noindent{\smallbf Abstract. }{\smallfont#1}}\vskip20pt}}


\def\narrowmargins{\hoffset-7mm
                   \hsize178mm
                   \voffset0mm
                   \vsize230mm}

\def\startbibliography{\bigskip\centerline{\largefont References}\medskip}
\def\bibentry#1#2#3#4{\item{\smallbf[#1]\ }\smallfont#2, {\smallit#3, }%
         #4\par}
\def\cite#1{{\bf[#1]}}

\newbox\endaddresses
\setbox\endaddresses=\vbox{}
\def\authoraddress#1#2{\setbox\endaddresses=\vbox{\unvbox\endaddresses\bigskip
     {\noindent\smallfont#1}\par
     {\noindent\smallit E-mail address: }{\smallfont#2}\par }}

\def\section#1{\bigbreak\centerline{\bf#1}\medskip}
\def\chunk#1{\medbreak\noindent{\bf#1}}
\def\subchunk#1{\medbreak\noindent{\it#1}}

\def\Bazaraa{1}
\def\BramsFishburn{2}
\def\BergLepelley{3}
\def\Chamberlin{4}
\def\Durrett{5}
\def\Gibbard{6}
\def\KimRoush{7}
\def\LepelleyMbih{8}
\def\Nitzan{9}
\def\PapadimitriouSteiglitz{10}
\def\ParkerRardin{11}
\def\PritchardSlinko{12}
\def\PritchardWilsonExactResults{13}
\def\Saari{14}
\def\Satterthwaite{15}
\def\WilsonPritchardVolumes{16}

\def\CLTN{1}
\def\CLTS{2}
\def\ignoreties{3}
\def\ignoreNlimit{4}
\def\ignoreintegrality{5}
\def\ignoreall{6}
\def\manipbonly{7}
\def\Tionly{8}
\def\zsufficient{9}
\def\zprog{10}
\def\dualprog{11}
\def\pluralitydominated{12}
\def\antipluralityundominated{13}
\def\Bordadomination{14}

\def\fourcandidateexamples{1}
\def\gwgthreefour{2}
\def\gwgfivesix{3}
\def\gwgtentwenty{4}

\def\Tpref{\bar T}
\def\onev{{\bf 1}}
\def\Tpba{\Tpref_{\beta a}}
\def\Tbeta{T_{\beta}}
\def\Tba{T_{ba}}
\def\Tb{T_b}
\def\Vp{V_{\scriptstyle\rm plurality}}
\def\Vap{V_{\scriptstyle\rm antiplurality}}
\def\Vhalfmappr{V_{\scriptstyle\grint{m/2}{\rm -approval}}}
\def\subborda{{\scriptstyle\rm Borda}}

\narrowmargins 

\papertitle{Asymptotics of the minimum manipulating coalition size}
\papertitle{for positional voting rules under IC behaviour}
\bigskip

\vskip30pt
\paperauthor{Geoffrey Pritchard${}^1$, University of Auckland}
\authoraddress{Geoffrey Pritchard, Dept. of Statistics, University of Auckland, Private Bag 92019, Auckland, New Zealand}
              {g.pritchard@auckland.ac.nz}
\paperauthor{Mark C. Wilson, University of Auckland}
\authoraddress{Mark C. Wilson, Dept. of Computer Science, University of Auckland, Private Bag 92019, Auckland, New Zealand}
              {mcw@cs.auckland.ac.nz}

\footnote{}{${}^1$ Corresponding author. Dept. of Statistics, University of Auckland, Private Bag 92019, Auckland, New Zealand;
                g.pritchard@auckland.ac.nz; +64 9 373 7599 ext. 87400}

\paperabstract{
We consider the problem of manipulation of elections using positional voting rules
under Impartial Culture voter behaviour.
We consider both the logical possibility of coalitional manipulation, and the number
of voters that must be recruited to form a manipulating coalition.
It is shown that the manipulation problem may be well approximated by a
very simple linear program in two variables.
This permits a comparative analysis of the asymptotic (large-population)
manipulability of the various rules. 
It is seen that the manipulation resistance of positional rules
with 5 or 6 (or more) candidates is quite different from the more commonly analyzed
3- and 4-candidate cases.
              }

\bigskip
\bigskip
{\it Key words and phrases:} scoring rule, social choice, manipulation,
Impartial Culture, Borda, plurality, anti-plurality, asymptotic, probability, linear programming.

\bigskip
\bigskip
{\it JEL Classification Numbers:} D71, D72.


\section{1. Introduction}

In 1973--75 Gibbard and Satterthwaite published a fundamental impossibility
theorem which states that every non-dictatorial social choice function, whose range
contains at least three alternatives, at certain profiles can be manipulated by a
single individual voter \cite{\Gibbard, \Satterthwaite}. After that, the natural question arose: if there
are no perfect rules, which ones are the best, i.e. least manipulable? 
To this question there can be no absolute answer -- it depends both on the behaviour
of the voters, and on the measure used to quantify the term ``manipulability".

Among models of voter behaviour, the following two have gained the most attention
(\cite{\BergLepelley,\Saari}).
The Impartial Culture (IC) model assumes that voters are independent, and that each
voter is equally likely to vote for any candidate.
The Impartial Anonymous Culture (IAC) model assumes some degree of dependency. 
This paper concerns itself with the IC model.

Among measures of manipulability, the most popular is the probability that the votes fall
in such a way as to create the (coalitional or individual) ``logical possibility of manipulation"
(\cite{\Chamberlin,\KimRoush,\LepelleyMbih,\Nitzan,\PritchardSlinko,\PritchardWilsonExactResults,\Saari}).
This means that some coalition of voters (or individual voter) with incentive to do so
can change the election result by voting insincerely. 
Counterthreats are not considered -- we assume that the manipulator(s) are not
opposed by the other, naive voters.
The probability of manipulability has been especially well-studied for the important
class of positional (scoring) voting rules, and significant progress has been made in comparing them.
In his seminal paper \cite{\Saari}, Saari showed that in his ``geometric'' model, Borda's rule is
the  least manipulable for the three-alternative case in relation to micro manipulation,
but that this does not extend to the case of four alternatives.

In this paper, we further refine this notion of manipulability by considering the sizes of the
coalitions involved. Intuitively, a rule is more resistant to manipulation if
many voters must be recruited to assemble the manipulating coalition, and less resistant
if only a few voters are required. We may thus consider the probability that a
coalition {\it of at most $k$ voters} can manipulate ($k=1,2,\ldots$).
Equivalently, we study the probability distribution of the size of the smallest
manipulating coalition (a random variable). Similar ideas are explored, in a more
limited way, in \cite{\PritchardSlinko} and \cite{\PritchardWilsonExactResults}.

We use the following notation and assumptions throughout.
An election is held to choose one from among $m$ candidates ($m\geq3$).
There are $n$ voters, who hold opinions according to the IC model.
That is, each voter (independently) chooses one of the $m!$ possible types
(preference orders on the candidates), each type being equally likely to be chosen.
The election uses the positional voting rule with score vector
$w=(w_1,\ldots,w_m)$, where $1=w_1\geq w_2\geq\cdots\geq w_m=0$.
That is, a vote ranking candidate $\alpha$ in $i$th place contributes $w_i$ to the
score of $\alpha$, and the candidate with the greatest total score is declared the
winner. The possibility of a tie for first place will not be considered in this
paper, as Proposition \ignoreties\ makes it largely irrelevant; it is
discussed in detail in \cite{\PritchardWilsonExactResults}.
We aim to describe the limiting probability distribution of the minimum manipulating
coalition size as $n\to\infty$, and to use this as a criterion for comparing the rules.

The remainder of this paper is organized as follows. 
The next section records some basic results regarding IC behaviour in large populations.
The long section 3 is the theoretical core of the paper: it defines the manipulation problem
as an integer linear program and then, through a series of simplifications, shows
how this may be replaced by a much simpler linear program.
Readers impatient to reach the main results may wish to skip much of this section at a
first reading.
The desired limiting probability distributions are derived in section 4,
and used in section 5 to compare the rules. Section 6 contains some conclusions.

\section{2. Asymptotic results for large electorates of IC voters}

Let $C$ be the set of candidates, and $T$ the set of all voter types (i.e. all permutations of $C$).
Let $N=(N_t)_{t\in T}$ be the random vector giving the number of voters of each type
(so $\sum_{t\in T} N_t = n$). For IC voter behaviour, $N$ has a multinomial probability distribution
with mean $(n/m!)\onev$. (Here and subsequently, the notation $\onev$ is used to denote a vector
whose entries are all $1$.)
Under the asymptotic conditions of interest to us, this may be approximated by a
multivariate normal distribution.

\chunk{Proposition \CLTN.}
$$ {N - np\onev\over\sqrt{n}} \cvgindistn N(0,\Sigma) ,
  $$
where $p=1/m!$, and $\Sigma$ is the matrix with entries
$$ \Sigma_{st} = \cases{ p(1-p), & if $s=t$ \cr
                         -p^2, & if $s\neq t$ .\cr
                       } 
  $$

\subchunk{Remark.} 
Here and in the rest of this paper, the notation $\cvgindistn$ denotes convergence in
distribution (see \cite{\Durrett, Ch. 2}).
Note that the limiting multivariate normal distribution is degenerate, being
supported on $\theset{(q_t)_{t\in T}}{\sum_{t\in T} q_t =0}$.

\subchunk{Proof.} Since IC voters choose their types at random, and independently, we have
$N=\sum_{i=1}^n X_i$, where $X_1,\ldots,X_n$ are independent and have probability distribution
assigning probability $1/m!$ to each of the unit vectors of $\R^T$. Note that
$\Expect{X_1}=p\onev$ and the covariance matrix of $X_1$ is $\Sigma$. The result then follows by the
central limit theorem (\cite{\Durrett, p.170}).
\endofproof

\subchunk{}
We use the notation $\sigma_t(\alpha)$ for the contribution to
candidate $\alpha$'s score made by a vote of type $t$ (so if $t$ ranks $\alpha$ in $i$th
place, then $\sigma_t(\alpha)=w_i$). The total score of $\alpha$ is then
$$ \abs{\alpha} = \sum_{t\in T} N_t \sigma_t(\alpha) .
  $$
Let $S=(\abs{\alpha})_{\alpha\in C}$ be the vector of candidates' scores (the ``scoreboard").
Proposition \CLTN\ immediately gives a central limit result for $S$, too.

\chunk{Proposition \CLTS.}
$$ {S - n\bar w\onev\over\sqrt{n}} \cvgindistn \sigma_w \left({m\over m-1}\right)^{1/2} (Z - \bar Z \onev) ,
  $$
where $\bar w$, $\sigma_w$ are the mean and standard deviation of the score vector $w$
(i.e. $\bar w = (w_1+\cdots+w_m)/m$ and $\sigma_w^2 = (w_1^2+\cdots+w_m^2)/m - (\bar w)^2$);
$Z$ is a vector, indexed by $C$, of independent standard normal random variables;
and $\bar Z = {\scriptstyle 1\over\scriptstyle m} \sum_{\alpha} Z_{\alpha}$.

\subchunk{Proof.}
Let $Y=(Y_t)_{t\in T}\sim N(0,\Sigma)$. From Proposition \CLTN\ we have
$$ {S - n\bar w\onev\over\sqrt{n}} \cvgindistn U ,
  $$
where $U_{\alpha} = \sum_{t\in T} Y_t \sigma_t(\alpha)$. It only remains to show that $U$
and $\sigma_w \left({m\over m-1}\right)^{1/2} (Z - \bar Z \onev)$ have the same
multivariate normal distribution. For this, it suffices to observe that they have the
same mean (zero), variances, and covariances. It is routine to check that
$$ \Var{\sum_{t\in T} Y_t \sigma_t(\alpha)} \,=\, \sigma_w^2 \,=\,
   \sigma_w^2\left({m\over m-1}\right)\Var{Z_{\alpha} - \bar Z}
  $$
for $\alpha\in C$, and for distinct $\alpha,\beta\in C$
$$ \Cov{\sum_{t\in T} Y_t \sigma_t(\alpha)}{\sum_{t\in T} Y_t \sigma_t(\beta)}
   \,=\, {-\sigma_w^2\over m-1} \,=\,
   \sigma_w^2\left({m\over m-1}\right)\Cov{Z_{\alpha} - \bar Z}{Z_{\beta} - \bar Z} .
  $$
\endofproof

\subchunk{}
Proposition \CLTS\ implies that under IC behaviour in large electorates, the average candidate's score will be
of order $n$, but the variability among the scores will be of order only $\sqrt{n}$.
Consequently, most elections will result in all candidates receiving relatively similar
scores, and the margin of victory will be small.
However, exact ties in the scores become increasingly rare as the number of voters increases.
The following result establishes this formally.

\chunk{Proposition \ignoreties.}
$$ \Prob{\hbox{all candidates' scores are numerically distinct}} \to 1 \qquad\hbox{ as } n\to\infty .
  $$
\subchunk{Proof.}
For two distinct candidates $\alpha$ and $\beta$, we have
$$ \set{\abs{\alpha}\neq\abs{\beta}} = \set{{S - n\bar w 1\over\sqrt{n}} \in G} ,
  $$
where $G=\theset{s\in\R^C}{s_{\alpha}\neq s_{\beta}}$.
By Proposition \CLTS\ (and \cite{\Durrett, p.87}), since $G$ is an open set we have
$$ \liminf_n \Prob{{S - n\bar w 1\over\sqrt{n}} \in G} 
 = \Prob{\sigma_w \left({m\over m-1}\right)^{1/2}(Z-\bar Z\onev)\in G}
 = \Prob{Z_{\alpha}\neq Z_{\beta}}
 = 1 .
  $$
That is, $\Prob{\abs{\alpha}=\abs{\beta}} \to0$ as $n\to\infty$.
The result follows since
$$ \Prob{\Union_{\alpha\neq\beta}\set{\abs{\alpha}=\abs{\beta}}} \leq \sum_{\alpha\neq\beta} \Prob{\abs{\alpha}=\abs{\beta}} .
  $$
\endofproof

\section{3. Approximations of the minimum manipulating coalition size.}

In this section, we formulate the coalitional manipulation problem as an integer linear program.
We then show, via a series of simplifying steps, that this is well approximated by a
much simpler linear program in which there are only two variables, and in which the
constraint set does not depend on the voting situation, but only on the voting rule.
Since the proofs involved are fairly lengthy, we summarize the steps involved before
embarking on them:
\smallbreak
\item{$\bullet$} The problem of assembling the smallest possible manipulating coalition
can be expressed as an integer linear program, in which the variables are the numbers of
voters of each type ($x_t$) to recruit;
\item{$\bullet$} The integrality and upper bound ($x_t\leq N_t$) constraints of this program
may be ignored (Propositions \ignoreNlimit\ and \ignoreintegrality).
\item{$\bullet$} Only manipulations in favour of the second-placegetter need be considered
(Proposition \manipbonly).
\item{$\bullet$} Coalition recruiting may be limited to those voters who rank the two leading candidates
$a$ and $b$ adjacent, as these voters are best able to manipulate (Proposition \Tionly).
\item{$\bullet$} We may find the minimum coalition size by considering only the members'
(sincere) rankings of $a$ and $b$, without regard to how they rank other candidates
(Propositions \zsufficient\ and \zprog).
This reduces the problem to a mere linear program with $m-1$ variables and two constraints.
\item{$\bullet$} Replacing this linear program with its dual gives us two variables and
$m-1$ constraints.
\medbreak

To specify an attempted coalitional manipulation, we must specify for each $t\in T$ the number $x_t$
of coalition members of (sincere preference) type $t$, as well as the number $y_t$ of coalition
members who will insincerely vote $t$. Of course, we must have $\sum_{t\in T}x_t = \sum_{t\in T}y_t$.

Proposition \ignoreties\ suggests that it will be enough to consider manipulation
only at profiles (or voting situations) for which there is a clear winner $a$
(i.e. a candidate with $\abs{a}>\abs{\alpha}$ for each $\alpha\neq a$).
For such profiles, we will consider a manipulation attempt successful if it results in the score of another
candidate (or tied group of candidates) matching or exceeding the score of $a$.

For a coalition to successfully manipulate in favour of a candidate $\beta$, its members
must all be of types which prefer $\beta$ to $a$. Let $\Tpba\subseteq T$ be the set of such types.

Additionally, it is clear that the coalition members need only consider insincere
votes of types which rank $\beta$ in first place; this is a dominant strategy for
such manipulations. Let $\Tbeta$ be the set of such types.

Let $MCS$ be the minimum size of a successful manipulating coalition (or $\infty$ if no
manipulation is possible). Then in the light of Proposition \ignoreties\ and the above remarks
we have
$$ \Prob{MCS = \hbox{min}_{\beta\neq a} Q_1(\beta)} \to 1
  $$
where $Q_1(\beta)$ is the optimal value of the linear program
\def\LPQone{(1)}
$$ \eqalign{
  \hbox{min}\; & \sum_{t\in\Tpba} x_t \cr
  \hbox{s.t.}\; & \sum_{t\in\Tbeta} y_t (1-\sigma_t(\alpha)) - \sum_{t\in\Tpba} x_t (\sigma_t(\beta) - \sigma_t(\alpha))
          \geq |\alpha| - |\beta|  \qquad\forall\alpha\neq\beta \cr
  & \sum_{t\in\Tbeta} y_t = \sum_{t\in\Tpba} x_t , \cr
  & 0\leq x_t \leq N_t \qquad\forall t\in\Tpba \cr
  & y_t \geq0 \qquad\forall t\in\Tbeta \cr
  & x_t, y_t \in\Z 
           }
  \eqno{\LPQone}
  $$
(or $\infty$ if \LPQone\ is infeasible).

\subchunk{}
We now aim to show that, for IC asymptotic purposes, we may
drop the constraints $x_t \leq N_t$ and $x_t, y_t \in\Z$ in \LPQone. This is as we should
expect: each $N_t$ will be about $n/m!$, while the differences between
candidates' scores (and hence, presumably, coalition sizes) are only of order $\sqrt{n}$.
Similarly, the requirement that $x_t$ and
$y_t$ be integral should not be much of a hindrance when dealing with large numbers of voters.

\subchunk{}
To this end, consider $Q_2(\beta)$, defined in the same way as $Q_1(\beta)$, except that we drop the
integrality constraints $x_t, y_t \in\Z$ in \LPQone, and replace the constraint $x_t \leq N_t$
by $x_t \leq N_t - K$, where $K$ is a constant that depends only on the voting rule.
It will be convenient to choose
$$ K = \cases{ 2m!(1-w_{m-1})^{-1}, & if $w_{m-1}<1$ \cr 0, & if $w_{m-1}=1$. \cr}
  $$
Define $Q_3(\beta)$ in the same way as $Q_2(\beta)$, except that the upper bound constraint on $x_t$
is dropped entirely.
Note that we have $Q_3(\beta)\leq Q_1(\beta)$ and $Q_3(\beta)\leq Q_2(\beta)$.

\chunk{Proposition \ignoreNlimit.}
$$ \Prob{Q_2(\beta) = Q_3(\beta)} \to 1 \qquad\hbox{ as } n\to\infty .
  $$
\subchunk{Proof.}
Define $f_0:\R^T\to[0,\infty]$ as follows: for $q\in\R^T$, $f_0(q)$ is the optimal value of the
linear program
$$ \eqalign{
  \hbox{min}\; & \sum_{t\in\Tpba} x_t \cr
  \hbox{s.t.}\; & \sum_{t\in\Tbeta} y_t (1-\sigma_t(\alpha)) - \sum_{t\in\Tpba} x_t (\sigma_t(\beta) - \sigma_t(\alpha))
          \geq \sum_{t\in T} q_t (\sigma_t(\alpha) - \sigma_t(\beta))  \qquad\forall\alpha\neq\beta \cr
  & \sum_{t\in\Tbeta} y_t = \sum_{t\in\Tpba} x_t , \cr
  & x_t \geq0 \qquad\forall t\in\Tpba \cr
  & y_t \geq0 \qquad\forall t\in\Tbeta \cr
           }
  $$
Then $Q_3(\beta)=f_0(N)$. Note that
$f_0(\lambda q + \mu e) = \lambda f_0(q)$ for any $q\in\R^T$, $\lambda\geq0$, and $\mu\in\R$.
Also, $f_0$ is continuous on the closed subset $K=\theset{q}{f_0(q)<\infty}$ of $\R^T$. 
Let $\pi:\R^T\to K$ be the projection which maps each point $q\in\R^T$ to the nearest point
of $K$ to $q$; then $\pi$ is continuous and $\pi(\lambda q + \mu e) = \lambda\pi(q)+ \mu e$
for any $q\in\R^T$, $\lambda\geq0$, and $\mu\in\R$.

Let $f:\R^T\to [0,\infty)$ be given by $f(q) = f_0(\pi(q))$. Then $f$ is continuous; has
$f(\lambda q + \mu e) = \lambda f(q)$ for any $q\in\R^T$, $\lambda\geq0$, and $\mu\in\R$;
and $Q_3(\beta)=f(N)$ whenever $Q_3(\beta)<\infty$.

When $Q_3(\beta)=\infty$, we have $Q_2(\beta)=\infty$ too. On the other hand, when
$Q_3(\beta)\leq\hbox{min}_t N_t - K$, the corresponding optimal point of the linear program
for $Q_3(\beta)$ is also feasible for the linear program for $Q_2(\beta)$, so $Q_2(\beta)=Q_3(\beta)$.
Hence
$$ \set{Q_2(\beta)\neq Q_3(\beta)} \subseteq \set{\hbox{min}_t N_t - K < Q_3(\beta) < \infty}
   \subseteq \set{\hbox{min}_t N_t - K < f(N)} ,
  $$
and so it suffices to show that this last event has probability converging to 0.

Now define $h:\R^T\to(-\infty,\infty)$ by $h(q)=(\hbox{min}_t q_t) - f(q)$.
By Proposition \CLTN\ and the continuity of $h$, we have
$$ h\left({N - ne\over\sqrt{n}}\right) \cvgindistn h(X), \qquad\qquad\hbox{ where } X\sim N(0,\Sigma).
  $$
This yields
$$ { (\hbox{min}_t N_t) - f(N) - n/m! \over\sqrt{n}} \cvgindistn h(X) ,
  $$
from which
$$ \limsup_n \Prob{{(\hbox{min}_t N_t) - f(N) - n/m! \over\sqrt{n}} \leq -\lambda} \leq
       \Prob{h(X) \leq -\lambda }
  $$
for any $\lambda>0$. That is,
$$ \limsup_n \Prob{ (\hbox{min}_t N_t) - K - f(N) \leq n/m! -\lambda\sqrt{n} - K } \leq
       \Prob{h(X) \leq -\lambda }.
  $$
We have $n/m! -\lambda\sqrt{n} - K >0$ for sufficiently large $n$, so
$$ \limsup_n \Prob{ (\hbox{min}_t N_t) - K \leq f(N) } \leq \Prob{h(X) \leq -\lambda }.
  $$
Since $\lambda$ was arbitrary,
$$ \Prob{ (\hbox{min}_t N_t) - K \leq f(N) } \to 0 \qquad\qquad\hbox{ as }n\to\infty,
  $$
and the result follows.
\endofproof

\chunk{Proposition \ignoreintegrality.}
With probability 1,
$$ Q_1(\beta) \leq Q_2(\beta) + K .
  $$
\subchunk{Proof of Proposition \ignoreintegrality\ for the case $w_{m-1}<1$.}
Let $(x_t)_{t\in\Tpba}$, $(y_t)_{t\in\Tbeta}$ be optimal for the problem defining
$Q_2(\beta)$. Then
$$ \eqalign{
  & \sum_{t\in\Tbeta} y_t (1-\sigma_t(\alpha)) - \sum_{t\in\Tpba} x_t (\sigma_t(\beta) - \sigma_t(\alpha))
          \geq |\alpha| - |\beta|  \qquad\forall\alpha\neq\beta \cr
  & \sum_{t\in\Tbeta} y_t = \sum_{t\in\Tpba} x_t , \cr
  & 0\leq x_t \leq N_t - K \qquad\forall t\in\Tpba \cr
  & y_t \geq0 . \qquad\forall t\in\Tbeta \cr
           }
  $$
Choose types $t_0\in\Tpba$, $t_1\in\Tbeta$ such that $t_0$ ranks $a$ last, $\beta$ next-to-last,
and some $\gamma$ first, while $t_1$ is obtained from $t_0$ by transposing the rankings of
$\beta$ and $\gamma$.
Define $(x'_t)_{t\in\Tpba}$, $(y'_t)_{t\in\Tbeta}$ by $x'_{t_0}=x_{t_0}+K$, $y'_{t_1}=y_{t_1}+K$,
and $x'_t=x_t$ and $y'_t=y_t$ for all other $t$. Then $(x'_t)$, $(y'_t)$ satisfy
$$ \eqalign{
 \sum_{t\in\Tbeta} y'_t &(1-\sigma_t(\alpha)) - \sum_{t\in\Tpba} x'_t (\sigma_t(\beta) - \sigma_t(\alpha)) \cr
 &= \sum_{t\in\Tbeta} y_t (1-\sigma_t(\alpha)) - \sum_{t\in\Tpba} x_t (\sigma_t(\beta) - \sigma_t(\alpha)) 
    + K(1-w_{m-1}) + K(\sigma_{t_0}(\alpha) - \sigma_{t_1}(\alpha)) \cr
 &\geq |\alpha| - |\beta| + 2m! \cr
           }
  $$
for each $\alpha\neq\beta$. Also, $ \sum_{t\in\Tbeta} y'_t = \sum_{t\in\Tpba} x'_t$
and $0\leq x'_t \leq N_t$ for each $t$.

Now obtain $(x^{''}_t)_{t\in\Tpba}$, $(y^{''}_t)_{t\in\Tbeta}$ by rounding each $x'_t$, $y'_t$
to an integral value. The choice between rounding up and rounding down (i.e. between
$x^{''}_t=\grint{x'_t}$ and $x^{''}_t=\lceil x'_t\rceil$) can be made arbitrarily, but should
be done in such a way that $\sum_{t\in\Tbeta} y^{''}_t = \sum_{t\in\Tpba} x^{''}_t$.
Since $\abs{x^{''}_t-x'_t}\leq1$ and $\abs{y^{''}_t-y'_t}\leq1$, we obtain
$$ \sum_{t\in\Tbeta} y^{''}_t (1-\sigma_t(\alpha)) - \sum_{t\in\Tpba} x^{''}_t (\sigma_t(\beta) - \sigma_t(\alpha))
   \geq |\alpha| - |\beta| ,
  $$
and so $(x^{''}_t)$, $(y^{''}_t)$ are feasible for \LPQone.
We have $\sum_{t\in\Tpba} x^{''}_t = \sum_{t\in\Tpba} x_t + K$, from which it follows that
$Q_1(\beta) \leq Q_2(\beta) + K$.
\endofproof

\subchunk{Proof of Proposition \ignoreintegrality\ for the case $w_{m-1}=1$.}
A separate proof is required for this case (the anti-plurality rule $w=(1,\ldots,1,0)$).
We can show that $Q_1(\beta)\leq Q_2(\beta)$ by showing that the optimal $(x_t)$, $(y_t)$
for the problem defining $Q_2(\beta)$ are always integral (and hence give a feasible
solution to \LPQone). To establish this, we will use a well-known result in linear
programming (see, e.g. \cite{\PapadimitriouSteiglitz} or \cite{\ParkerRardin}), which
assures us that the optimal solution
of a linear program will always be integral when the constraint coefficient matrix $A$
is totally unimodular (i.e. every square submatrix has determinant $1$, $-1$, or $0$).

A useful sufficient condition for total unimodularity is given in \cite{\PapadimitriouSteiglitz}
(Theorem 13.3) as follows: a matrix whose entries are all $1$, $-1$, or $0$ is
totally unimodular if each column has at most two non-zero entries, and if the rows
can be partitioned into two sets $I_1$ and $I_2$ such that: (i) if a column has two entries
of the same sign, their rows are in different sets; (ii) if a column has two entries
of different signs, their rows are in the same set.

For this problem $A$ has columns corresponding to the variables $x_t$ ($t\in\Tpba$) and
$y_t$ ($t\in\Tbeta$). There is one row corresponding to each candidate $\alpha\neq\beta$,
in which the entry in the column corresponding to $x_t$ is $-1$ if $t$ ranks $\alpha$ last,
and $0$ otherwise; the entry in the column corresponding to $y_t$ is $1$ if $t$ ranks
$\alpha$ last, and $0$ otherwise. Let these rows constitute the set $I_1$.
There is also a further row corresponding to the constraint
$\sum_{t\in\Tbeta} y_t - \sum_{t\in\Tpba} x_t =0$; in this row, the entries in the
columns corresponding to the $x_t$ are all $-1$, and those corresponding to the $y_t$
are all $1$. Let this row constitute the set $I_2$. It is clear that this matrix
satisfies the sufficient condition above; the result follows.
\endofproof

\chunk{Corollary \ignoreall.}
$$ \Prob{\abs{MCS - \hbox{min}_{\beta\neq a} Q_3(\beta)} \leq K} \to 1 \qquad\qquad\hbox{as $n\to\infty$}
  $$

\subchunk{Proof.}
From Propositions \ignoreNlimit\ and \ignoreintegrality\ we have
$$ \Prob{Q_1(\beta)-K \leq Q_2(\beta) \leq Q_3(\beta) \leq Q_1(\beta)} \to 1 \qquad\qquad\hbox{as $n\to\infty$} ,
  $$
from which
$$ \Prob{\abs{Q_1(\beta)-Q_3(\beta)}\leq K} \to1 \qquad\qquad\hbox{as $n\to\infty$,}
  $$
and so
$$ \Prob{\abs{\hbox{min}_{\beta\neq a} Q_1(\beta) - \hbox{min}_{\beta\neq a} Q_3(\beta)}\leq K} \to1
     \qquad\qquad\hbox{as $n\to\infty$.}
  $$
The result then follows from the earlier observation that 
$\Prob{MCS = \hbox{min}_{\beta\neq a} Q_1(\beta)} \to 1$ .
\endofproof

\subchunk{}
Now that we have effectively eliminated, for our purposes, the integrality and upper-bound
constraints in \LPQone, we can further simplify the description of the minimum manipulating
coalition size.

Suppose $b$ is the candidate with second-highest score after $a$. The next result consists
of the observation that only manipulations in favour of $b$ need now be considered.

\chunk{Proposition \manipbonly.}
$\hbox{min}_{\beta\neq a} Q_3(\beta) = Q_3(b)$.

\subchunk{Proof.}
Let $(x_t)_{t\in\Tpba}$, $(y_t)_{t\in\Tbeta}$ be optimal for the problem defining $Q_3(\beta)$,
some $\beta\neq b$. Form $(x'_t)_{t\in\Tpba}$, $(y'_t)_{t\in\Tbeta}$ by transposing $\beta$
and $b$ in all the voter types involved. (So if types $s$ and $t$ are related by
transposition of the ranks of $\beta$ and $b$, then $x'_s=x_t$.) We have, for any $\alpha$,
$$ \sum_{t\in\Tb} y'_t (1-\sigma_t(\alpha)) - \sum_{t\in\Tpref_{ba}} x'_t (\sigma_t(b) - \sigma_t(\alpha))
  = \sum_{t\in\Tbeta} y_t (1-\sigma_t(\alpha)) - \sum_{t\in\Tpba} x_t (\sigma_t(\beta) - \sigma_t(\alpha)) .
  $$
Since $\abs{\beta}\leq\abs{b}$, we see that $(x'_t)$, $(y'_t)$ are feasible for the problem
defining $Q_3(b)$, and give the same objective value. The result follows.
\endofproof

\subchunk{}
We next show that a manipulating coalition may always be formed by recruiting
only voters who sincerely rank $b$ and $a$ adjacent.
Let $T_i\subset T$ consist of those types which rank $b$ in $i$th place and $a$ in $(i+1)$st place,
and $\Tba=\cup_{i=1}^{m-1} T_i\subset \Tpref_{ba}$.
If we replace $\Tpref_{ba}$ by $\Tba$ in the linear program defining $Q_3(b)$,
we obtain the linear program

\def\LPa{(2)}
$$ \eqalign{
  \hbox{min}\; & \sum_{t\in\Tba} x_t \cr
  \hbox{s.t.}\; & \sum_{t\in\Tb} y_t (1-\sigma_t(\alpha)) - \sum_{t\in\Tba} x_t (\sigma_t(b) - \sigma_t(\alpha))
          \geq |\alpha| - |b|  \qquad\forall\alpha\neq b \cr
  & \sum_{t\in\Tb} y_t = \sum_{t\in\Tba} x_t , \cr
  & x_t \geq0 \qquad\forall t\in\Tba \cr
  & y_t \geq0 \qquad\forall t\in\Tb \cr
           }
  \eqno{\LPa}
  $$
Let $Q$ denote the optimal value of \LPa\ (or $\infty$ if \LPa\ is infeasible).

\chunk{Proposition \Tionly.}
$Q_3(b)=Q$.

\subchunk{Proof.}
For $t\in T_i$, let $\tau(t)\subseteq\Tpref_{ba}$ consist of those types
which (i) agree with $t$ in ranking positions $i+1,\ldots,m$; and (ii)
rank the remaining candidates, other than $b$, in either the same position
as $t$ does, or one place lower. Observe that $\theset{\tau(t)}{t\in\Tba}$
is a partition of $\Tpref_{ba}$, and that for $s\in\tau(t)$, $\alpha\neq b$
$$ \sigma_t(b) - \sigma_t(\alpha) \leq \sigma_s(b) - \sigma_s(\alpha) .
  $$
(This implies that a voter of type $s\in\tau(t)$ can always be dismissed from
the manipulating coalition and replaced with one of type $t$.)

Let $(x_t)_{t\in\Tpref_{ba}}$, $(y_t)_{t\in\Tb}$ be optimal for the problem defining $Q_3(b)$.
Form $(x'_t)_{t\in\Tba}$ by $x'_t = \sum_{s\in\tau(t)} x_s$. Then
$\sum_{t\in\Tba} x'_t = \sum_{s\in\Tpref_{ba}} x_s$, and
$$ \sum_{t\in\Tba} x'_t (\sigma_t(b) - \sigma_t(\alpha) \leq
     \sum_{s\in\Tpref_{ba}} x_s (\sigma_s(b) - \sigma_s(\alpha)) .
  $$
Hence $(x'_t)$ is feasible for \LPa. The result follows.
\endofproof

\subchunk{}
Our efforts thus far have established that $\Prob{\abs{MCS - Q}\leq K}\to 1$ as $n\to\infty$.
This will ensure that for IC asymptotic purposes, we may replace our original
description of the minimum manipulating coalition size with the more tractable
problem \LPa. Our next step results in a considerable further simplification of
the linear program.

\chunk{Theorem \zsufficient.} Suppose non-negative numbers $z_1,\ldots,z_{m-1}$ satisfy
\def\zconds{(3)}
$$ \eqalign{
  \sum_{i=1}^{m-1} z_i (1-w_i+w_{i+1}) &\geq |a| - |b| \cr
  \sum_{i=1}^{m-1} z_i (1-w_i) &\geq n\bar w - |b| \cr
           }
  \eqno{\zconds}
  $$
where $\bar w = (w_1+\cdots+w_m)/m$.
Then there exists $(x_t)_{t\in\Tba}$, $(y_t)_{t\in\Tb}$ feasible for \LPa, with
$\sum_{t\in T_i} x_t = z_i$.

\subchunk{Remark.} 
The conclusion of this result asserts that we may arrange a successful manipulation in which
$z_i$ members of the coalition are of types in $T_i$, for each $i=1,\ldots,m-1$.

Note that if $z_i$ voters of types in $T_i$ all cast insincere votes which rank $b$ first,
the score of $b$ will be increased by $z_i(1-w_i)$. If they all cast insincere votes which rank
$a$ last, the score of $a$ will be decreased by $z_i w_{i+1}$. Thus, the first inequality
of \zconds\ makes it possible for $b$ to catch up to $a$. The second inequality of \zconds\ makes
it possible for $b$ to catch up to the average candidate's score.
It is remarkable that these two (apparently rather weak) linear conditions are sufficient
to make manipulation possible.

\subchunk{Proof.}
Let $D$ be the set of candidates other than $a$ or $b$.
For any type $t\in T$, we denote by $t(i)$ the candidate ranked in $i$th place by $t$.

We will write our proposed $(x_t)$, $(y_t)$ in terms of parameters $r$ and $(u_\alpha)_{\alpha\in D}$,
which are to be determined later in such a way that $0\leq r\leq 1$, $u_\alpha\geq0$ $\forall\alpha$,
and $\sum_{\alpha\in D} u_{\alpha} = 1$.

For each $\alpha\in D$, let
$$ v_{\alpha} = \cases{ {1-u_{\alpha}\over m-3} &, if $m\geq4$ \cr
                         1 &, if $m=3$; \cr
                      }
  $$
note $\sum_{\alpha\in D} v_{\alpha} = 1$.

For each $t$, let $x_t = \sum_{k=1}^4 x_t^{(k)}$ and $y_t = \sum_{k=1}^4 t_t^{(k)}$, where
$$ \eqalign{
   x_t^{(1)} &= \cases{ {ru_{t(m)} z_i \over (m-3)!} &, if $t\in T_i$, $i=1,\ldots,m-2$ \cr
                        0 &, if $t\in T_{m-1}$ \cr
                      } \cr
\noalign{\smallskip}
   y_t^{(1)} &= \cases{ \sum_{i=1}^{m-2} {ru_{t(i+1)} z_i \over (m-3)!} &, if $t\in\Tb$ with $t(m)=a$ \cr
                        0 &, otherwise \cr
                      } \cr
           }
  $$
(this corresponds to $ru_{\alpha}z_i$ voters of types in $T_i$ changing sincere votes of form
$\ldots ba\ldots\alpha$ to insincere ones $b\ldots\alpha\ldots a$, for each $i=1,\ldots,m-2$);
$$ \eqalign{
   x_t^{(2)} &= \cases{ {(1-r)v_{t(1)} z_i \over (m-3)!} &, if $t\in T_i$, $i=2,\ldots,m-2$ \cr
                        0 &, if $t\in T_1\union T_{m-1}$ \cr
                      } \cr
\noalign{\smallskip}
   y_t^{(2)} &= \cases{ {(1-r)v_{t(i)} z_i \over (m-3)!} &, if $t\in\Tb$ with $t(i+1)=a$, $i=2,\ldots,m-2$ \cr
                        0 &, otherwise \cr
                      } \cr
           }
  $$
(this corresponds to $(1-r)v_{\alpha}z_i$ voters of types in $T_i$ changing sincere votes of form
$\alpha\ldots ba\ldots$ to insincere ones $b\ldots\alpha a\ldots$, for each $i=2,\ldots,m-2$);
$$
   x_t^{(3)} = y_t^{(3)} = \cases{ {(1-r)z_1\over (m-2)!} &, if $t\in T_1$ \cr
                                   0 &, otherwise; \cr
                                 }
  $$
(this corresponds to $(1-r)z_1$ voters of types in $T_1$ leaving their votes unchanged);
$$ \eqalign{
   x_t^{(4)} &= \cases{ {v_{t(1)} z_{m-1} \over (m-3)!} &, if $t\in T_{m-1}$ \cr
                        0 &, if $t\notin T_{m-1}$ \cr
                      } \cr
\noalign{\smallskip}
   y_t^{(4)} &= \cases{ {v_{t(m-1)} z_{m-1} \over (m-3)!} &, if $t\in\Tb$ with $t(m)=a$ \cr
                        0 &, otherwise \cr
                      } \cr
           }
  $$
(this corresponds to $v_{\alpha}z_{m-1}$ voters of types in $T_{m-1}$ changing sincere votes of form
$\alpha\ldots ba$ to insincere ones $b\ldots\alpha a$).

The reader may verify that
$$ \sum_{t\in T_i} x_t = \sum_{k=1}^4 \sum_{t\in T_i} x_t^{(k)} = z_i
   \qquad\hbox{ and }\qquad
   \sum_{t\in\Tb} y_t = \sum_{k=1}^4 \sum_{t\in\Tb} y_t^{(k)} = \sum_{i=1}^{m-1} z_i .
  $$
We now verify that the inequality constraints of \LPa\ hold.
Note that
$$ \sum_{t\in\Tba} x_t \sigma_t(b) = \sum_{i=1}^{m-1} z_i w_i ,
   \qquad\qquad\qquad
   \sum_{t\in\Tba} x_t \sigma_t(a) = \sum_{i=1}^{m-1} z_i w_{i+1} ,
  $$
and
$$ \eqalign{
 \sum_{t\in\Tb} y_t \sigma_t(a) &= \sum_{k=1}^4 \sum_{t\in\Tb} y_t^{(k)} \sigma_t(a) \cr
   &= 0 \;+\; (1-r)\sum_{i=2}^{m-2} z_i w_{i+1} \;+\; (1-r)z_1 w_2 \;+\; 0 \cr
   &= (1-r)\sum_{i=1}^{m-1} z_i w_{i+1} . \cr
           }
  $$
Hence
$$ \sum_{t\in\Tb} y_t (1-\sigma_t(a)) - \sum_{t\in\Tba} x_t (\sigma_t(b)-\sigma_t(a)) = B + rA ,
  $$
where $A=\sum_{i=1}^{m-1} z_i w_{i+1}$ and $B=\sum_{i=1}^{m-1} z_i (1-w_i)$.
The inequality constraint for $\alpha=a$ in \LPa\ thus reduces to
$$ B + rA \geq \abs{a} - \abs{b} .
  $$
Now consider the other inequality constraints. We use the notation
$\bar w_{-i,j,k}$ to denote the average of all elements of $w$ other than the $i$th,
$j$th, and $k$th (i.e. $(-w_i-w_j-w_k+\sum_{\ell=1}^m w_{\ell})/(m-3)$), or
$0$ if $m=3$. Similarly $\bar w_{-i,j}=(-w_i-w_j+\sum_{\ell=1}^m w_{\ell})/(m-2)$.
For $\alpha\in D$ we have
$$ \eqalign{
 \sum_{t\in\Tba} x_t \sigma_t(\alpha) &= \sum_{k=1}^4 \sum_{t\in\Tba} x_t^{(k)} \sigma_t(\alpha) \cr
 = &\quad r\sum_{i=1}^{m-2} z_i \left(u_{\alpha}\cdot 0 + 
          \left(\sum_{\gamma\neq\alpha} u_{\gamma}\right)\bar w_{-i,i+1,m}\right) \cr
 &+ (1-r)\sum_{i=2}^{m-2} z_i \left(v_{\alpha}\cdot 1 + 
          \left(\sum_{\gamma\neq\alpha} v_{\gamma}\right)\bar w_{-1,i,i+1}\right) \cr
 &+ (1-r)z_1 \bar w_{-1,2} \cr
 &+ z_{m-1} \left(v_{\alpha}\cdot 1 +
          \left(\sum_{\gamma\neq\alpha} v_{\gamma}\right)\bar w_{-1,m-1,m}\right) \cr
           }
  $$
and
$$ \eqalign{
 \sum_{t\in\Tb} y_t \sigma_t(\alpha) &= \sum_{k=1}^4 \sum_{t\in\Tb} y_t^{(k)} \sigma_t(\alpha) \cr
 = &\quad r\sum_{i=1}^{m-2} z_i \left(u_{\alpha}w_{i+1} +
          \left(\sum_{\gamma\neq\alpha} u_{\gamma}\right)\bar w_{-1,i+1,m}\right) \cr
 &+ (1-r)\sum_{i=2}^{m-2} z_i \left(v_{\alpha}w_i +
          \left(\sum_{\gamma\neq\alpha} v_{\gamma}\right)\bar w_{-1,i,i+1}\right) \cr
 &+ (1-r)z_1 \bar w_{-1,2} \cr
 &+ z_{m-1} \left(v_{\alpha}w_{m-1} +
          \left(\sum_{\gamma\neq\alpha} v_{\gamma}\right)\bar w_{-1,m-1,m}\right) . \cr
           }
  $$
It follows that
$$ \eqalign{
 \sum_{t\in\Tb} y_t (1 - \sigma_t(\alpha)) &- \sum_{t\in\Tba} x_t (\sigma_t(b) - \sigma_t(\alpha)) = \cr
 &\quad \sum_{i=1}^{m-1} z_i - ru_{\alpha} \sum_{i=1}^{m-2} z_i w_{i+1}
     + r(1-u_{\alpha})\sum_{i=1}^{m-2} z_i (\bar w_{-i,i+1,m} - \bar w_{-1,i+1,m}) \cr
 &\qquad\qquad\qquad + (1-r)v_{\alpha} \sum_{i=2}^{m-2} z_i (1-w_i) + z_{m-1} v_{\alpha} (1-w_{m-1})
     - \sum_{i=1}^{m-1} z_i w_i \cr
 &= \sum_{i=1}^{m-1} z_i (1-w_i) - ru_{\alpha} \sum_{i=1}^{m-1} z_i w_{i+1}
     + r v_{\alpha} \sum_{i=1}^{m-2} z_i (1-w_i) \cr
 &\qquad\qquad\qquad + (1-r)v_{\alpha} \sum_{i=1}^{m-2} z_i (1-w_i) + z_{m-1} v_{\alpha} (1-w_{m-1}) \cr
 &= (1+v_{\alpha})B - r u_{\alpha} A .\cr
           }
  $$
The $(x_t),(y_t)$ that we have constructed will thus satisfy the constraints required to be
feasible for \LPa, provided $r$ and $(u_\alpha)$ can be chosen in such a way that
$$ \eqalign{
 B + rA &\geq \abs{a}-\abs{b} \cr
 (1+v_{\alpha})B - ru_{\alpha}A &\geq \abs{\alpha} - \abs{b} \qquad\qquad\forall\alpha\in D \cr
 \sum_{\alpha\in D} u_{\alpha} &= 1 \cr
 u_{\alpha} &\geq 0 \qquad\qquad\forall\alpha\in D \cr
 0\leq r &\leq 1 . \cr
           }
  $$

Suppose for the moment that $m\geq4$. Then the inequality required of $u_\alpha$ may be
written
$$ \left(1 + {1-u_{\alpha}\over m-3}\right)B - ru_{\alpha}A \geq \abs{\alpha} - \abs{b} ,
  $$
or
$$ \left(rA + {B\over m-3}\right) u_{\alpha} \leq \abs{b} - \abs{\alpha} + \left({m-2\over m-3}\right)B ,
  $$
a simple upper bound on $u_{\alpha}$. Note that the upper bound is non-negative.
For given $r$, it will be possible to choose non-negative $(u_\alpha)$ such that
$\sum_{\alpha\in D} u_{\alpha}=1$ while complying with these upper bounds if and only if
$$ \sum_{\alpha\in D} \left(\abs{b} - \abs{\alpha} + \left({m-2\over m-3}\right)B\right)
   \geq rA + {B\over m-3} ,
  $$
or
$$ rA \leq (m-2)\abs{b} - \sum_{\alpha\in D}\abs{\alpha} + \left({(m-2)^2-1\over m-3}\right)B ,
  $$
that is
$$ rA \leq (m-1)(\abs{b} + B - n\bar w) + \abs{a} - n\bar w ,
   \eqno{(*)}
  $$
using the fact that the sum of all the candidates' scores is $mn\bar w$.

Consider now the case $m=3$. Then the sole $\alpha\in D$ has $u_{\alpha}=v_{\alpha}=1$, and the
inequality required of $u_\alpha$ reduces to
$$ 2B - rA \geq \abs{\alpha} - \abs{b} ,
  $$
which is $(*)$.

We have thus reduced our requirements to a condition on $r\in[0,1]$:
$$ \abs{a}-\abs{b}-B \;\leq\;  rA \;\leq\; (m-1)(\abs{b} + B - n\bar w) + (\abs{a} - n\bar w) .
  $$
To see that there exists $r\in[0,1]$ satisfying this condition we note that, firstly,
$$ (m-1)(\abs{b} + B - n\bar w) + (\abs{a} - n\bar w) \geq 0 ;
  $$
secondly,
$$ \abs{a}-\abs{b}-B \leq A
  $$
(from the first condition of \zconds); and thirdly,
$$ \abs{a}-\abs{b}-B \leq (m-1)(\abs{b} + B - n\bar w) + (\abs{a} - n\bar w) .
  $$
This last condition can be simplified to
$$ m(\abs{b} + B - n\bar w) \geq 0 ,
  $$
the second condition of \zconds.
\endofproof

\medskip
\chunk{Corollary \zprog.} The linear program \LPa\ has the same optimal value as the following one.
\def\LPb{(4)}
$$ \eqalign{
  \hbox{min}\; & \sum_{i=1}^{m-1} z_i \cr
  \hbox{s.t.}\; & \sum_{i=1}^{m-1} (1-w_i+w_{i+1})z_i \geq |a| - |b| \cr
  & \sum_{i=1}^{m-1} (1-w_i)z_i \geq n\bar w - |b| \cr
  & z_i \geq0 \qquad\hbox{ for } i=1,\ldots,m-1 \cr
           }
  \eqno{\LPb}
  $$
\medskip
{\it Proof.}
Theorem \zsufficient\ shows that for any feasible point of \LPb, there corresponds a feasible point
of \LPa\ with the same objective value. Hence the optimal value of \LPa\ cannot be greater than that of \LPb.

To show that the optimal value of \LPa\ cannot be less than that of \LPb, suppose we have
$(x_t)_{t\in\Tba}$, $(y_t)_{t\in\Tb}$ feasible for \LPa. Let $z_i=\sum_{t\in T_i} x_t$.
Then $\sum_{i=1}^{m-1} z_i = \sum_{t\in\Tba} x_t$, and we will show that $(z_i)_{i=1}^{m-1}$
is feasible for \LPb. The inequality for $\alpha=a$ in \LPa\ says that
$$ \sum_{t\in\Tb} y_t (1-\sigma_t(a)) \geq |a| - |b| + \sum_{t\in\Tba} x_t (\sigma_t(b) - \sigma_t(a)) .
  $$
Noting that $\sum_{t\in\Tb} y_t \sigma_t(a) \geq0$, we obtain
$$ \sum_{t\in\Tb} y_t \geq |a| - |b| + \sum_{i=1}^{m-1} \sum_{t\in T_i} x_t (w_i - w_{i+1})
  $$
and so
$$ \sum_{i=1}^{m-1} z_i \geq |a| - |b| + \sum_{i=1}^{m-1} (w_i - w_{i+1}) z_i
  $$
from which the first constraint of \LPb\ can be seen to hold.

If we add the inequalities in \LPa\ for all $\alpha\neq b$, we obtain
$$ \sum_{t\in\Tb} y_t \sum_{\alpha\neq b} (1-\sigma_t(\alpha)) \geq \left(\sum_{\alpha\neq b} |\alpha|\right)
   - (m-1)|b| + \sum_{t\in\Tba} x_t \sum_{\alpha\neq b} (\sigma_t(b) - \sigma_t(\alpha)) ,
  $$
or
$$ \sum_{t\in\Tb} y_t \sum_{i=1}^{m-1} (1-w_i) \geq nm\bar w - m|b| + 
   \sum_{i=1}^{m-1} \sum_{t\in T_i} x_t \sum_{j\neq i} (w_i - w_j) ,
  $$
which gives
$$ (1-\bar w) \sum_{i=1}^{m-1} z_i \geq n\bar w - |b| + \sum_{i=1}^{m-1} z_i (w_i - \bar w) ,
  $$
from which the second constraint of \LPb\ can be seen to hold.
\endofproof

\chunk{Theorem \dualprog.}
The linear programs \LPa\ and \LPb\ have the same objective value as the following one.
\def\LPc{(5)}
$$ \eqalign{
  \hbox{max}\; & (\abs{a}-n\bar w)\lambda + (n\bar w - \abs{b})\mu  \cr
  \hbox{s.t.}\; & w_{i+1} \lambda + (1-w_i)\mu \leq 1 \qquad\hbox{ for } i=1,\ldots,m-1 \cr
  & 0 \leq \lambda \leq \mu . \cr
           }
  \eqno{\LPc}
  $$
Also, \LPc\ is unbounded if and only if \LPa\ and \LPb\ are infeasible.

\subchunk{Proof.}
The linear program \LPb\ has the same objective value as its dual program (\cite{\Bazaraa, p. 251}):
$$ \eqalign{
  \hbox{max}\; & (\abs{a}-\abs{b})\lambda + (n\bar w - \abs{b})\lambda' \cr
  \hbox{s.t.}\; & (1-w_i+w_{i+1})\lambda + (1-w_i)\lambda' \leq 1 \qquad\hbox{ for } i=1,\ldots,m-1 \cr
  & \lambda\geq0 , \quad \lambda'\geq0 . \cr
           }
  $$
Substituting $\mu=\lambda+\lambda'$ yields \LPc.
\endofproof

\subchunk{Remark.}
We have now replaced our original description of the minimum manipulating coalition size
with the simple two-variable linear program \LPc. This will be exploited in the remaining
sections of the paper.

\section{4. IC asymptotics of the minimum manipulating coalition size}

The results of the previous section have established that
$$ \Prob{\abs{MCS - Q}\leq K} \to 1 \qquad\qquad\hbox{as $n\to\infty$,}
  $$
where
$$ Q = \hbox{max}\theset{\lambda(\abs{a}-n\bar w)+\mu(n\bar w-\abs{b})}{(\lambda,\mu)\in M_w}
  $$
and
$$ M_w = \theset{(\lambda,\mu)}{0\leq\lambda\leq\mu\hbox{ and }w_{i+1}\lambda + (1-w_i)\mu\leq1\hbox{ for }i=1,\ldots,m-1} .
  $$
Note in particular that the constraint set $M_w$ does not depend on the
voting situation, but only on the voting rule. For a given rule, it is possible to identify the
corresponding $M_w$ (a two-dimensional linear polytope), and then to identify
the (finitely many) vertices of $M_w$ which may achieve the optimum $Q$. This often leads
to an explicit expression for $Q$ in terms of $\abs{a}$ and $\abs{b}$.

Furthermore, we have $(\abs{a},\abs{b})=(\rho_1(S),\rho_2(S))$, where
$\rho_j(x)$ denotes the $j$th largest element of a vector $x$.
Proposition \CLTS\ gives us
$$ {(\rho_1(S) - n\bar w, n\bar w - \rho_2(S))\over\sqrt{n}} \cvgindistn
   \sigma_w \left({m\over m-1}\right)^{1/2}(\rho_1(Z)-\bar Z, \bar Z - \rho_2(Z)) ,
  $$
and so it follows that
$$ {Q\over\sqrt{n}} \cvgindistn V_w ,
  $$
where
$$ V_w =
  \hbox{max}\theset{\lambda(\rho_1(Z)-\bar Z)+\mu(\bar Z-\rho_2(Z))}{(\lambda,\mu)\in \sigma_w \left({m\over m-1}\right)^{1/2} M_w} .
  $$
Hence, too,
$$ {MCS\over\sqrt{n}} \cvgindistn V_w ,
  $$
by the converging-together lemma (\cite{\Durrett, p.91}).
Consequently,
$$ \Prob{MCS\leq v\sqrt{n}} \to g_w(v) := \Prob{V_w\leq v} \qquad\hbox{ as }n \to\infty.
  $$
That is, the asymptotic probability that the voting situation is manipulable by a
coalition of $v\sqrt{n}$ or fewer voters is computable as a (non-decreasing) function
of $v$. This function depends only on the voting rule.

\subchunk{}
A further observation will be helpful in determining which vertices of $M_w$ may
achieve the above maximum. If $x\in\R^m$ has mean element $\bar x = (x_1+\cdots+x_m)/m$, then
$$ \rho_1(x) - \bar x \geq 0 
   \qquad\qquad\hbox{and}\qquad\qquad
   -(\rho_1(x) - \bar x) \;\leq\; \bar x - \rho_2(x) \;\leq\; {1\over m-1}(\rho_1(x) - \bar x) .
  $$

We note in passing that for all rules other than the anti-plurality rule $w=(1,\ldots,1,0)$,
$M_w$ is a bounded set, since the constraints defining it include $(1-w_{m-1})\mu\leq1$.
Hence $V_w$ is a finite-valued random variable. It follows that $\Prob{MCS=\infty}\to0$
as $n\to\infty$. For the anti-plurality rule itself,
$M_w$ is unbounded and $V_w$ may take $\infty$ as a value; see below for further detail.

\subchunk{}
The remainder of this section will be devoted to carrying out the above analysis for 
some common positional voting rules.

\subchunk{Borda's rule.} $w_i=(m-i)/(m-1)$ for $i=1,\ldots,m$.
The constraints defining $M_w$ for this rule are
$$ \eqalign{
 (m-2)\lambda &\leq m-1 \cr
 (m-3)\lambda + \mu &\leq m-1 \cr
 \vdots \quad & \cr
 \lambda + (m-3)\mu &\leq m-1 \cr
 (m-2)\mu &\leq m-1 \cr
 0\leq\lambda&\leq\mu \cr
           }
  $$
Note that the point with $\lambda=\mu=(m-1)/(m-2)$ satisfies each of the constraints
with equality. So $M_w$ is the triangle with vertices at $(0,0)$,
$(0, {m-1\over m-2})$, and $({m-1\over m-2},{m-1\over m-2})$; the last of these always achieves the optimum $Q$.
Thus, a good approximation to the minimum manipulating coalition size is
$$ Q = \left({m-1\over m-2}\right) (\abs{a}-\abs{b}) .
  $$
For the corresponding asymptotic result, note that
$$ \sigma_w^2 = {m+1 \over 12(m-1)} .
  $$
It then transpires that
$$ {MCS\over\sqrt{n}} \cvgindistn \left({m(m+1)\over 12(m-2)^2}\right)^{1/2}(\rho_1(Z)-\rho_2(Z)) .
  $$

\subchunk{Anti-plurality rule:} $w=(1,\ldots,1,0)$.
The constraints defining $M_w$ reduce to $0\leq\lambda\leq1$, $\mu\geq\lambda$.
Thus
$$ Q = \cases{ \abs{a}-\abs{b}, & if $\abs{b}\geq n\bar w$ \cr
             \infty, & otherwise. \cr
             }
  $$
The corresponding asymptotic result is
$$ {MCS\over\sqrt{n}} \cvgindistn V_w \;=\;
   \cases{m^{-1/2} (\rho_1(Z)-\rho_2(Z)), & if $\rho_2(Z)\geq\bar Z$ \cr
          \infty, & otherwise. \cr
         }
  $$
As a consequence of this, we can find the limiting probability that an
anti-plurality election is invulnerable to manipulation:
$$ \lim_n \Prob{MCS=\infty} = \Prob{V_w=\infty} = \Prob{\rho_2(Z) < \bar Z} .
  $$
The limit is positive, a property unique to anti-plurality. For all other
positional rules, $\lim_n \Prob{MCS=\infty} = 0$, a result essentially
contained in \cite{\KimRoush}.

\subchunk{Plurality and $k$-approval rules:} $w=(1,\ldots,1,0,\ldots0)$ (with $k$ 1s), where $1\leq k\leq m-2$.
The simple plurality rule is included here as the case $k=1$.
The constraints defining $M_w$ reduce to $0\leq\lambda\leq\mu\leq1$, so $M_w$
is the triangle with vertices at $(0,0)$, $(0,1)$, and $(1,1)$; the last of
these always achieves the optimum. Thus, our approximation of the
minimum manipulating coalition size is simply
$$ Q =\abs{a}-\abs{b} .
  $$
While this expression is valid for all $k$-approval rules ($1\leq k\leq m-2$),
different values of $k$ will give rise to different probability distributions for
$(\abs{a},\abs{b})$, and so different asymptotic results for $MCS$.
We have $\bar w=k/m$ and $\sigma_w^2=k(m-k)/m^2$, giving
$$ {MCS\over\sqrt{n}} \cvgindistn \left({k(m-k)\over m(m-1)}\right)^{1/2}(\rho_1(Z)-\rho_2(Z)) .
  $$

\subchunk{Three-candidate ``easy case" rules:} $w=(1,1-p,0)$ where $1/2\leq p\leq 1$.
This family includes the 3-candidate versions of plurality voting ($p=1$) and Borda's rule ($p=1/2$).
The constraints defining $M_w$ are $0\leq\lambda\leq\mu$, $(1-p)\lambda\leq1$, and
$p\mu\leq1$. Thus $M_w$ is the triangle with vertices at $(0,0)$,
$(0,p^{-1})$, and $(p^{-1},p^{-1})$; the last of these always achieves the optimum.
Thus
$$ Q=p^{-1}(\abs{a}-\abs{b}) .
  $$
We have $\sigma_w^2=2(1-p+p^2)/9$; it follows from this that
$$ {MCS\over\sqrt{n}} \cvgindistn \left({1-p+p^2\over 3p^2}\right)^{1/2}(\rho_1(Z)-\rho_2(Z)) .
  $$

\subchunk{Three-candidate ``hard case" rules:} $w=(1,1-p,0)$ where $0<p\leq 1/2$.
This family includes the remaining 3-candidate positional rules not already considered.
For these rules, $M_w$ is the quadrilateral with vertices at $(0,0)$,
$(0,p^{-1})$, $((1-p)^{-1},(1-p)^{-1})$, and $(1-p)^{-1},p^{-1})$.
In voting situations with $\abs{b}\geq n\bar w$, the optimum for $Q$ is achieved at
$((1-p)^{-1},(1-p)^{-1})$; otherwise, $(1-p)^{-1},p^{-1})$ is optimal. Thus
$$ Q=\cases{ (1-p)^{-1}(\abs{a}-\abs{b}), & if $\abs{b}\geq n\bar w$ \cr
      (1-p)^{-1}(\abs{a}-n\bar w) + p^{-1}(n\bar w-\abs{b}), & if $\abs{b}\leq n\bar w$ \cr
         }
  $$
that is
$$ Q=(1-p)^{-1}(\abs{a}-\abs{b}) + \left({1\over p} - {1\over 1-p}\right)(n\bar w - \abs{b})_+ ,
  $$
where $x_+$ denotes max$(x,0)$. The corresponding asymptotic result is
$$ {MCS\over\sqrt{n}} \cvgindistn
  \left({1-p+p^2\over 3}\right)^{1/2} \left(\left({1\over1-p}\right)(\rho_1(Z)-\rho_2(Z))
    + \left({1\over p} - {1\over 1-p}\right)(\bar Z - \rho_2(Z))_+ \right) .
  $$

\subchunk{Four-candidate rules.} For $m\geq4$ it would be tedious to reduce all possible cases
to asymptotic expressions of the kind above. Instead, we have simply illustrated
the sets $\sigma_w M_w$ for a variety of rules $w$ in Figure \fourcandidateexamples.
Note that $\sigma_w M_w$ may have up to $3$ optimal vertices.
(For general $m$, $\sigma_w M_w$ might have up to $m-1$ optimal vertices.)

\pageinsert
\vbox{
\centerline{
\beginpicture
\setcoordinatesystem units <52mm, 52mm>
\setplotarea x from 0 to 0.6, y from 0 to 0.9
\arrow < 2mm>  [ .2679492, .7279404] from 0 0 to 0.66 0
\arrow < 2mm>  [ .2679492, .7279404] from 0 0 to 0 0.99
\put{$\lambda$} [lt] <0pt, -5pt> at 0.66 0
\put{$\mu$} [rB] <-5pt,0pt> at 0 0.99
\put{$w=(1,1,1,0)$ anti-plurality} at 0.3 -0.16
\put{$\bullet$} at 0.433 0.433
\plot 0 0 0.46 0.46 /
\plot 0.433 0.4 0.433 0.9 /
\endpicture
\hfill
\beginpicture
\setcoordinatesystem units <52mm, 52mm>
\setplotarea x from 0 to 0.6, y from 0 to 0.9
\arrow < 2mm>  [ .2679492, .7279404] from 0 0 to 0.66 0
\arrow < 2mm>  [ .2679492, .7279404] from 0 0 to 0 0.99
\put{$\lambda$} [lt] <0pt, -5pt> at 0.66 0
\put{$\mu$} [rB] <-5pt,0pt> at 0 0.99
\put{$w=(1,1,{1\over2},0)$} at 0.3 -0.16
\put{$\bullet$} at 0.415 0.415
\put{$\bullet$} at 0.415 0.829
\plot 0 0 0.44 0.44 /
\plot 0.415 0.38 0.415 0.86 /
\plot 0.44 0.829 0 0.829 /
\endpicture
\hfill
\beginpicture
\setcoordinatesystem units <52mm, 52mm>
\setplotarea x from 0 to 0.6, y from 0 to 0.9
\arrow < 2mm>  [ .2679492, .7279404] from 0 0 to 0.66 0
\arrow < 2mm>  [ .2679492, .7279404] from 0 0 to 0 0.99
\put{$\lambda$} [lt] <0pt, -5pt> at 0.66 0
\put{$\mu$} [rB] <-5pt,0pt> at 0 0.99
\put{$w=(1,{7\over9},{1\over2},0)$} at 0.3 -0.16
\put{$\bullet$} at 0.480 0.480
\put{$\bullet$} at 0.480 0.600
\put{$\bullet$} at 0.415 0.747
\plot 0 0 0.51 0.51 /
\plot 0.480 0.45 0.480 0.63 /
\plot 0 0.747 0.44 0.747 /
\plot 0.395 0.792 0.500 0.555 /
\endpicture
           }
\vfill    
\centerline{
\beginpicture
\setcoordinatesystem units <52mm, 52mm>
\setplotarea x from 0 to 0.6, y from 0 to 0.9
\arrow < 2mm>  [ .2679492, .7279404] from 0 0 to 0.66 0
\arrow < 2mm>  [ .2679492, .7279404] from 0 0 to 0 0.99
\put{$\lambda$} [lt] <0pt, -5pt> at 0.66 0
\put{$\mu$} [rB] <-5pt,0pt> at 0 0.99
\put{$w=(1,{3\over5},{1\over2},0)$} at 0.3 -0.16
\put{$\bullet$} at 0.396 0.396
\plot 0 0 0.43 0.43 /
\plot 0 0.712 0.18 0.712 /
\plot 0.123 0.737 0.416 0.371 /
\endpicture
\hfill
\beginpicture
\setcoordinatesystem units <52mm, 52mm>
\setplotarea x from 0 to 0.6, y from 0 to 0.9
\arrow < 2mm>  [ .2679492, .7279404] from 0 0 to 0.66 0
\arrow < 2mm>  [ .2679492, .7279404] from 0 0 to 0 0.99
\put{$\lambda$} [lt] <0pt, -5pt> at 0.66 0
\put{$\mu$} [rB] <-5pt,0pt> at 0 0.99
\put{$w=(1,{1\over2},{1\over4},0)$} at 0.3 -0.16
\put{$\bullet$} at 0.370 0.370
\put{$\bullet$} at 0.370 0.671
\plot 0 0 0.4 0.4 /
\plot 0.370 0.34 0.370 0.7 /
\plot 0 0.740 0.4 0.666 /
\endpicture
\hfill
\beginpicture
\setcoordinatesystem units <52mm, 52mm>
\setplotarea x from 0 to 0.6, y from 0 to 0.9
\arrow < 2mm>  [ .2679492, .7279404] from 0 0 to 0.66 0
\arrow < 2mm>  [ .2679492, .7279404] from 0 0 to 0 0.99
\put{$\lambda$} [lt] <0pt, -5pt> at 0.66 0
\put{$\mu$} [rB] <-5pt,0pt> at 0 0.99
\put{$w=(1,{1\over2},{1\over2},0)$} at 0.3 -0.16
\put{$\bullet$} at 0.354 0.354
\plot 0 0 0.39 0.39 /
\plot 0 0.707 0.374 0.334 /
\endpicture
           }
\vfill    
\centerline{
\beginpicture
\setcoordinatesystem units <52mm, 52mm>
\setplotarea x from 0 to 0.6, y from 0 to 0.9
\arrow < 2mm>  [ .2679492, .7279404] from 0 0 to 0.66 0
\arrow < 2mm>  [ .2679492, .7279404] from 0 0 to 0 0.99
\put{$\lambda$} [lt] <0pt, -5pt> at 0.66 0
\put{$\mu$} [rB] <-5pt,0pt> at 0 0.99
\put{$w=(1,{2\over3},{1\over3},0)$ Borda} at 0.3 -0.16
\put{$\bullet$} at 0.559 0.559
\plot 0 0 0.59 0.59 /
\plot 0 0.559 0.59 0.559 /
\plot 0.559 0.52 0.559 0.59 /
\plot 0.539 0.579 0.579 0.539 /
\endpicture
\hfill
\beginpicture
\setcoordinatesystem units <52mm, 52mm>
\setplotarea x from 0 to 0.6, y from 0 to 0.9
\arrow < 2mm>  [ .2679492, .7279404] from 0 0 to 0.66 0
\arrow < 2mm>  [ .2679492, .7279404] from 0 0 to 0 0.99
\put{$\lambda$} [lt] <0pt, -5pt> at 0.66 0
\put{$\mu$} [rB] <-5pt,0pt> at 0 0.99
\put{$w=(1,1,0,0)$ 2-approval} at 0.3 -0.16
\put{$\bullet$} at 0.5 0.5
\plot 0 0 0.53 0.53 /
\plot 0.5 0.47 0.5 0.53 /
\plot 0 0.5 0.53 0.5 /
\endpicture
\hfill
\beginpicture
\setcoordinatesystem units <52mm, 52mm>
\setplotarea x from 0 to 0.6, y from 0 to 0.9
\arrow < 2mm>  [ .2679492, .7279404] from 0 0 to 0.66 0
\arrow < 2mm>  [ .2679492, .7279404] from 0 0 to 0 0.99
\put{$\lambda$} [lt] <0pt, -5pt> at 0.66 0
\put{$\mu$} [rB] <-5pt,0pt> at 0 0.99
\put{$w=(1,0,0,0)$ plurality} at 0.3 -0.16
\put{$\bullet$} at 0.433 0.433
\plot 0 0 0.46 0.46 /
\plot 0.433 0.4 0.433 0.46 /
\plot 0 0.433 0.46 0.433 /
\endpicture
           }
\vskip20pt
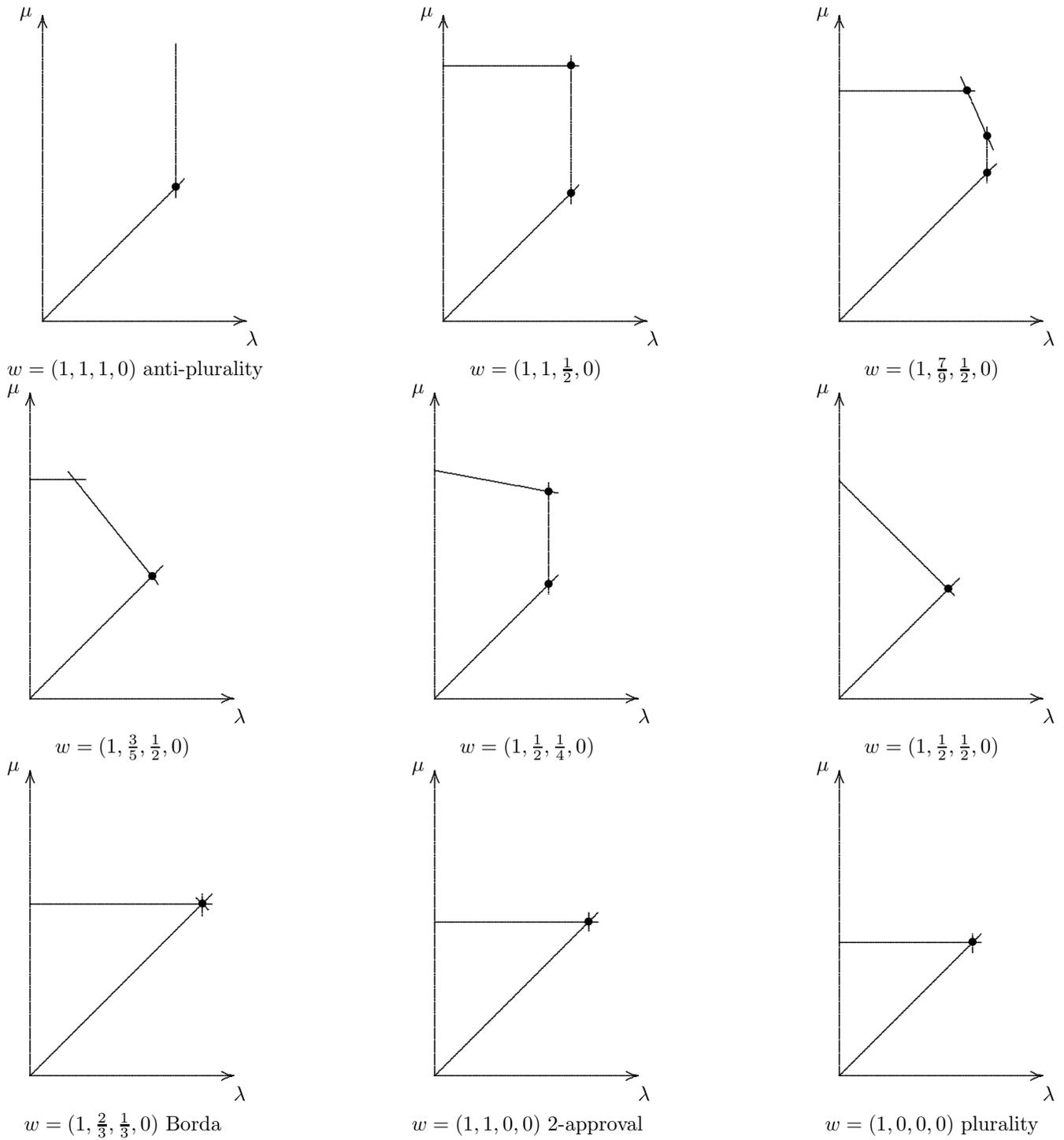
Figure \fourcandidateexamples: the sets $\sigma_w M_w$ for several four-candidate positional rules,
depicted at consistent scales. The black dots show which points may be optimal.
     }
\vfill
\endinsert

The functions $g_w(v)=\Prob{V_w\leq v}$ for some common voting rules are shown in Figures
\gwgthreefour, \gwgfivesix, and \gwgtentwenty.

\topinsert 
\vbox{
\centerline{
\beginpicture
\setcoordinatesystem units <28mm, 60mm>
\setplotarea x from 0 to 2.5, y from 0 to 1
\axis left ticks numbered from 0 to 1 by 0.5 /
\axis bottom ticks numbered from  0 to 2 by 1  /
\arrow < 2mm>  [ .2679492, .7279404] from 0 0 to 2.5 0
\arrow < 2mm>  [ .2679492, .7279404] from 0 0 to 0 1.1
\put{$v$} [lt] <0pt, -5pt> at 2.5 0
\put{$g_w(v)$} [rB] <-5pt,0pt> at 0 1.1
\plot 0 0 0.05 0.04196 0.10 0.08318 0.15 0.1236 0.20 0.1632 0.25 0.2019 0.30 0.2397 0.35 0.2766 0.40 0.3125 0.45 0.3474
      0.50 0.3812 0.55 0.4140 0.60 0.4457 0.65 0.4764 0.70 0.5059 0.75 0.5344 0.80 0.5618 0.85 0.5880 0.90 0.6132 0.95 0.6374
      1.00 0.6604 1.05 0.6824 1.10 0.7033 1.15 0.7232 1.20 0.7422 1.25 0.7601 1.30 0.7771 1.35 0.7932 1.40 0.8083 1.45 0.8226
      1.50 0.8360 1.55 0.8487 1.60 0.8605 1.65 0.8716 1.70 0.8820 1.75 0.8916 1.80 0.9007 1.85 0.9091 1.90 0.9169 1.95 0.9241
      2.00 0.9308 2.05 0.9370 2.10 0.9427 2.15 0.9480 2.20 0.9529 2.25 0.9573 2.30 0.9614 2.35 0.9652 2.40 0.9686 2.45 0.9717
      2.50 0.9746 /
\setdashes
\plot 0 0 0.02887 0.04196 0.05773 0.08318 0.08660 0.1236 0.1155 0.1632 0.1443 0.2019 0.1732 0.2397 0.2021 0.2766 0.2309 0.3125
      0.2598 0.3474 0.2887 0.3812 0.3175 0.4140 0.3464 0.4457 0.3753 0.4764 0.4041 0.5059 0.4330 0.5344 0.4619 0.5618
      0.4907 0.5880 0.5196 0.6132 0.5485 0.6374 0.5773 0.6604 0.6062 0.6824 0.6351 0.7033 0.6639 0.7232 0.6928 0.7422
      0.7217 0.7601 0.7505 0.7771 0.7794 0.7932 0.8083 0.8083 0.8371 0.8226 0.8660 0.8360 0.8949 0.8487 0.9237 0.8605
      0.9526 0.8716 0.9815 0.8820 1.010 0.8916 1.039 0.9007 1.068 0.9091 1.097 0.9169 1.126 0.9241 1.155 0.9308 1.183 0.9370
      1.212 0.9427 1.241 0.9480 1.270 0.9529 1.299 0.9573 1.328 0.9614 1.357 0.9652 1.386 0.9686 1.415 0.9717 1.443 0.9746
      1.472 0.9772 1.501 0.9796 1.530 0.9817 1.559 0.9837 1.588 0.9854 1.616 0.9870 1.645 0.9884 1.674 0.9897 1.703 0.9909
      1.732 0.9919 2.5 1 /
\setdots
\plot 0 0 0.02887 0.04127 0.05773 0.08043 0.08660 0.1175 0.1155 0.1523 0.1443 0.1850 0.1732 0.2156 0.2021 0.2441 0.2309 0.2705
      0.2598 0.2948 0.2887 0.3173 0.3175 0.3378 0.3464 0.3566 0.3753 0.3736 0.4041 0.3890 0.4330 0.4029 0.4619 0.4154
      0.4907 0.4266 0.5196 0.4365 0.5485 0.4453 0.5773 0.4530 0.6062 0.4599 0.6351 0.4658 0.6639 0.4710 0.6928 0.4755
      0.7217 0.4794 0.7505 0.4827 0.7794 0.4856 0.8083 0.4880 0.8371 0.4901 0.8660 0.4918 0.8949 0.4933 0.9237 0.4945
      0.9526 0.4955 0.9815 0.4964 1.010 0.4971 1.039 0.4977 1.068 0.4981 1.097 0.4985 1.126 0.4988 1.155 0.4991 1.183 0.4993
      1.212 0.4994 1.241 0.4996 1.270 0.4997 1.299 0.4997 1.328 0.4998 1.357 0.4998 1.386 0.4999 1.415 0.4999 1.443 0.4999
      1.472 0.4999 2.5 0.5 /
\endpicture
\hfill
\beginpicture
\setcoordinatesystem units <48mm, 60mm>
\setplotarea x from 0 to 1.5, y from 0 to 1
\axis left ticks numbered from 0 to 1 by 0.5 /
\axis bottom ticks numbered from  0 to 1 by 1  /
\arrow < 2mm>  [ .2679492, .7279404] from 0 0 to 1.5 0
\arrow < 2mm>  [ .2679492, .7279404] from 0 0 to 0 1.1
\put{$v$} [lt] <0pt, -5pt> at 1.5 0
\put{$g_w(v)$} [rB] <-5pt,0pt> at 0 1.1
\plot 0 0 0.03227 0.05077 0.06456 0.1001 0.09682 0.1480 0.1291 0.1944 0.1614 0.2393 0.1936 0.2826 0.2259 0.3243 0.2582 0.3645
  0.2905 0.4031 0.3227 0.4401 0.3550 0.4756 0.3873 0.5094 0.4194 0.5418 0.4520 0.5726 0.4841 0.6018 0.5163 0.6296 0.5488 0.6560
  0.5810 0.6809 0.6131 0.7045 0.6456 0.7267 0.6778 0.7476 0.7099 0.7673 0.7425 0.7857 0.7746 0.8030 0.8067 0.8191 0.8393 0.8342
  0.8714 0.8482 0.9036 0.8613 0.9361 0.8734 0.9682 0.8846 1.000 0.8950 1.033 0.9047 1.065 0.9135 1.097 0.9217 1.130 0.9292
  1.162 0.9361 1.194 0.9424 1.227 0.9481 1.259 0.9534 1.291 0.9582 1.323 0.9625 1.356 0.9665 1.388 0.9701 1.420 0.9733
  1.452 0.9762 1.485 0.9789 /
\setdashes
\plot 0 0 0.02500 0.05077 0.05000 0.1001 0.07500 0.1480 0.1000 0.1944 0.1250 0.2393 0.1500 0.2826 0.1750 0.3243 0.2000 0.3645
  0.2250 0.4031 0.2500 0.4401 0.2750 0.4756 0.3000 0.5094 0.3250 0.5418 0.3500 0.5726 0.3750 0.6018 0.4000 0.6296 0.4250 0.6560
  0.4500 0.6809 0.4750 0.7045 0.5000 0.7267 0.5250 0.7476 0.5500 0.7673 0.5750 0.7857 0.6000 0.8030 0.6250 0.8191 0.6500 0.8342
  0.6750 0.8482 0.7000 0.8613 0.7250 0.8734 0.7500 0.8846 0.7750 0.8950 0.8000 0.9047 0.8250 0.9135 0.8500 0.9217 0.8750 0.9292
  0.9000 0.9361 0.9250 0.9424 0.9500 0.9481 0.9750 0.9534 1.000 0.9582 1.025 0.9625 1.050 0.9665 1.075 0.9701 1.100 0.9733
  1.125 0.9762 1.150 0.9789 1.175 0.9813 1.200 0.9834 1.225 0.9853 1.250 0.9870 1.275 0.9886 1.300 0.9899 1.325 0.9912
  1.350 0.9922 1.375 0.9932 1.400 0.9941 1.425 0.9948 1.450 0.9955 1.475 0.9961 1.500 0.9966 /
\setdashpattern <6pt,3pt,1pt,3pt>
\plot 0 0 0.02887 0.05077 0.05773 0.1001 0.08660 0.1480 0.1155 0.1944 0.1443 0.2393 0.1732 0.2826 0.2021 0.3243 0.2309 0.3645
  0.2598 0.4031 0.2887 0.4401 0.3175 0.4756 0.3464 0.5094 0.3753 0.5418 0.4041 0.5726 0.4330 0.6018 0.4619 0.6296 0.4907 0.6560
  0.5196 0.6809 0.5485 0.7045 0.5773 0.7267 0.6062 0.7476 0.6351 0.7673 0.6639 0.7857 0.6928 0.8030 0.7217 0.8191 0.7505 0.8342
  0.7794 0.8482 0.8083 0.8613 0.8371 0.8734 0.8660 0.8846 0.8949 0.8950 0.9237 0.9047 0.9526 0.9135 0.9815 0.9217 1.010 0.9292
  1.039 0.9361 1.068 0.9424 1.097 0.9481 1.126 0.9534 1.155 0.9582 1.183 0.9625 1.212 0.9665 1.241 0.9701 1.270 0.9733
  1.299 0.9762 1.328 0.9789 1.357 0.9813 1.386 0.9834 1.415 0.9853 1.443 0.9870 1.472 0.9886 1.501 0.9899 /
\setdots
\plot 0 0 0.02500 0.050762071 0.05000 0.10003675 0.07500 0.1477876 0.1000 0.1938998 0.1250 0.2383316 0.1500 0.280944
  0.1750 0.321703 0.2000 0.360679 0.2250 0.397747 0.2500 0.432888 0.2750 0.466189 0.3000 0.49744 0.3250 0.52695 0.3500 0.55451
  0.3750 0.58014 0.4000 0.60405 0.4250 0.62626 0.4500 0.64669 0.4750 0.66557 0.5000 0.68284 0.5250 0.69862 0.5500 0.71305
  0.5750 0.72606 0.6000 0.73789 0.6250 0.74847 0.6500 0.75804 0.6750 0.76653 0.7000 0.77417 0.7250 0.78089 0.7500 0.78681
  0.7750 0.7921 0.8000 0.7968 0.8250 0.8007 0.8500 0.8043 0.8750 0.8074 0.9000 0.8100 0.9250 0.8123 0.9500 0.8142
  0.9750 0.8159 1.000 0.8173 1.025 0.8185 1.050 0.8196 1.075 0.8204 1.100 0.8211 1.125 0.8217 1.150 0.8223 1.175 0.8227
  1.200 0.8230 1.225 0.8233 1.250 0.8235 1.275 0.8238 1.300 0.8238 1.325 0.8241 1.350 0.8241 1.375 0.8242 1.400 0.8243
  1.425 0.8243 1.450 0.8244 1.475 0.8244 1.500 0.8244 /
\endpicture
}
\vskip12pt
\centerline{\hfill\smash{$m=3$}\hfill\hfill\smash{$m=4$}\hfill}
\vskip10pt
\centerline{
\hfill
\raise3pt\hbox{\beginpicture
\setcoordinatesystem units <10mm, 2mm>
\setplotarea x from 0 to 1, y from 0 to 1
\setsolid
\plot 0 0 1 0 /
      \endpicture}
Borda
\hfill
\raise3pt\hbox{\beginpicture
\setcoordinatesystem units <10mm, 2mm>
\setplotarea x from 0 to 1, y from 0 to 1
\setdashes
\plot 0 0 1 0 /
      \endpicture}
plurality
\hfill
\raise3pt\hbox{\beginpicture
\setcoordinatesystem units <10mm, 2mm>
\setplotarea x from 0 to 1, y from 0 to 1
\setdashpattern <6pt,3pt,1pt,3pt>
\plot 0 0 1 0 /
      \endpicture}
2-approval
\hfill
\raise3pt\hbox{\beginpicture
\setcoordinatesystem units <10mm, 2mm>
\setplotarea x from 0 to 1, y from 0 to 1
\setdots
\plot 0 0 1 0 /
      \endpicture}
anti-plurality
\hfill
           }
\vskip10pt
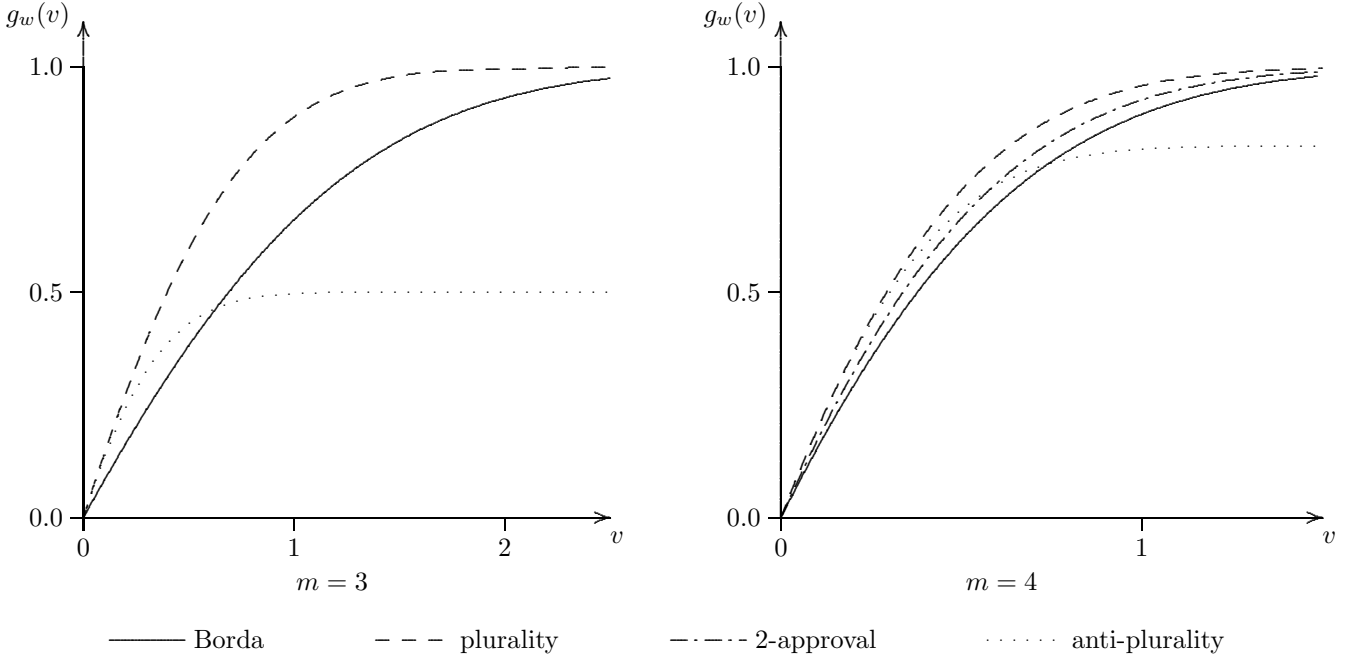
\centerline{Figure \gwgthreefour: the functions $g_w(v)=\Prob{V_w\leq v}$ for some three- and four-candidate voting rules.}
}
\endinsert

\section{5. Comparisons between positional voting rules.}

In this section we compare the various positional voting rules with respect to their manipulability under
IC asymptotic conditions. 

It is apparent from Figure \gwgthreefour\ that the susceptibility of voting rules to
coalitional manipulation depends on the size of the coalition involved. The graph for $m=3$, for example,
shows that elections using the plurality or Borda rules are highly likely to be manipulable by
a large coalition (at least $2\sqrt{n}$ voters), whereas only about half of anti-plurality elections
are so manipulable. But, if one is more concerned about manipulation by small groups,
the graph shows that plurality and anti-plurality elections are about equally susceptible to
manipulation by coalitions of less than $0.25\sqrt{n}$ voters, while Borda is rather less susceptible.

This suggests that there will be no single rule which is clearly superior to all others with respect to
IC coalitional manipulation.
However, some are clearly inferior to others. The asymptotic manipulation probability $g_w(v)$ for the
plurality rule, for example, is greater than for the other two 3-candidate rules shown for all
values of $v$.

More formally, given two $m$-candidate positional voting rules $w$, $w'$,
we will say that $w$ {\it dominates} $w'$ (and write $w'\preceq w$) if $g_w(v)\leq g_{w'}(v)$ for all $v\geq0$.
That is, $w$ is less susceptible than $w'$ to manipulation by coalitions of any given size.
Alternatively, the asymptotic minimum coalition size $V_w$ is larger than $V_{w'}$
in the sense of first-order stochastic dominance. 

Note that $\preceq$ gives a partial order on the rules. Although the present paper is concerned only
with positional rules, the partial order $\preceq$ could be defined in the same way for
any voting rules.

If $w'\preceq w$, then the rule $w$ is to be preferred to $w'$. More generally, the best rules to use
(at least from the point of view of manipulation by IC populations) are those not dominated by any other.

\chunk{Proposition \pluralitydominated.}
The plurality rule is always dominated.

\subchunk{Proof.} Indeed, plurality is dominated both by anti-plurality and by Borda.
This is apparent from the asymptotic results of the previous section:
$$ \Vap \;=\;
   \cases{m^{-1/2} (\rho_1(Z)-\rho_2(Z)), & if $\rho_2(Z)\geq\bar Z$ \cr
          \infty, & otherwise , \cr
         }
  $$
while
$$ \Vp \;=\; m^{-1/2}(\rho_1(Z)-\rho_2(Z)) ,
  $$
giving $\Vp\leq\Vap$ and so (plurality) $\preceq$ (anti-plurality).
Similarly for Borda.
\endofproof

\chunk{Proposition \antipluralityundominated.}
The anti-plurality rule is never dominated.

\subchunk{Proof.}
No other rule $w$ may dominate the anti-plurality rule, because
$\lim_{v\to\infty} g_{\scriptstyle\rm antiplurality}(v) = \Prob{\rho_2(Z)\geq \bar Z} < 1$,
whereas $\lim_{v\to\infty} g_w(v) = 1$.
\endofproof

\subchunk{}
Proposition \antipluralityundominated\ is true because the anti-plurality rule
is resistant to manipulation by very large coalitions, in a way that every other
positional rule is not. However, as we shall see later, this advantage becomes
very slight in elections with 6 or more candidates.

\chunk{Proposition \Bordadomination.}
Borda's rule is undominated for $m\in\set{3,4}$, but dominated for $m\geq5$.

\subchunk{Proof.}
From the asymptotic results of the previous section:
$$ V_\subborda \;=\; \left({m(m+1)\over 12(m-2)^2}\right)^{1/2} (\rho_1(Z)-\rho_2(Z)) ,
  $$
while
$$ \Vhalfmappr \;=\; \left({\grint{m/2}(m-\grint{m/2})\over m(m-1)}\right)^{1/2} (\rho_1(Z)-\rho_2(Z)) .
  $$
Note that for $m\geq5$,
$$ {\grint{m/2}(m-\grint{m/2})\over m(m-1)} \;\geq\; {(m-1)(m+1)\over 4m(m-1)}
   \;=\; 3\left(1-{2\over m}\right)^2 {m(m+1)\over 12(m-2)^2}
   \;\geq\; {m(m+1)\over 12(m-2)^2} .
  $$
This gives $V_\subborda\leq\Vhalfmappr$ and so (Borda) $\preceq$ ($\grint{m/2}$-approval).

Now let $m\in\set{3,4}$ and let $w$ be any positional rule for $m$-candidate elections.
We may observe that
$$ M_w \subseteq M'_w := \theset{(\lambda,\mu)\in\R^2}{0\leq\lambda\leq b_w, \lambda\leq\mu} ,
  $$
where
$$ b_w = \hbox{min}_{i} (1-w_i+w_{i+1})^{-1} .
  $$
It is shown in \cite{\PritchardSlinko} that the quantity $\sigma_w b_w$ is maximized,
over all $m$-candidate positional rules, by Borda's rule.
(This is true only when $m\in\set{3,4}$; for $m\geq5$, the maximum value is achieved by
the $\grint{m/2}$-approval rule.)
From this we obtain
$$ \eqalign{
  g_w(v)
  &= \Prob{V_w \leq v} \cr
  &\geq \Prob{V_w \leq v\hbox{ and } \rho_2(Z)\geq\bar Z} \cr
  &\geq \Prob{\hbox{max}\theset{\lambda(\rho_1(Z)-\bar Z) + \mu(\bar Z - \rho_2(Z))}
            {(\lambda,\mu)\in\left({m\over m-1}\right)^{1/2}\sigma_w M_w} \leq v\hbox{ and } \rho_2(Z)\geq\bar Z} \cr
  &\geq \Prob{\hbox{max}\theset{\lambda(\rho_1(Z)-\bar Z) + \mu(\bar Z - \rho_2(Z))}
            {(\lambda,\mu)\in\left({m\over m-1}\right)^{1/2}\sigma_w M'_w} \leq v\hbox{ and } \rho_2(Z)\geq\bar Z} \cr
  &= \Prob{ \left({m\over m-1}\right)^{1/2} \sigma_w b_w (\rho_1(Z) - \rho_2(Z)) \leq v\hbox{ and } \rho_2(Z)\geq\bar Z} \cr
  &\geq \Prob{ \left({m\over m-1}\right)^{1/2} \sigma_{\subborda} b_{\subborda}
            (\rho_1(Z) - \rho_2(Z)) \leq v\hbox{ and } \rho_2(Z)\geq\bar Z} \cr
  &= \Prob{ V_{\subborda} \leq v\hbox{ and } \rho_2(Z)\geq\bar Z} \cr
  &= g_{\subborda}(v) \;-\; \Prob{V_{\subborda}\leq v\hbox{ and } \rho_2(Z)<\bar Z} . \cr
           }
  $$
Let $c=(m/(m-1))^{1/2}\sigma_{\subborda} b_{\subborda}$.
Since all orderings of $Z=(Z_1,\ldots,Z_m)$ are equally likely,
$$ \eqalign{
  \Prob{V_{\subborda}\leq v\hbox{ and } \rho_2(Z)<\bar Z}
  &= m! \Prob{Z_1\geq Z_2\geq\cdots\geq Z_m, \; c(Z_1-Z_2)\leq v,\; \hbox{ and }
         Z_2< \bar Z} \cr
  &\leq m! \Prob{0\leq \bar Z - Z_2 \leq Z_1-Z_2 \leq v/c} \cr
  &= m! \int_0^{v/c} \int_0^x f(x,y)\;dxdy \cr
  &\leq {1\over2} m! (\sup f) \left({v\over c}\right)^2 \cr
  &= O(v^2) \qquad\qquad\hbox{ as }v\to0 ,
           }
  $$
where $f$ is the (non-degenerate) bivariate normal probability density of $(Z_1-Z_2,\bar Z-Z_2)$.
It follows that
$$ g'_w(0) \geq g'_{\subborda}(0) ,
  $$
and hence that $w$ cannot dominate Borda.
\endofproof

The argument used in Proposition \Bordadomination\ shows that among all positional rules for
three- or four- candidate elections, Borda's rule has the $g$ function with the smallest
derivative at the origin. This means that it is the most resistant to
manipulation by very small coalitions. However, this property does not hold when there
are five or more candidates. In that case, it is the $\grint{m/2}$-approval rule
which enjoys maximal resistance to manipulation by very small coalitions.
These results are similar to those of \cite{\PritchardSlinko}, although the criterion
considered there is the ``average threshhold coalition size" rather than the minimum
manipulating coalition size of the present paper.

\pageinsert 
\vbox{
\centerline{
\beginpicture
\setcoordinatesystem units <58mm, 60mm>
\setplotarea x from 0 to 1.2, y from 0 to 1
\axis left ticks numbered from 0 to 1 by 0.5 /
\axis bottom ticks numbered from  0 to 1 by 1  /
\arrow < 2mm>  [ .2679492, .7279404] from 0 0 to 1.2 0
\arrow < 2mm>  [ .2679492, .7279404] from 0 0 to 0 1.1
\put{$v$} [lt] <0pt, -5pt> at 1.2 0
\put{$g_w(v)$} [rB] <-5pt,0pt> at 0 1.1
\plot
 0.0000 0.0000  0.0091 0.0200  0.0184 0.0400  0.0280 0.0600  0.0371 0.0800  0.0467 0.1000  0.0566 0.1200
 0.0665 0.1400  0.0763 0.1600  0.0862 0.1800  0.0969 0.2000  0.1074 0.2200  0.1181 0.2400  0.1289 0.2600
 0.1395 0.2800  0.1508 0.3000  0.1622 0.3200  0.1739 0.3400  0.1860 0.3600  0.1981 0.3800  0.2105 0.4000
 0.2232 0.4200  0.2361 0.4400  0.2497 0.4600  0.2635 0.4800  0.2774 0.5000  0.2920 0.5200  0.3073 0.5400
 0.3232 0.5600  0.3399 0.5800  0.3569 0.6000  0.3738 0.6200  0.3919 0.6400  0.4107 0.6600  0.4302 0.6800
 0.4512 0.7000  0.4729 0.7200  0.4969 0.7400  0.5214 0.7600  0.5479 0.7800  0.5772 0.8000  0.6081 0.8200
 0.6400 0.8400  0.6773 0.8600  0.7198 0.8800  0.7676 0.9000  0.8240 0.9200  0.8947 0.9400  0.9871 0.9600
 0.9927 0.9610  1.0000 0.9620  1.0057 0.9630  1.0115 0.9640  1.0169 0.9650
 1.0228 0.9660  1.0292 0.9670  1.0354 0.9680  1.0416 0.9690  1.0490 0.9700  1.0558 0.9710  1.0644 0.9720
 1.0727 0.9730  1.0804 0.9740  1.0881 0.9750  1.0976 0.9760  1.1052 0.9770  1.1146 0.9780  1.1233 0.9790
 1.1341 0.9800  1.1454 0.9810  1.1572 0.9820  1.1680 0.9830  1.1813 0.9840  1.1949 0.9850 / 
\setdashes
\plot
 0.0000 0.0000  0.0078 0.0200  0.0156 0.0400  0.0237 0.0600  0.0315 0.0800  0.0397 0.1000  0.0480 0.1200
 0.0564 0.1400  0.0647 0.1600  0.0732 0.1800  0.0822 0.2000  0.0912 0.2200  0.1002 0.2400  0.1093 0.2600
 0.1184 0.2800  0.1280 0.3000  0.1377 0.3200  0.1475 0.3400  0.1579 0.3600  0.1681 0.3800  0.1786 0.4000
 0.1894 0.4200  0.2003 0.4400  0.2119 0.4600  0.2236 0.4800  0.2353 0.5000  0.2478 0.5200  0.2608 0.5400
 0.2743 0.5600  0.2884 0.5800  0.3028 0.6000  0.3172 0.6200  0.3325 0.6400  0.3485 0.6600  0.3650 0.6800
 0.3829 0.7000  0.4013 0.7200  0.4216 0.7400  0.4425 0.7600  0.4649 0.7800  0.4898 0.8000  0.5160 0.8200
 0.5431 0.8400  0.5747 0.8600  0.6107 0.8800  0.6514 0.9000  0.6992 0.9200  0.7592 0.9400  0.8376 0.9600
 0.8424 0.9610  0.8485 0.9620  0.8533 0.9630  0.8583 0.9640  0.8629 0.9650
 0.8679 0.9660  0.8733 0.9670  0.8786 0.9680  0.8839 0.9690  0.8901 0.9700  0.8959 0.9710  0.9032 0.9720
 0.9102 0.9730  0.9168 0.9740  0.9233 0.9750  0.9314 0.9760  0.9378 0.9770  0.9458 0.9780  0.9531 0.9790
 0.9623 0.9800  0.9719 0.9810  0.9819 0.9820  0.9911 0.9830  1.0023 0.9840  1.0139 0.9850  1.0285 0.9860
 1.0404 0.9870  1.0534 0.9880  1.0682 0.9890  1.0829 0.9900  1.1026 0.9910  1.1205 0.9920  1.1442 0.9930
 1.1710 0.9940  1.1982 0.9950 /  
\setdashpattern <6pt,3pt,1pt,3pt>
\plot
 0.0000 0.0000  0.0095 0.0200  0.0191 0.0400  0.0291 0.0600  0.0386 0.0800  0.0486 0.1000  0.0588 0.1200
 0.0691 0.1400  0.0793 0.1600  0.0896 0.1800  0.1007 0.2000  0.1116 0.2200  0.1227 0.2400  0.1339 0.2600
 0.1450 0.2800  0.1567 0.3000  0.1686 0.3200  0.1807 0.3400  0.1933 0.3600  0.2059 0.3800  0.2187 0.4000
 0.2320 0.4200  0.2453 0.4400  0.2595 0.4600  0.2738 0.4800  0.2882 0.5000  0.3035 0.5200  0.3194 0.5400
 0.3359 0.5600  0.3533 0.5800  0.3709 0.6000  0.3885 0.6200  0.4072 0.6400  0.4269 0.6600  0.4470 0.6800
 0.4689 0.7000  0.4915 0.7200  0.5164 0.7400  0.5419 0.7600  0.5694 0.7800  0.5999 0.8000  0.6320 0.8200
 0.6651 0.8400  0.7038 0.8600  0.7480 0.8800  0.7978 0.9000  0.8563 0.9200  0.9298 0.9400  1.0258 0.9600
 1.0317 0.9610  1.0392 0.9620  1.0451 0.9630  1.0512 0.9640  1.0568 0.9650
 1.0629 0.9660  1.0696 0.9670  1.0761 0.9680  1.0825 0.9690  1.0902 0.9700  1.0972 0.9710  1.1062 0.9720
 1.1148 0.9730  1.1228 0.9740  1.1308 0.9750  1.1407 0.9760  1.1486 0.9770  1.1583 0.9780  1.1674 0.9790
 1.1786 0.9800  1.1904 0.9810  1.2026 0.9820  /
\setdots
\plot
 0.0000 0.0000  0.0021 0.0050  0.0039 0.0100  0.0058 0.0150  0.0078 0.0200  0.0097 0.0250  0.0118 0.0300
 0.0137 0.0350  0.0156 0.0400  0.0177 0.0450  0.0197 0.0500  0.0217 0.0550  0.0237 0.0600  0.0257 0.0650
 0.0276 0.0700  0.0295 0.0750  0.0315 0.0800  0.0334 0.0850  0.0354 0.0900  0.0376 0.0950  0.0397 0.1000
 0.0417 0.1050  0.0438 0.1100  0.0459 0.1150  0.0480 0.1200  0.0502 0.1250  0.0521 0.1300  0.0542 0.1350
 0.0564 0.1400  0.0585 0.1450  0.0606 0.1500  0.0627 0.1550  0.0647 0.1600  0.0668 0.1650  0.0688 0.1700
 0.0710 0.1750  0.0732 0.1800  0.0754 0.1850  0.0776 0.1900  0.0800 0.1950  0.0822 0.2000  0.0843 0.2050
 0.0866 0.2100  0.0889 0.2150  0.0912 0.2200  0.0935 0.2250  0.0958 0.2300  0.0980 0.2350  0.1002 0.2400
 0.1024 0.2450  0.1046 0.2500  0.1070 0.2550  0.1093 0.2600  0.1117 0.2650  0.1139 0.2700  0.1161 0.2750
 0.1184 0.2800  0.1207 0.2850  0.1233 0.2900  0.1256 0.2950  0.1280 0.3000  0.1305 0.3050  0.1329 0.3100
 0.1353 0.3150  0.1377 0.3200  0.1403 0.3250  0.1429 0.3300  0.1452 0.3350  0.1476 0.3400  0.1500 0.3450
 0.1527 0.3500  0.1554 0.3550  0.1579 0.3600  0.1604 0.3650  0.1630 0.3700  0.1656 0.3750  0.1682 0.3800
 0.1707 0.3850  0.1733 0.3900  0.1761 0.3950  0.1788 0.4000  0.1814 0.4050  0.1842 0.4100  0.1870 0.4150
 0.1896 0.4200  0.1922 0.4250  0.1949 0.4300  0.1977 0.4350  0.2006 0.4400  0.2035 0.4450  0.2066 0.4500
 0.2094 0.4550  0.2122 0.4600  0.2151 0.4650  0.2182 0.4700  0.2211 0.4750  0.2240 0.4800  0.2268 0.4850
 0.2297 0.4900  0.2329 0.4950  0.2358 0.5000  0.2389 0.5050  0.2421 0.5100  0.2451 0.5150  0.2483 0.5200
 0.2516 0.5250  0.2549 0.5300  0.2580 0.5350  0.2613 0.5400  0.2649 0.5450  0.2684 0.5500  0.2719 0.5550
 0.2752 0.5600  0.2786 0.5650  0.2823 0.5700  0.2860 0.5750  0.2895 0.5800  0.2931 0.5850  0.2969 0.5900
 0.3005 0.5950  0.3042 0.6000  0.3078 0.6050  0.3116 0.6100  0.3152 0.6150  0.3188 0.6200  0.3227 0.6250
 0.3265 0.6300  0.3305 0.6350  0.3343 0.6400  0.3383 0.6450  0.3424 0.6500  0.3465 0.6550  0.3509 0.6600
 0.3550 0.6650  0.3595 0.6700  0.3635 0.6750  0.3679 0.6800  0.3727 0.6850  0.3771 0.6900  0.3818 0.6950
 0.3867 0.7000  0.3916 0.7050  0.3964 0.7100  0.4008 0.7150  0.4059 0.7200  0.4113 0.7250  0.4166 0.7300
 0.4220 0.7350  0.4272 0.7400  0.4324 0.7450  0.4380 0.7500  0.4438 0.7550  0.4494 0.7600  0.4549 0.7650
 0.4613 0.7700  0.4677 0.7750  0.4743 0.7800  0.4811 0.7850  0.4878 0.7900  0.4939 0.7950  0.5012 0.8000
 0.5085 0.8050  0.5156 0.8100  0.5232 0.8150  0.5305 0.8200  0.5382 0.8250  0.5463 0.8300  0.5543 0.8350
 0.5627 0.8400  0.5720 0.8450  0.5816 0.8500  0.5915 0.8550  0.6015 0.8600  0.6120 0.8650  0.6229 0.8700
 0.6344 0.8750  0.6465 0.8800  0.6603 0.8850  0.6748 0.8900  0.6896 0.8950  0.7052 0.9000  0.7230 0.9050
 0.7417 0.9100  0.7634 0.9150  0.7863 0.9200  0.8136 0.9250  0.8478 0.9300  0.8846 0.9350  0.9389 0.9400
 1.0329 0.9450  1.2 0.947 /
\endpicture
\hfill
\beginpicture
\setcoordinatesystem units <58mm, 60mm>
\setplotarea x from 0 to 1.2, y from 0 to 1
\axis left ticks numbered from 0 to 1 by 0.5 /
\axis bottom ticks numbered from  0 to 1 by 1  /
\arrow < 2mm>  [ .2679492, .7279404] from 0 0 to 1.2 0
\arrow < 2mm>  [ .2679492, .7279404] from 0 0 to 0 1.1
\put{$v$} [lt] <0pt, -5pt> at 1.2 0
\put{$g_w(v)$} [rB] <-5pt,0pt> at 0 1.1
\plot
 0.0000 0.0000  0.0072 0.0200  0.0149 0.0400  0.0224 0.0600  0.0304 0.0800  0.0380 0.1000  0.0461 0.1200
 0.0542 0.1400  0.0625 0.1600  0.0708 0.1800  0.0794 0.2000  0.0879 0.2200  0.0968 0.2400  0.1059 0.2600
 0.1151 0.2800  0.1242 0.3000  0.1340 0.3200  0.1434 0.3400  0.1534 0.3600  0.1635 0.3800  0.1739 0.4000
 0.1840 0.4200  0.1951 0.4400  0.2060 0.4600  0.2175 0.4800  0.2290 0.5000  0.2408 0.5200  0.2536 0.5400
 0.2670 0.5600  0.2805 0.5800  0.2948 0.6000  0.3090 0.6200  0.3242 0.6400  0.3395 0.6600  0.3562 0.6800
 0.3737 0.7000  0.3920 0.7200  0.4117 0.7400  0.4322 0.7600  0.4542 0.7800  0.4784 0.8000  0.5044 0.8200
 0.5326 0.8400  0.5634 0.8600  0.5985 0.8800  0.6385 0.9000  0.6868 0.9200  0.7473 0.9400  0.8258 0.9600
 0.8306 0.9610  0.8352 0.9620  0.8390 0.9630  0.8441 0.9640  0.8496 0.9650
 0.8556 0.9660  0.8612 0.9670  0.8660 0.9680  0.8723 0.9690  0.8785 0.9700  0.8839 0.9710  0.8899 0.9720
 0.8978 0.9730  0.9048 0.9740  0.9116 0.9750  0.9202 0.9760  0.9274 0.9770  0.9347 0.9780  0.9428 0.9790
 0.9522 0.9800  0.9601 0.9810  0.9698 0.9820  0.9802 0.9830  0.9891 0.9840  1.0023 0.9850  1.0131 0.9860
 1.0280 0.9870  1.0428 0.9880  1.0576 0.9890  1.0733 0.9900  1.0898 0.9910  1.1079 0.9920  1.1293 0.9930
 1.1537 0.9940  1.1858 0.9950  1.2 0.9954 /  
\setdashes
\plot
 0.0000 0.0000  0.0063 0.0200  0.0130 0.0400  0.0196 0.0600  0.0266 0.0800  0.0331 0.1000  0.0402 0.1200
 0.0473 0.1400  0.0546 0.1600  0.0618 0.1800  0.0693 0.2000  0.0767 0.2200  0.0845 0.2400  0.0924 0.2600
 0.1005 0.2800  0.1084 0.3000  0.1169 0.3200  0.1252 0.3400  0.1339 0.3600  0.1427 0.3800  0.1518 0.4000
 0.1606 0.4200  0.1703 0.4400  0.1798 0.4600  0.1899 0.4800  0.1999 0.5000  0.2102 0.5200  0.2214 0.5400
 0.2331 0.5600  0.2448 0.5800  0.2573 0.6000  0.2697 0.6200  0.2830 0.6400  0.2963 0.6600  0.3109 0.6800
 0.3262 0.7000  0.3422 0.7200  0.3594 0.7400  0.3772 0.7600  0.3965 0.7800  0.4176 0.8000  0.4403 0.8200
 0.4649 0.8400  0.4918 0.8600  0.5225 0.8800  0.5574 0.9000  0.5995 0.9200  0.6523 0.9400  0.7208 0.9600
 0.7250 0.9610  0.7290 0.9620  0.7323 0.9630  0.7368 0.9640  0.7416 0.9650
 0.7468 0.9660  0.7517 0.9670  0.7559 0.9680  0.7614 0.9690  0.7668 0.9700  0.7715 0.9710  0.7768 0.9720
 0.7836 0.9730  0.7898 0.9740  0.7957 0.9750  0.8032 0.9760  0.8095 0.9770  0.8159 0.9780  0.8229 0.9790
 0.8311 0.9800  0.8380 0.9810  0.8465 0.9820  0.8556 0.9830  0.8634 0.9840  0.8749 0.9850  0.8843 0.9860
 0.8973 0.9870  0.9102 0.9880  0.9232 0.9890  0.9368 0.9900  0.9512 0.9910  0.9670 0.9920  0.9857 0.9930
 1.0070 0.9940  1.0350 0.9950  1.0643 0.9960  1.1029 0.9970  1.1578 0.9980  1.2 0.9985 /  
\setdashpattern <6pt,3pt,1pt,3pt>
\plot
 0.0000 0.0000  0.0084 0.0200  0.0175 0.0400  0.0263 0.0600  0.0357 0.0800  0.0445 0.1000  0.0539 0.1200
 0.0634 0.1400  0.0732 0.1600  0.0829 0.1800  0.0930 0.2000  0.1029 0.2200  0.1133 0.2400  0.1240 0.2600
 0.1348 0.2800  0.1454 0.3000  0.1569 0.3200  0.1679 0.3400  0.1797 0.3600  0.1915 0.3800  0.2037 0.4000
 0.2154 0.4200  0.2285 0.4400  0.2413 0.4600  0.2547 0.4800  0.2682 0.5000  0.2820 0.5200  0.2970 0.5400
 0.3127 0.5600  0.3285 0.5800  0.3452 0.6000  0.3618 0.6200  0.3797 0.6400  0.3976 0.6600  0.4171 0.6800
 0.4376 0.7000  0.4590 0.7200  0.4822 0.7400  0.5061 0.7600  0.5319 0.7800  0.5603 0.8000  0.5907 0.8200
 0.6237 0.8400  0.6598 0.8600  0.7009 0.8800  0.7478 0.9000  0.8043 0.9200  0.8751 0.9400  0.9670 0.9600
 0.9727 0.9610  0.9781 0.9620  0.9825 0.9630  0.9886 0.9640  0.9949 0.9650
 1.0019 0.9660  1.0085 0.9670  1.0141 0.9680  1.0215 0.9690  1.0287 0.9700  1.0351 0.9710  1.0422 0.9720
 1.0514 0.9730  1.0596 0.9740  1.0675 0.9750  1.0776 0.9760  1.0861 0.9770  1.0946 0.9780  1.1041 0.9790
 1.1151 0.9800  1.1243 0.9810  1.1357 0.9820  1.1479 0.9830  1.1584 0.9840  1.1738 0.9850  1.1864 0.9860
 1.2038 0.9870  /
\setdots
\plot
 0.0000 0.0000  0.0017 0.0050  0.0033 0.0100  0.0048 0.0150  0.0063 0.0200  0.0079 0.0250  0.0096 0.0300
 0.0112 0.0350  0.0130 0.0400  0.0147 0.0450  0.0163 0.0500  0.0180 0.0550  0.0196 0.0600  0.0213 0.0650
 0.0231 0.0700  0.0248 0.0750  0.0266 0.0800  0.0282 0.0850  0.0298 0.0900  0.0314 0.0950  0.0331 0.1000
 0.0350 0.1050  0.0368 0.1100  0.0386 0.1150  0.0402 0.1200  0.0419 0.1250  0.0437 0.1300  0.0455 0.1350
 0.0473 0.1400  0.0490 0.1450  0.0509 0.1500  0.0527 0.1550  0.0546 0.1600  0.0562 0.1650  0.0581 0.1700
 0.0598 0.1750  0.0618 0.1800  0.0636 0.1850  0.0656 0.1900  0.0675 0.1950  0.0693 0.2000  0.0712 0.2050
 0.0731 0.2100  0.0749 0.2150  0.0767 0.2200  0.0786 0.2250  0.0806 0.2300  0.0824 0.2350  0.0845 0.2400
 0.0865 0.2450  0.0886 0.2500  0.0905 0.2550  0.0924 0.2600  0.0944 0.2650  0.0965 0.2700  0.0984 0.2750
 0.1005 0.2800  0.1024 0.2850  0.1045 0.2900  0.1065 0.2950  0.1084 0.3000  0.1105 0.3050  0.1128 0.3100
 0.1147 0.3150  0.1169 0.3200  0.1190 0.3250  0.1211 0.3300  0.1231 0.3350  0.1252 0.3400  0.1273 0.3450
 0.1295 0.3500  0.1317 0.3550  0.1339 0.3600  0.1362 0.3650  0.1385 0.3700  0.1406 0.3750  0.1427 0.3800
 0.1450 0.3850  0.1472 0.3900  0.1495 0.3950  0.1518 0.4000  0.1541 0.4050  0.1562 0.4100  0.1584 0.4150
 0.1606 0.4200  0.1629 0.4250  0.1655 0.4300  0.1681 0.4350  0.1703 0.4400  0.1729 0.4450  0.1750 0.4500
 0.1774 0.4550  0.1799 0.4600  0.1822 0.4650  0.1848 0.4700  0.1873 0.4750  0.1899 0.4800  0.1924 0.4850
 0.1949 0.4900  0.1974 0.4950  0.2000 0.5000  0.2025 0.5050  0.2051 0.5100  0.2076 0.5150  0.2102 0.5200
 0.2131 0.5250  0.2157 0.5300  0.2184 0.5350  0.2214 0.5400  0.2244 0.5450  0.2273 0.5500  0.2303 0.5550
 0.2331 0.5600  0.2360 0.5650  0.2389 0.5700  0.2419 0.5750  0.2450 0.5800  0.2480 0.5850  0.2512 0.5900
 0.2543 0.5950  0.2574 0.6000  0.2605 0.6050  0.2636 0.6100  0.2667 0.6150  0.2698 0.6200  0.2730 0.6250
 0.2763 0.6300  0.2797 0.6350  0.2831 0.6400  0.2864 0.6450  0.2895 0.6500  0.2930 0.6550  0.2964 0.6600
 0.3001 0.6650  0.3037 0.6700  0.3073 0.6750  0.3111 0.6800  0.3151 0.6850  0.3187 0.6900  0.3224 0.6950
 0.3264 0.7000  0.3301 0.7050  0.3338 0.7100  0.3381 0.7150  0.3424 0.7200  0.3467 0.7250  0.3510 0.7300
 0.3555 0.7350  0.3598 0.7400  0.3641 0.7450  0.3687 0.7500  0.3730 0.7550  0.3776 0.7600  0.3822 0.7650
 0.3872 0.7700  0.3920 0.7750  0.3972 0.7800  0.4023 0.7850  0.4076 0.7900  0.4125 0.7950  0.4184 0.8000
 0.4239 0.8050  0.4296 0.8100  0.4354 0.8150  0.4415 0.8200  0.4475 0.8250  0.4538 0.8300  0.4600 0.8350
 0.4663 0.8400  0.4728 0.8450  0.4798 0.8500  0.4867 0.8550  0.4939 0.8600  0.5013 0.8650  0.5090 0.8700
 0.5171 0.8750  0.5251 0.8800  0.5333 0.8850  0.5419 0.8900  0.5517 0.8950  0.5617 0.9000  0.5715 0.9050
 0.5818 0.9100  0.5939 0.9150  0.6062 0.9200  0.6194 0.9250  0.6337 0.9300  0.6477 0.9350  0.6632 0.9400
 0.6800 0.9450  0.6989 0.9500  0.7197 0.9550  0.7423 0.9600  0.7697 0.9650  0.8043 0.9700  0.8460 0.9750
 0.9077 0.9800  1.0178 0.9850  1.2 0.9875 /   
\endpicture
}
\vskip12pt
\centerline{\hfill\smash{$m=5$}\hfill\hfill\smash{$m=6$}\hfill}
\vskip10pt
\centerline{
\hfill
\raise3pt\hbox{\beginpicture
\setcoordinatesystem units <10mm, 2mm>
\setplotarea x from 0 to 1, y from 0 to 1
\setsolid
\plot 0 0 1 0 /
      \endpicture}
Borda
\hfill
\raise3pt\hbox{\beginpicture
\setcoordinatesystem units <10mm, 2mm>
\setplotarea x from 0 to 1, y from 0 to 1
\setdashes
\plot 0 0 1 0 /
      \endpicture}
plurality
\hfill
\raise3pt\hbox{\beginpicture
\setcoordinatesystem units <10mm, 2mm>
\setplotarea x from 0 to 1, y from 0 to 1
\setdashpattern <6pt,3pt,1pt,3pt>
\plot 0 0 1 0 /
      \endpicture}
3-approval
\hfill
\raise3pt\hbox{\beginpicture
\setcoordinatesystem units <10mm, 2mm>
\setplotarea x from 0 to 1, y from 0 to 1
\setdots
\plot 0 0 1 0 /
      \endpicture}
anti-plurality
\hfill
           }
\vskip10pt
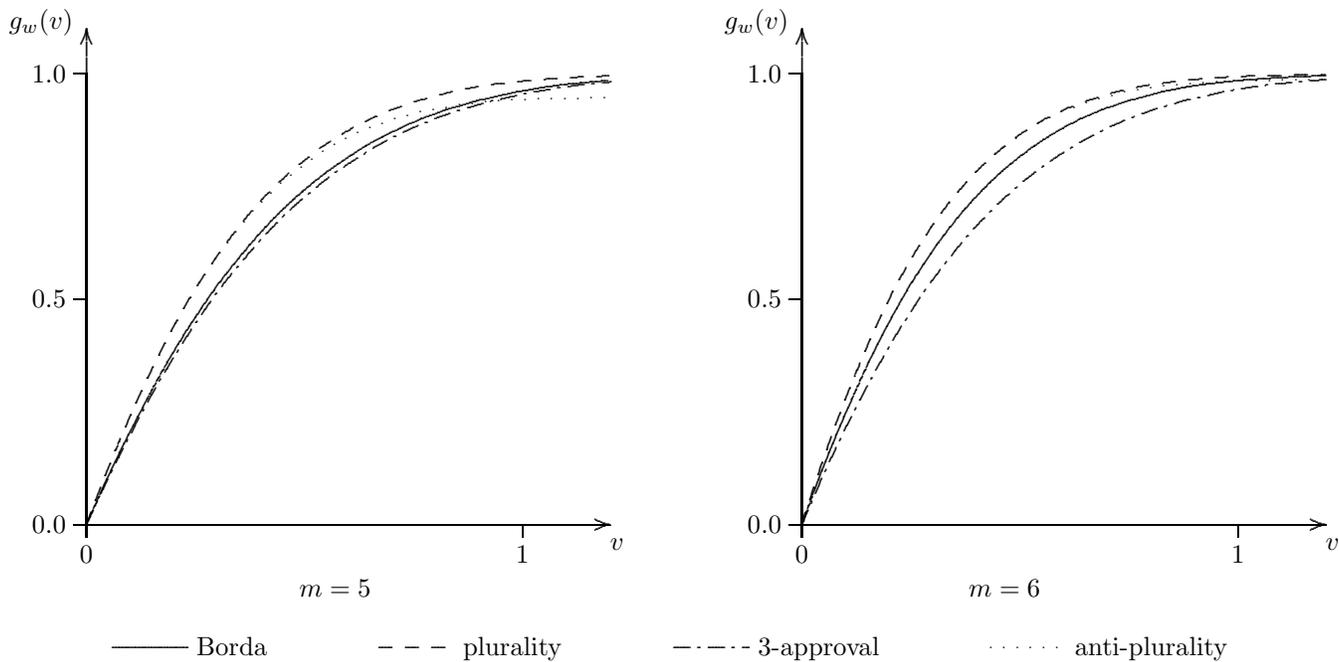
\centerline{Figure \gwgfivesix: the functions $g_w(v)=\Prob{V_w\leq v}$ for some five- and six-candidate voting rules.}
}
\vfill
\vbox{
\centerline{
\beginpicture
\setcoordinatesystem units <100mm, 60mm>
\setplotarea x from 0 to 0.72, y from 0 to 1
\axis left ticks numbered from 0 to 1 by 0.5 /
\axis bottom ticks numbered from  0 to 0.5 by 0.5  /
\arrow < 2mm>  [ .2679492, .7279404] from 0 0 to 0.72 0
\arrow < 2mm>  [ .2679492, .7279404] from 0 0 to 0 1.1
\put{$v$} [lt] <0pt, -5pt> at 0.72 0
\put{$g_w(v)$} [rB] <-5pt,0pt> at 0 1.1
\plot
 0.0000 0.0000  0.0025 0.0100  0.0049 0.0200  0.0074 0.0300  0.0098 0.0400  0.0123 0.0500  0.0149 0.0600
 0.0174 0.0700  0.0200 0.0800  0.0227 0.0900  0.0253 0.1000  0.0279 0.1100  0.0306 0.1200  0.0334 0.1300
 0.0361 0.1400  0.0387 0.1500  0.0415 0.1600  0.0442 0.1700  0.0472 0.1800  0.0500 0.1900  0.0529 0.2000
 0.0558 0.2100  0.0587 0.2200  0.0618 0.2300  0.0648 0.2400  0.0678 0.2500  0.0709 0.2600  0.0738 0.2700
 0.0770 0.2800  0.0802 0.2900  0.0835 0.3000  0.0869 0.3100  0.0901 0.3200  0.0934 0.3300  0.0969 0.3400
 0.1004 0.3500  0.1038 0.3600  0.1073 0.3700  0.1109 0.3800  0.1144 0.3900  0.1180 0.4000  0.1215 0.4100
 0.1255 0.4200  0.1290 0.4300  0.1328 0.4400  0.1366 0.4500  0.1406 0.4600  0.1445 0.4700  0.1484 0.4800
 0.1527 0.4900  0.1571 0.5000  0.1614 0.5100  0.1657 0.5200  0.1699 0.5300  0.1743 0.5400  0.1786 0.5500
 0.1832 0.5600  0.1879 0.5700  0.1924 0.5800  0.1972 0.5900  0.2022 0.6000  0.2073 0.6100  0.2129 0.6200
 0.2182 0.6300  0.2235 0.6400  0.2288 0.6500  0.2344 0.6600  0.2402 0.6700  0.2459 0.6800  0.2518 0.6900
 0.2581 0.7000  0.2644 0.7100  0.2709 0.7200  0.2775 0.7300  0.2845 0.7400  0.2919 0.7500  0.2991 0.7600
 0.3067 0.7700  0.3146 0.7800  0.3226 0.7900  0.3311 0.8000  0.3402 0.8100  0.3499 0.8200  0.3598 0.8300
 0.3698 0.8400  0.3814 0.8500  0.3929 0.8600  0.4056 0.8700  0.4182 0.8800  0.4329 0.8900  0.4476 0.9000
 0.4641 0.9100  0.4820 0.9200  0.5019 0.9300  0.5251 0.9400  0.5536 0.9500  0.5878 0.9600  0.6279 0.9700
 0.6803 0.9800  0.72 0.9844 /  
\setdashes
\plot
 0.0000 0.0000  0.0021 0.0100  0.0041 0.0200  0.0062 0.0300  0.0082 0.0400  0.0103 0.0500  0.0124 0.0600
 0.0145 0.0700  0.0167 0.0800  0.0190 0.0900  0.0211 0.1000  0.0234 0.1100  0.0256 0.1200  0.0279 0.1300
 0.0301 0.1400  0.0324 0.1500  0.0347 0.1600  0.0369 0.1700  0.0394 0.1800  0.0418 0.1900  0.0442 0.2000
 0.0466 0.2100  0.0491 0.2200  0.0516 0.2300  0.0541 0.2400  0.0566 0.2500  0.0592 0.2600  0.0617 0.2700
 0.0643 0.2800  0.0670 0.2900  0.0698 0.3000  0.0726 0.3100  0.0753 0.3200  0.0781 0.3300  0.0809 0.3400
 0.0839 0.3500  0.0867 0.3600  0.0897 0.3700  0.0926 0.3800  0.0956 0.3900  0.0986 0.4000  0.1015 0.4100
 0.1048 0.4200  0.1078 0.4300  0.1110 0.4400  0.1142 0.4500  0.1175 0.4600  0.1207 0.4700  0.1240 0.4800
 0.1276 0.4900  0.1312 0.5000  0.1349 0.5100  0.1384 0.5200  0.1420 0.5300  0.1456 0.5400  0.1493 0.5500
 0.1531 0.5600  0.1570 0.5700  0.1608 0.5800  0.1648 0.5900  0.1689 0.6000  0.1732 0.6100  0.1779 0.6200
 0.1823 0.6300  0.1867 0.6400  0.1912 0.6500  0.1958 0.6600  0.2007 0.6700  0.2054 0.6800  0.2104 0.6900
 0.2157 0.7000  0.2209 0.7100  0.2263 0.7200  0.2319 0.7300  0.2377 0.7400  0.2439 0.7500  0.2499 0.7600
 0.2562 0.7700  0.2629 0.7800  0.2696 0.7900  0.2767 0.8000  0.2843 0.8100  0.2924 0.8200  0.3007 0.8300
 0.3090 0.8400  0.3187 0.8500  0.3283 0.8600  0.3389 0.8700  0.3495 0.8800  0.3617 0.8900  0.3740 0.9000
 0.3878 0.9100  0.4027 0.9200  0.4194 0.9300  0.4388 0.9400  0.4626 0.9500  0.4912 0.9600  0.5246 0.9700
 0.5685 0.9800  0.6432 0.9900  0.6532 0.9910  0.6652 0.9920  0.6768 0.9930  0.6922 0.9940  0.7107 0.9950
 0.72 0.9954 /    
\setdashpattern <6pt,3pt,1pt,3pt>
\plot
 0.0000 0.0000  0.0035 0.0100  0.0069 0.0200  0.0103 0.0300  0.0137 0.0400  0.0172 0.0500  0.0207 0.0600
 0.0242 0.0700  0.0279 0.0800  0.0316 0.0900  0.0352 0.1000  0.0389 0.1100  0.0426 0.1200  0.0465 0.1300
 0.0502 0.1400  0.0539 0.1500  0.0578 0.1600  0.0615 0.1700  0.0657 0.1800  0.0696 0.1900  0.0737 0.2000
 0.0777 0.2100  0.0818 0.2200  0.0861 0.2300  0.0902 0.2400  0.0944 0.2500  0.0987 0.2600  0.1028 0.2700
 0.1072 0.2800  0.1117 0.2900  0.1163 0.3000  0.1210 0.3100  0.1255 0.3200  0.1301 0.3300  0.1349 0.3400
 0.1398 0.3500  0.1446 0.3600  0.1495 0.3700  0.1544 0.3800  0.1593 0.3900  0.1643 0.4000  0.1692 0.4100
 0.1747 0.4200  0.1797 0.4300  0.1849 0.4400  0.1903 0.4500  0.1958 0.4600  0.2012 0.4700  0.2066 0.4800
 0.2126 0.4900  0.2187 0.5000  0.2248 0.5100  0.2307 0.5200  0.2367 0.5300  0.2427 0.5400  0.2488 0.5500
 0.2551 0.5600  0.2617 0.5700  0.2680 0.5800  0.2747 0.5900  0.2816 0.6000  0.2886 0.6100  0.2965 0.6200
 0.3038 0.6300  0.3112 0.6400  0.3187 0.6500  0.3264 0.6600  0.3345 0.6700  0.3424 0.6800  0.3506 0.6900
 0.3595 0.7000  0.3682 0.7100  0.3772 0.7200  0.3865 0.7300  0.3962 0.7400  0.4065 0.7500  0.4165 0.7600
 0.4271 0.7700  0.4382 0.7800  0.4493 0.7900  0.4611 0.8000  0.4738 0.8100  0.4873 0.8200  0.5011 0.8300
 0.5150 0.8400  0.5311 0.8500  0.5471 0.8600  0.5649 0.8700  0.5825 0.8800  0.6028 0.8900  0.6234 0.9000
 0.6463 0.9100  0.6712 0.9200  0.6990 0.9300  0.72 0.9365 /
\endpicture
\hfill
\beginpicture
\setcoordinatesystem units <120mm, 60mm>
\setplotarea x from 0 to 0.6, y from 0 to 1
\axis left ticks numbered from 0 to 1 by 0.5 /
\axis bottom ticks numbered from  0 to 0.5 by 0.5  /
\arrow < 2mm>  [ .2679492, .7279404] from 0 0 to 0.6 0
\arrow < 2mm>  [ .2679492, .7279404] from 0 0 to 0 1.1
\put{$v$} [lt] <0pt, -5pt> at 0.6 0
\put{$g_w(v)$} [rB] <-5pt,0pt> at 0 1.1
\plot
 0.0000 0.0000  0.0017 0.0100  0.0035 0.0200  0.0052 0.0300  0.0071 0.0400  0.0089 0.0500  0.0107 0.0600
 0.0125 0.0700  0.0143 0.0800  0.0161 0.0900  0.0181 0.1000  0.0201 0.1100  0.0219 0.1200  0.0238 0.1300
 0.0258 0.1400  0.0278 0.1500  0.0298 0.1600  0.0317 0.1700  0.0339 0.1800  0.0360 0.1900  0.0379 0.2000
 0.0400 0.2100  0.0422 0.2200  0.0443 0.2300  0.0465 0.2400  0.0486 0.2500  0.0508 0.2600  0.0532 0.2700
 0.0554 0.2800  0.0578 0.2900  0.0601 0.3000  0.0625 0.3100  0.0648 0.3200  0.0672 0.3300  0.0696 0.3400
 0.0721 0.3500  0.0745 0.3600  0.0771 0.3700  0.0795 0.3800  0.0821 0.3900  0.0849 0.4000  0.0877 0.4100
 0.0904 0.4200  0.0930 0.4300  0.0959 0.4400  0.0989 0.4500  0.1016 0.4600  0.1046 0.4700  0.1075 0.4800
 0.1106 0.4900  0.1135 0.5000  0.1166 0.5100  0.1197 0.5200  0.1230 0.5300  0.1264 0.5400  0.1300 0.5500
 0.1335 0.5600  0.1370 0.5700  0.1407 0.5800  0.1445 0.5900  0.1481 0.6000  0.1521 0.6100  0.1559 0.6200
 0.1599 0.6300  0.1638 0.6400  0.1678 0.6500  0.1721 0.6600  0.1766 0.6700  0.1811 0.6800  0.1856 0.6900
 0.1901 0.7000  0.1950 0.7100  0.2001 0.7200  0.2050 0.7300  0.2100 0.7400  0.2155 0.7500  0.2212 0.7600
 0.2270 0.7700  0.2328 0.7800  0.2390 0.7900  0.2456 0.8000  0.2525 0.8100  0.2596 0.8200  0.2668 0.8300
 0.2746 0.8400  0.2835 0.8500  0.2925 0.8600  0.3020 0.8700  0.3123 0.8800  0.3238 0.8900  0.3357 0.9000
 0.3493 0.9100  0.3636 0.9200  0.3801 0.9300  0.3974 0.9400  0.4182 0.9500  0.4437 0.9600  0.4766 0.9700
 0.5190 0.9800  0.5877 0.9900 /
\setdashes
\plot
 0.0000 0.0000  0.0012 0.0100  0.0024 0.0200  0.0035 0.0300  0.0048 0.0400  0.0060 0.0500  0.0073 0.0600
 0.0085 0.0700  0.0098 0.0800  0.0110 0.0900  0.0123 0.1000  0.0136 0.1100  0.0149 0.1200  0.0162 0.1300
 0.0175 0.1400  0.0189 0.1500  0.0203 0.1600  0.0216 0.1700  0.0230 0.1800  0.0245 0.1900  0.0258 0.2000
 0.0272 0.2100  0.0287 0.2200  0.0301 0.2300  0.0316 0.2400  0.0330 0.2500  0.0346 0.2600  0.0362 0.2700
 0.0377 0.2800  0.0393 0.2900  0.0409 0.3000  0.0425 0.3100  0.0441 0.3200  0.0457 0.3300  0.0473 0.3400
 0.0490 0.3500  0.0507 0.3600  0.0525 0.3700  0.0541 0.3800  0.0559 0.3900  0.0577 0.4000  0.0597 0.4100
 0.0615 0.4200  0.0633 0.4300  0.0653 0.4400  0.0673 0.4500  0.0691 0.4600  0.0711 0.4700  0.0731 0.4800
 0.0753 0.4900  0.0772 0.5000  0.0794 0.5100  0.0814 0.5200  0.0837 0.5300  0.0860 0.5400  0.0885 0.5500
 0.0908 0.5600  0.0932 0.5700  0.0957 0.5800  0.0983 0.5900  0.1008 0.6000  0.1035 0.6100  0.1061 0.6200
 0.1088 0.6300  0.1115 0.6400  0.1142 0.6500  0.1171 0.6600  0.1201 0.6700  0.1232 0.6800  0.1263 0.6900
 0.1293 0.7000  0.1327 0.7100  0.1361 0.7200  0.1395 0.7300  0.1429 0.7400  0.1466 0.7500  0.1505 0.7600
 0.1544 0.7700  0.1584 0.7800  0.1626 0.7900  0.1671 0.8000  0.1718 0.8100  0.1766 0.8200  0.1815 0.8300
 0.1868 0.8400  0.1929 0.8500  0.1990 0.8600  0.2055 0.8700  0.2125 0.8800  0.2203 0.8900  0.2284 0.9000
 0.2376 0.9100  0.2474 0.9200  0.2586 0.9300  0.2704 0.9400  0.2845 0.9500  0.3018 0.9600  0.3243 0.9700
 0.3531 0.9800  0.3998 0.9900  0.4066 0.9910  0.4140 0.9920  0.4232 0.9930  0.4340 0.9940  0.4449 0.9950
 0.4581 0.9960  0.4792 0.9970  0.5071 0.9980  0.5477 0.9990 0.6 1 /
\setdashpattern <6pt,3pt,1pt,3pt>
\plot
 0.0000 0.0000  0.0027 0.0100  0.0055 0.0200  0.0081 0.0300  0.0110 0.0400  0.0139 0.0500  0.0167 0.0600
 0.0195 0.0700  0.0224 0.0800  0.0252 0.0900  0.0282 0.1000  0.0313 0.1100  0.0342 0.1200  0.0372 0.1300
 0.0402 0.1400  0.0433 0.1500  0.0465 0.1600  0.0495 0.1700  0.0528 0.1800  0.0561 0.1900  0.0592 0.2000
 0.0625 0.2100  0.0659 0.2200  0.0691 0.2300  0.0725 0.2400  0.0758 0.2500  0.0793 0.2600  0.0831 0.2700
 0.0865 0.2800  0.0902 0.2900  0.0939 0.3000  0.0976 0.3100  0.1011 0.3200  0.1048 0.3300  0.1086 0.3400
 0.1125 0.3500  0.1163 0.3600  0.1204 0.3700  0.1241 0.3800  0.1282 0.3900  0.1324 0.4000  0.1369 0.4100
 0.1410 0.4200  0.1452 0.4300  0.1497 0.4400  0.1543 0.4500  0.1586 0.4600  0.1632 0.4700  0.1678 0.4800
 0.1727 0.4900  0.1772 0.5000  0.1821 0.5100  0.1868 0.5200  0.1920 0.5300  0.1973 0.5400  0.2030 0.5500
 0.2083 0.5600  0.2138 0.5700  0.2196 0.5800  0.2255 0.5900  0.2312 0.6000  0.2374 0.6100  0.2433 0.6200
 0.2496 0.6300  0.2557 0.6400  0.2620 0.6500  0.2687 0.6600  0.2756 0.6700  0.2826 0.6800  0.2897 0.6900
 0.2967 0.7000  0.3044 0.7100  0.3122 0.7200  0.3200 0.7300  0.3278 0.7400  0.3364 0.7500  0.3453 0.7600
 0.3543 0.7700  0.3633 0.7800  0.3731 0.7900  0.3833 0.8000  0.3942 0.8100  0.4052 0.8200  0.4164 0.8300
 0.4285 0.8400  0.4425 0.8500  0.4565 0.8600  0.4713 0.8700  0.4875 0.8800  0.5054 0.8900  0.5240 0.9000
 0.5452 0.9100  0.5675 0.9200  0.5932 0.9300  /
\endpicture
}
\vskip12pt
\centerline{\hfill\smash{$m=10$}\hfill\hfill\smash{$m=20$}\hfill}
\vskip10pt
\centerline{
\hfill
\raise3pt\hbox{\beginpicture
\setcoordinatesystem units <10mm, 2mm>
\setplotarea x from 0 to 1, y from 0 to 1
\setsolid
\plot 0 0 1 0 /
      \endpicture}
Borda
\hfill
\raise3pt\hbox{\beginpicture
\setcoordinatesystem units <10mm, 2mm>
\setplotarea x from 0 to 1, y from 0 to 1
\setdashes
\plot 0 0 1 0 /
      \endpicture}
plurality
\hfill
\raise3pt\hbox{\beginpicture
\setcoordinatesystem units <10mm, 2mm>
\setplotarea x from 0 to 1, y from 0 to 1
\setdashpattern <6pt,3pt,1pt,3pt>
\plot 0 0 1 0 /
      \endpicture}
$m/2$-approval
\hfill
           }
\vskip10pt
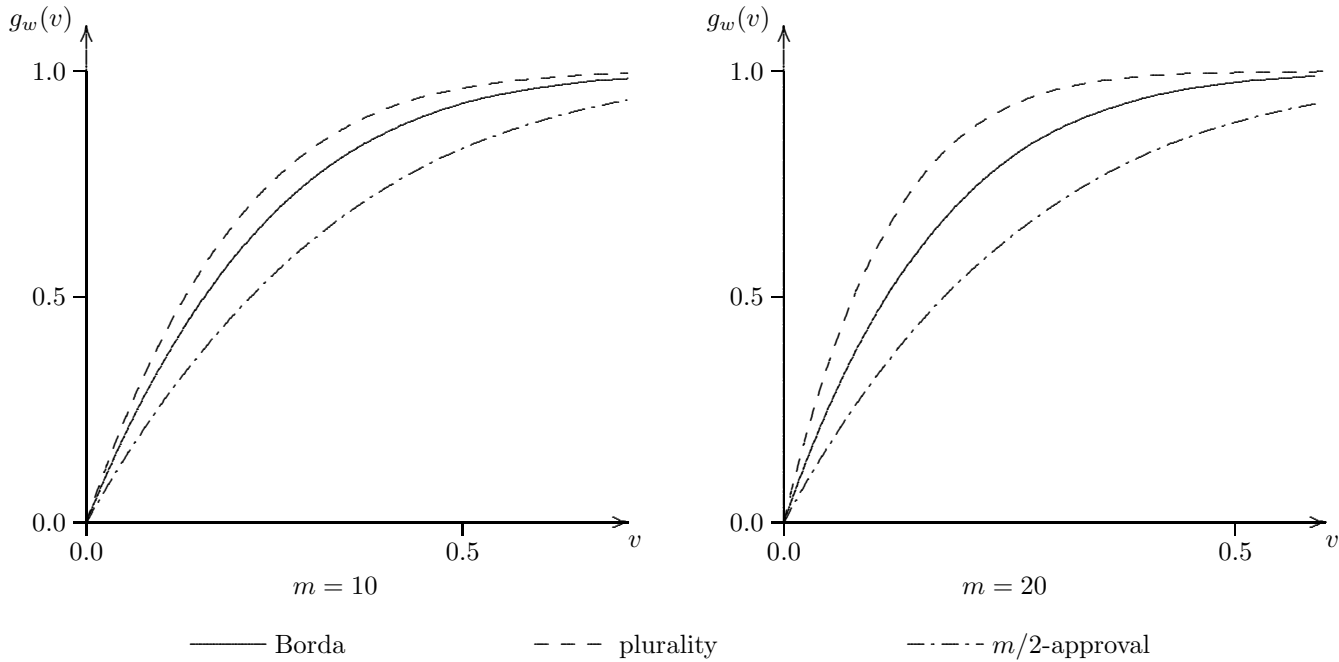
\centerline{Figure \gwgtentwenty: the functions $g_w(v)=\Prob{V_w\leq v}$ for some ten- and twenty-candidate voting rules.}
}
\vfill
\endinsert

\subchunk{}
Figures \gwgthreefour--\gwgtentwenty\ show
the manipulability of the rules for particular numbers of candidates.
We see from Figures \gwgthreefour\ and \gwgfivesix\ that for $m=4,5,6$ there is not much difference
in the manipulability of the common rules, at least by comparison with the differences evident
when $m=3$.
In the graph for $m=5$ we can see the Borda rule being dominated (by 2-approval) for the first time.
Also worth noting is that for $m\geq5$ there is very little difference between plurality and
anti-plurality from a manipulation point of view. Anti-plurality has a slight additional chance
of resisting attack by a large coalition -- and on this basis dominates plurality -- but this
advantage has become almost imperceptible by $m=6$.

Figure \gwgtentwenty\ shows the behaviour of the rules when there are many candidates.
In this figure, the curves for the anti-plurality rule have been left out, as they are
indistinguishable from those for plurality.
We see the $\grint{m/2}$-approval rule dominating the others.
It should be noted, though, that the IC hypothesis is possibly unconvincing when applied
to elections with many candidates, as it assumes in effect that all candidates are about
equally popular.

\section{6. Conclusions.}

The technique presented in this paper makes it possible to compute the (IC) limiting probability
distribution of the minimum manipulating coalition size for any positional rule, with any
number of candidates.

The consideration of coalition sizes is especially useful when comparing the rules.
Some rules are especially resistant to manipulation by small coalitions,
while others fare better with respect to manipulation by large coalitions.
Previous work has made this distinction in a rather limited way, by considering
``individual" and ``coalitional" manipulation (the latter meaning that the coalition
may be of any size). But these extremes may be somewhat uninformative.
Given a large voter population, all positional rules (except anti-plurality) are highly
likely to be manipulable by some coalition, and highly unlikely to be manipulable by
any individual. Studying coalitions of intermediate sizes starts to reveal more
differences between the rules.

The picture also changes when the number of candidates is varied. Much previous work
has concentrated on the three-candidate case, for which the behaviour of the rules is
quite different (Borda is least susceptible to small-coalition manipulation, anti-plurality
to large-coalition manipulation). The four-candidate case is similar, except that
the differences between rules are smaller. But when there are five candidates, it appears that
all positional rules are about equally manipulable, across the whole range of coalition sizes,
and there is not much to choose between them. With six or more candidates, the
$\grint{m/2}$-approval rules emerge as favourites.

Another, perhaps surprising, conclusion is that for $m\geq5$ there is very little difference
between plurality and anti-plurality from a manipulation point of view.
The approximate symmetry between these rules does not appear when $m=3$, and so
does not appear to have been noticed before (although it was recognized in a more limited
sense already in \cite{\Saari}).

It would be possible in principle to produce results like those in this paper for the IAC
voter behaviour model. However, the technique of analysis would have to be quite different.
Rather than reducing the probabilities to those involving normal distributions, the
calculations would entail the computation of convex volumes, as outlined in
\cite{\WilsonPritchardVolumes}.

It would also be of interest to produce limiting distributions similar to those in this
paper (or at least, graphs like those in Figures \gwgthreefour--\gwgtentwenty)
for voting rules other than positional rules.
Approval voting (\cite{\BramsFishburn}) would be an especially attractive target.
However, this too would require new techniques.

\startbibliography
\bibentry{\Bazaraa}{M. Bazaraa, J. Jarvis, and H. Sherali}{Linear programming and network flows}{Wiley, 2005.}
\bibentry{\BramsFishburn}{S. Brams and P. Fishburn}{Approval voting}{Springer, 2007.}
\bibentry{\BergLepelley}{S. Berg and D. Lepelley}{On probability models in voting theory}
 {Statistica Neerlandica 48 (1994), 133--146.}
\bibentry{\Chamberlin}{J. Chamberlin}{Investigation into the relative manipulability of
 four voting systems}{Behavior Science 30 (1985), 195--203.}
\bibentry{\Durrett}{R. Durrett}{Probability: Theory and examples}{Duxbury, 1996.}
\bibentry{\Gibbard}{A. Gibbard}{Manipulation of voting schemes: a general result}
 {Econometrica 41 (1973), 587--601.}
\bibentry{\KimRoush}{K. H. Kim and F. W. Roush}{Statistical manipulability of social choice functions}
 {Group Decision and Negotiation 5 (1996), 263--282.}
\bibentry{\LepelleyMbih}{D. Lepelley and B. Mbih}
 {The proportion of coalitionally unstable situations under the plurality rule}
 {Economic Letters 24 (1987), 311--315.}
\bibentry{\Nitzan}{S. Nitzan}{The vulnerability of point-voting schemes to
 preference variation and strategic manipulation}{Public Choice 47 (1985), 349--370.}
\bibentry{\PapadimitriouSteiglitz}{C. Papadimitriou and K. Steiglitz}{Combinatorial optimization}{Prentice-Hall, 1982.}
\bibentry{\ParkerRardin}{G. Parker and R. Rardin}{Discrete Optimization}{Academic Press, 1988.}
\bibentry{\PritchardSlinko}{G. Pritchard and A. Slinko}{On the average minimum size of a manipulating coalition}
 {Social Choice and Welfare 27 (2006), 263--277.}
\bibentry{\PritchardWilsonExactResults}{G. Pritchard and M. Wilson}
 {Exact results on manipulability of positional voting rules}{Social Choice and Welfare, in press.}
\bibentry{\Saari}{D. Saari}{Susceptibility to manipulation}{Public Choice 64 (1990), 21--41.}
\bibentry{\Satterthwaite}{M. A. Satterthwaite}{Strategy-proofness and Arrow's conditions: 
 existence and correspondence theorems for voting procedures and social welfare functions}
 {J. Economic Theory 10 (1975), 187--217.}
\bibentry{\WilsonPritchardVolumes}{M. Wilson and G. Pritchard}
 {Probability calculations under the IAC hypothesis}{Mathematical Social Sciences, in press.}

\end